\documentclass[letterpaper,fleqn]{cas-sc}

\usepackage[authoryear,longnamesfirst]{natbib}
\setcitestyle{numbers,square}
\pdfoutput=1
\usepackage{amsmath}
\numberwithin{equation}{section}
\usepackage{float}
\usepackage{graphicx}
\usepackage{clrscode}
\usepackage{subfigure}
\usepackage{amsfonts}
\usepackage{mathrsfs}
\usepackage{latexsym}
\usepackage{amssymb}
\usepackage{color}
\usepackage{caption}
\usepackage[title]{appendix}
\usepackage[colorinlistoftodos,disable]{todonotes}
\renewcommand{\Im}[0]{\mathfrak{Im}}
\renewcommand{\Re}[0]{\mathfrak{Re}}

\usepackage{array}
\usepackage{algorithm}
\usepackage{algpseudocode}
\usepackage{amsmath}
\usepackage{graphics}
\usepackage{epsfig}
\begin{document}
%bkdef.tex (file name)
\newcommand{\bul}{\bu^\lambda}
\newcommand{\bum}{\bu^{\mu}}
\newcommand{\buleta}{\bul(\bet)}
\newcommand{\bldn}{\bld{n}}
\newcommand{\bxzero}{\bld{x}_0}
\newcommand{\bxone}{\bld{x}_1}
\newcommand{\bldv}{\bld{v}}
\newcommand{\bldgam}{\mbox{\boldmath$\gamma$}}
\newcommand{\bgam}{\mbox{\boldmath$\gamma$}}
\newcommand{\bal}{\mbox{\boldmath$\alpha$}}
\newcommand{\bet}{\bld{\eta}}
\newcommand{\bldeta}{\bld{\eta}}
\newcommand{\bnu}{\bld{\nu}}
\newcommand{\bom}{\bld{\omega}}
\newcommand{\A}{ {\alpha} }
\newcommand{\uA}{ {\underline{A}} }
\newcommand{\Ab}{ {\overline{A}} }
\newcommand{\ab}{\bar{\alpha}}
\newcommand{\ah}{\hat a(0)}
\newcommand{\ahat}{\hat a}
\newcommand{\Ah}{ \hat{\alpha} }
\newcommand{\ahp}{\hat a\prime(0)}
\newcommand{\ahpp}{\hat a{\prime\prime}(0)}
\newcommand{\Ai}{{\rm Ai}}
\newcommand{\AjN}{A_j^{(N)}}
\newcommand{\AkN}{A_k^{(N)}}
\newcommand{\AkphN}{A_{k-\hf}^{(N)}}
\newcommand{\al}{\alpha}
\newcommand{\ang}[1]{\langle #1 \rangle }
\newcommand{\AN}[1]{\stackrel{N}{A\,}_{#1}\!}
\newcommand{\ANt}[1]{\stackrel{{\tilde N}}{A\,}_{#1}\!}
\newcommand{\atz}{{\tilde a}_0}
\newcommand{\B}{ {\beta} }
\newcommand{\ba}{\begin{array}}
\newcommand{\bA}{{A}}
\newcommand{\bAN}{\stackrel{N}{\bld{A}\,}\!\!}
\newcommand{\bANt}{\stackrel{{\tilde N}}{\bld{A}\,}\!\!}
\newcommand{\bApN}{\stackrel{N}{\bld{A}}\!\!{}^T\!}
\newcommand{\bb}{\bar{\beta}}
\newcommand{\bBM}{\stackrel{M}{\bld{B}\,}\!\!\!}
\newcommand{\bc}{{\bld c}}
\newcommand{\bchi}{{ \mbox{\boldmath $\chi$} }}
\newcommand{\bd}{\bld{d}}
\newcommand{\bD}{\bld{D}}
\newcommand{\be}{\beta}
\newcommand{\bE}{\bld{E}}
\newcommand{\beq}[1]{ 
%      \marginpar{[#1]}
        \begin{equation} \label{#1} }
\newcommand{\beql}[1]{
\begin{equation} \label{#1}}

\newcommand{\beqa}{\begin{eqnarray}}
\newcommand{\blde}{{\bld e}}
\newcommand{\beps}{{\bld \epsilon}}
\newcommand{\bldf}{{\bld f}}
\newcommand{\BF}[1]{ \par \indent \begin{fgr} \refstepcounter{fgr}
        \label{#1} {\bf \thefgr.} }
\newcommand{\bh}{{\hat \beta}}
\newcommand{\bld}[1]{{ \mbox{\boldmath $#1$} }}
\newcommand{\bnab}{\mbox{\boldmath$\nabla$}}
\newcommand{\bcdot}{{\bld{\cdot}}}
\newcommand{\barn}{{\bar n}}
\newcommand{\barbn}{{\bar \bn}}
\newcommand{\barbx}{{\bar \bx}}
\newcommand{\bhtal}{\hat{\bld{\alpha}}}
\newcommand{\bldscr}[1]{{ \mbox{\boldmath ${\scriptstyle #1}$} }}
\newcommand{\bF}{\bld{F}}
\newcommand{\bG}{\bld{G}}
\newcommand{\bH}{\bld{H}}
\newcommand{\bI}{\bld{I}}
\newcommand{\bk}{\bld{k}}
\newcommand{\bl}{\bld{l}}
\newcommand{\BM}[1]{\stackrel{M}{B\,}_{#1}\!}
\newcommand{\BMt}[1]{\stackrel{{\tilde M}}{B\,}_{#1}\!}
\newcommand{\bn}{{\bld n}}
\newcommand{\bp}{\bld{ p}}
\newcommand{\bP}{{\bld P}}
\newcommand{\bptl}{\bld{\partial}}
\newcommand{\bq}{\bld{q}}
\newcommand{\bldr}{\bld{r}}
\newcommand{\Brm}{{\rm B}}
\newcommand{\bs}{\bld{s}}
\newcommand{\bS}{\bld{S}}
\newcommand{\bsigma}{\bld{\sigma}}
\newcommand{\bsig}{\bld{\sigma}}
\newcommand{\bt}{\bld{t}}
\newcommand{\bT}{\bld{T}}
\newcommand{\btau}{\bld{\tau}}
\newcommand{\btdal}{\tilde{\bld{\alpha}}}
\newcommand{\bu}{\bld{u}}
\newcommand{\bU}{\bld{U}}
\newcommand{\bup}{\bld{\upsilon}}
\newcommand{\bv}{\bld{v}}
\newcommand{\bV}{\bld{V}}
\newcommand{\bw}{{\bld w}}
\newcommand{\bW}{{\bld W}}
\newcommand{\bx}{\bld{x}}
\newcommand{\bxi}{\bld{\xi}}
\newcommand{\bxidotx}{\bxi\cdotb\bx}
\newcommand{\bX}{\bld{X}}
\newcommand{\by}{\bld{y}}
\newcommand{\bz}{\bld{z}}
\newcommand{\bY}{\bld{Y}}
\newcommand{\bzeta}{\bld{\zeta}}
\newcommand{\bzero}{{\bf 0}}
\newcommand{\cdotb}{{\bld{\cdot}}}
\newcommand{\Chapter}[1]{\setcounter{equation}{0} \chapter{#1}}
\newcommand{\cth}{{\cos\theta}}
\newcommand{\Db}{ \bar{\delta} }
\newcommand{\dd}{\,{\rm d}}
\newcommand{\deldel}[2]{ \fr{\partial #1}{\partial #2} }
\newcommand{\del}[1]{ {\partial_#1} }
\newcommand{\de}{\delta}
\newcommand{\dee}[1]{ \,{\rm d}#1 }
\newcommand{\drm}[1]{ \,{\rm d}#1 }
\newcommand{\dlt}{\delta}
\newcommand{\Dh}{ \hat{\delta} }
\newcommand{\Dlt}{\Delta}
\newcommand{\dol}[1]{\overline{\overline{#1}}}
\newcommand{\e}[1]{\bld{e}_{#1}}
\newcommand{\edot}[1]{\dot{\bld{e}}_{#1}}
\newcommand{\ue}{ {\underline{e}} }
\newcommand{\ea}{\end{array}}
\newcommand{\eeq}{ \end{equation} }
\newcommand{\eeqa}{ \end{eqnarray} }
\newcommand{\EF}{ \end{fgr} }
\newcommand{\eps}{ {\epsilon} }
\newcommand{\eq}[2]{
        \vspace{24pt}
        \marginpar
        [\hfill{\footnotesize #1}]
        {\footnotesize #1\hfill} \vspace{-24pt}
        \begin{equation} \label{#1} 
        #2
        \eeq}
\newcommand{\erm}{{\rm e}}
\newcommand{\erfc}{{\rm erfc}}
\newcommand{\eiop}{{\rm e}^{\im\om\psi}}
\newcommand{\eiot}{{\rm e}^{{\rm i}\omega t}}
\newcommand{\F}[1]{ {\cal F}\{#1\} }
\newcommand{\fb}{\bld{f}}
\newcommand{\uf}{\ul{f}}
\newcommand{\fr}[2]{{ \displaystyle \frac{#1}{#2} }}
\newcommand{\frr}[2]{{\textstyle \frac{#1}{#2} }}
\newcommand{\frscr}[2]{{\scriptstyle \frac{#1}{#2} }}
\newcommand{\frtxt}[2]{{\textstyle \frac{#1}{#2} }}
\newcommand{\g}[1]{ {\gamma_{#1}} }
\newcommand{\ga}{ {\gamma_{\bar{\alpha}}} }
\newcommand{\gam}{\gamma}
\newcommand{\gams}{\gamma^2}
\newcommand{\gb}{ {\gamma_{\bar{\beta}}} }
\newcommand{\Gb}{ \bar{\gamma} }
\newcommand{\gef}{\hat{\gamma}^{\rm eff}}
\newcommand{\geff}{  { \gamma_{{\rm eff}} }  }
\newcommand{\Gh}{ \hat{\gamma} }
\newcommand{\hbxi}{ \hat{\bxi} }
\newcommand{\hbxidotx}{\hbxi\cdotb\bx}
\newcommand{\hem}{\mbox{\hspace{1em}}}
\newcommand{\hf}{{ \textstyle \frac{1}{2} }}
\newcommand{\Hf}{\frac{1}{2}}
\newcommand{\hff}{{ \scriptstyle \frac{1}{2} }}
\newcommand{\hftxt}{{ \textstyle \frac{1}{2} }}
\newcommand{\hfdis}{{ \displaystyle \frac{1}{2} }}
\newcommand{\hfscr}{{ \scriptstyle \frac{1}{2} }}
\newcommand{\rh}{{\rm h}}
\newcommand{\Hrm}{{\rm H}}
\newcommand{\hw}{\hat w}
\newcommand{\hxi}{\hat \xi}
\newcommand{\im}{{\rm i}}
\newcommand{\Imag}[1]{ \Im\{{#1}\} }
\newcommand{\imm}{{\cal im}}
\newcommand{\INT}{{\displaystyle \int}}
\newcommand{\itf}{$\mbox{in the form}$}
\newcommand{\io}{{{\rm i}\omega}}
\newcommand{\dint}{ \displaystyle{\int} }
\newcommand{\Int}[2]{{\displaystyle \int_{#1}^{#2}}}
\newcommand{\LA}{ {\lambda} }
\newcommand{\la}{ {\langle} }
\newcommand{\Lb}{ \bar{\lambda} }
\newcommand{\Lh}{ \hat{\lambda} }
\newcommand{\Lim}{{\displaystyle \lim_{N\to\infty}}}
\newcommand{\lm}{\lim_{N\to\infty}}
\newcommand{\M}{ {\mu} }
\newcommand{\Mb}{ \bar{\mu} }
\newcommand{\Mh}{ \hat{\mu} }
\newcommand{\n}{\bld{n}}
\newcommand{\N}{\nu}
\newcommand{\nab}{{\nabla}}
\newcommand{\nabs}{{\nabla^2}}
\newcommand{\NM}{^{(\tilde{N}\hat{M})}}
\newcommand{\npvf}{ \vspace*{\fill}\newpage}
\newcommand{\Nth}{ {\frac{1}{N}} }
\newcommand{\om}{\omega}
\newcommand{\omb}{\bar{\omega}}
\newcommand{\oms}{\omega^2}
\newcommand{\pq}{ (p^2+q^2) }
\newcommand{\Prm}{{\rm P}}
\newcommand{\Pro}[1]{{\displaystyle \prod_{#1}}}
\newcommand{\Prod}[2]{{\displaystyle \prod_{#1}^{#2}}}
\newcommand{\PS}{ {\psi} }
\newcommand{\psib}{ \bar{\psi} }
\newcommand{\qt}{ {\frac{1}{4}} }
\newcommand{\qttxt}{ \textstyle{\frac{1}{4}} }
\newcommand{\Qt}{ \frac{1}{4} }
\newcommand{\re}{{\rm re}}
\newcommand{\ra}{{\rangle}}
\newcommand{\rH}{{\rm H}}
\newcommand{\rJ}{{\rm J}}
\newcommand{\Real}[1]{\Re{\{ {#1} \}} }
\newcommand{\real}[1]{\Re{\{ {#1} \}} }
\newcommand{\rf}[1]{(\ref{#1})}
\newcommand{\refrf}[1]{[\ref{#1}]}
\newcommand{\rF}[1]{\ref{#1}}
\renewcommand{\r}{\rho}
\newcommand{\rhob}{\bar{\rho}}
\newcommand{\rhoh}{\hat{\rho}}
\newcommand{\rkN}{r_k^{(N)}}
\newcommand{\rkhN}{r_{k-\hf}^{(N)}}
\newcommand{\rkphN}{r_{k+\hf}^{(N)}}
\newcommand{\rmbox}[1]{{\mbox{\rm #1}}}
\newcommand{\rotn}{{\mbox{\boldmath $\Omega$}}}
\newcommand{\s}{\sigma}
\newcommand{\us}{ \underline{s}}
\newcommand{\Section}[1]{\setcounter{equation}{0} \section{#1}}
\newcommand{\sgn}{{\rm sgn}}
\newcommand{\Sm}[2]{{\displaystyle \sum_{#1}^{#2}}}
\newcommand{\Srm}{{\rm S}}
\newcommand{\ssf}[1]{{\sf #1}}
\newcommand{\sth}{{\sin\theta}}
\newcommand{\Strain}{{\mbox{\boldmath $\varepsilon$}}}
\newcommand{\strain}{{\varepsilon}}
\newcommand{\Stress}{{\mbox{\boldmath $\tau$}}}
\newcommand{\stress}{\tau}
\newcommand{\PKStress}{\bld{T}}
\newcommand{\PKstress}{T}
\newcommand{\ptl}[1]{ {\partial_#1} }
\newcommand{\Sum}[1]{{\displaystyle \sum_{#1}}}
\newcommand{\T}{\bld{T}}
\newcommand{\taub}{\bar{\tau}}
\newcommand{\taupm}{\tau\prime}
\newcommand{\tb}{\bld{t}}
\renewcommand{\theequation}{\arabic{section}.\arabic{equation}}
\newtheorem{theorem}{Theorem}
\newcommand{\ti}{\theta_i}
\newcommand{\tu}{{\tilde u}}
\newcommand{\tmu}{{\tilde \mu}}
\newcommand{\tnu}{{\tilde \nu}}
\newcommand{\Twmw}{$\mbox{Thus we may write}$}
\newcommand{\txN}{{\tilde x}^{(N)}}
\newcommand{\ul}[1]{\und{#1}}
\newcommand{\und}[1]{\underline{#1}}
\newcommand{\uu}{\ul{u}}
\newcommand{\uxi}{\bld{\xi}}
\newcommand{\ux}{ \bld{x}}
\newcommand{\vb}{\bar{v}}
\newcommand{\vf}{\vfill}
\newcommand{\vsf}{\vspace*{\fill}}
\newcommand{\vh}{\hat{v}}
\newcommand{\Vh}{\hat{V}}
\newcommand{\wmw}{$\mbox{we may write}$}
\newcommand{\wrt}{$\mbox{with respect to}$}
\newcommand{\xb}{\bld{x}}
\newcommand{\x}{\xi}
\newcommand{\X}{\bld{X}}
\newcommand{\xki}{\xi_{k,i}}
\newcommand{\xN}[1]{x_{#1}^{(N)}}
\newcommand{\xNxi}[1]{{x_{#1}^{(N)} (\X)}}
\newcommand{\y}{\bld{y}}
\newcommand{\Y}{\bld{Y}}
\newcommand{\yb}{\bld{y}}
\newcommand{\Z}{ {\zeta} }
\newcommand{\Zb}{ {\overline{\zeta}} }
\newcommand{\ZkmN}{ {\zeta _{k-1}^{(N)}} }
\newcommand{\ZkN}{ {\zeta _k^{(N)}} }
\newcommand{\z}[1]{z_{#1}}

% Sans Serif
\newcommand{\steq}[1]{\stackrel{\scriptscriptstyle{\mathrm{#1}}}{=}}
\newcommand{\stlp}[1]{\stackrel{\scriptscriptstyle{\mathrm{#1}}}{(}}
\newcommand{\strp}[1]{\stackrel{\scriptscriptstyle{\mathrm{#1}}}{)}}
\newcommand{\stk}[2]{\hspace{1pt}\stackrel{\scriptscriptstyle{\mathrm{#1}}}{#2}\hspace{1pt}}
\newcommand{\stkrse}[2]{\!\raisebox{.3ex}{\mbox{$\stk{#1}{#2}$}}\!}
\newcommand{\szero}{\bld{\sf 0}}
\newcommand{\sA}{\bld{\sf A}}
\newcommand{\sB}{\bld{\sf B}}
\newcommand{\sC}{\bld{\sf C}}
\newcommand{\sE}{\bld{\sf E}}
\newcommand{\sG}{\bld{\sf G}}
\newcommand{\sI}{\bld{\sf I}}
\newcommand{\sM}{\bld{\sf M}}
\newcommand{\sR}{\bld{\sf R}}
\newcommand{\sS}{\bld{\sf S}}
\newcommand{\sU}{\bld{\sf U}}

\newcommand{\sub}[1]{_{,#1}}

%Calligraphic
\newcommand{\cA}{{\cal A}}
\newcommand{\calC}{{\cal C}}
\newcommand{\calE}{{\cal E}}
\newcommand{\calF}{{\cal F}}
\newcommand{\calP}{{\cal P}}
\newcommand{\calR}{{\cal R}}
\newcommand{\calS}{{\cal S}}
\newcommand{\calT}{{\cal T}}

%Figures

\newcommand{\FFig}[4]{
\newpage\centerline{\bf #3}
\begin{figure}[#1]
\centerline{\psfig{height=#2in,file=BurridgeFigure#3.pdf}}
\caption{ \label{BurridgeFigure#3} {\footnotesize #4} \newline  BurridgeFigure#3} 
\end{figure}
}

\newcommand{\cross}{\hspace{-2pt}\bld{\times}\hspace{-1pt}}                              %
\newcommand{\irm}{{\rm i}}                                                                                    %
\newcommand{\aps}{  {
\rotatebox{-90}{ \hspace{-10pt}$\backslash$\hspace{-8pt}$\supset$ }\hspace{-1pt}
}  }

\let\WriteBookmarks\relax
\def\floatpagepagefraction{0.5}
\def\textpagefraction{.001}
\let\printorcid\relax % 可去掉页面下方的ORCID(s)

% \shorttitle{<short title of the paper for running head>}
\shorttitle{Lagrangian Hadamard Integrator}

% \shortauthors{<short author list for running head>}
\shortauthors{Wei, Cheng, Burridge, and Qian}

\title[mode = title]{Hadamard integrator for time-dependent wave equations: Lagrangian formulation via ray tracing}

% \tnotemark[<tnote number>]
\tnotemark[1]

\author[1]{Yuxiao Wei}[type=editor]
\ead[URL]{19110180029@fudan.edu.cn}

\author[1]{Jin Cheng}[]
\ead[URL]{jcheng@fudan.edu.cn}
\author[2]{Robert Burridge}
\ead[URL]{burridge137@gmail.com}
\author[3]{Jianliang Qian}[]
\cormark[1]
\ead[URL]{jqian@msu.edu}
%\ead{rishi@sayahna.org}
%\ead[URL]{www.sayahna.org}

%\ead{karl@freefriends.org}
%\ead[URL]{www.tug.org}

\address[1]{School of Mathematical Sciences, Fudan University, Shanghai 200433, China }
\address[2]{Department of Mathematics and Statistics, University of New Mexico, Albuquerque, NM 87131, USA}
\address[3]{Department of Mathematics, Michigan State University, East Lansing, MI 48824, USA}

\cortext[1]{Corresponding author}

\begin{abstract}
Starting from the time-domain Kirchhoff-Huygens representation of wave solutions, we propose a novel Hadamard integrator for the self-adjoint time-dependent wave equation in an inhomogeneous medium. First, we create a new asymptotic series based on the Gelfand-Shilov function, dubbed Hadamard's ansatz, to approximate the Green's function of the time-dependent wave equation. Accordingly, the governing equations and related  initializations for the eikonal and Hadamard coefficients are derived using the properties of the Gelfand-Shilov generalized function. Second, incorporating the leading term of Hadamard's ansatz into the Kirchhoff-Huygens representation, we develop an original Hadamard integrator for the Cauchy problem of the time-dependent wave equation and derive the corresponding Lagrangian formulation in geodesic polar coordinates. Third, to construct the Hadamard integrator in the Lagrangian formulation efficiently, we use a short-time ray tracing method to obtain equal-time wavefront locations accurately, and we further develop fast algorithms to compute Chebyshev-polynomial based low-rank representations of both wavefront locations and variants of Hadamard coefficients. Fourth, equipped with these low-rank representations, we apply the Hadamard integrator to efficiently solve time-dependent wave equations with highly oscillatory initial conditions, where the time step size is independent of the initial conditions. By judiciously choosing the medium-dependent time step, our new Hadamard integrator can propagate wave field beyond caustics implicitly and advance spatially overturning waves in time naturally. Moreover, since the integrator is independent of initial conditions, the Hadamard integrator can be applied to many different initial conditions once it is constructed. Both two-dimensional and three-dimensional numerical examples illustrate the accuracy and performance of the proposed method.
\end{abstract}


\begin{highlights}
\item Novel asymptotic series for the Green's function of time-dependent wave equation in an inhomogeneous medium.
\item Hadamard integrator that accurately propagates highly oscillatory wave field beyond caustics.
\item Low-rank representation based fast ray tracing method.
\end{highlights}

% Keywords

\maketitle

% Main text

\section{Introduction}
We consider the Cauchy problem for the self-adjoint wave equation in $m$-dimensional space $\mathbb{R}^m$,
\begin{equation}\label{1.1}
  \rho u_{tt}-\nab\cdot(\nu \nabla u)=0,\; \boldsymbol{x}\in \mathbb{R}^m,\; t>0
\end{equation}
with initial conditions
\begin{equation}\label{1.2}
  u(0,\boldsymbol{x})=u^1(\boldsymbol{x}),\; u_t(0,\boldsymbol{x})=u^2(\boldsymbol{x}),
\end{equation}
where $t$ is time, the subscripts ${}_t$ and ${}_{tt}$ represent the first and second time derivative, respectively, position $\boldsymbol{x}=[x_1,x_2,\cdots,x_m]^T$, the gradient operator $\nab =[\partial_{x_1},\partial_{x_2},\cdots,\partial_{x_m}]^T$, both variables $\rho$ and $\nu$ are analytic and positive functions of position $\boldsymbol{x}$, characterizing certain physical parameters of the medium, and  $u^1(\boldsymbol{x})$ and $u^2(\boldsymbol{x})$ are compactly supported, highly oscillatory $L^2$ functions. Our wave equation is time independent in the sense that it is invariant under shifts in time, but time dependent in the sense that it contains $t$ as an independent variable. When the initial conditions are highly oscillatory, the wave equation propagates these oscillations in space and time; however, direct numerical methods, such as finite-difference or finite-element methods, for such problems may suffer from dispersion or pollution errors \cite{baygoltur85,babsau00}, so that such methods require an enormous computational grid to resolve these oscillations and are thus very costly in practice. Therefore, alternative methods, such as geometrical-optics based asymptotic methods, have been sought to resolve these highly-oscillatory wave phenomena. We first propose a novel Hadamard asymptotic ansatz based on Gelfand and Shilov’s family of functions $f_{+}^\lambda$. On top of this new ansatz, we further develop an original Hadamard integrator to compute highly oscillatory time-dependent wave phenomena in inhomogeneous media.

To start with, we use the Gelfand-Shilov generalized function as the basis to seek an asymptotic representation of the Green's function $G(t,\boldsymbol{x}_0;\bx)$ which satisfies the wave equation (\ref{1.1}) with the initial conditions $u^1(\boldsymbol{x})=0$ and $u^2(\boldsymbol{x})=\frac{1}{\rho(\boldsymbol{x}_0)}\delta(\boldsymbol{x}-\boldsymbol{x}_0)$. Inserting this asymptotic representation into the wave equation and taking into account orders of singularities of the Gelfand-Shilov function and its derivatives, we can obtain time-independent eikonal and transport equations for the phase function and Hadamard coefficients, respectively, where the fact that the coefficients $\rho$ and $\nu$ are independent of $t$ enables us to get away with not having any time dependence in the eikonal and Hadamard coefficients. Since these Hadamard ingredients are independent of time, we can first precompute these functions and then compress them into low-rank representations which can be used for further time evolution. In particular, these low-rank representations allow us to rapidly construct a short-time caustic-free asymptotic Green's function which is valid locally in time. To solve initial value problems of the corresponding time-dependent wave equation globally in time, we incorporate the short-time caustic-free asymptotic Green's function into the time-domain Kirchhoff-Huygens representation formula so that we can take multiple local in-time steps to achieve global in-time caustic-friendly wave propagation, leading to the novel Hadamard integrator.

One of the essential difficulties in applying geometrical optics to construct Green's functions for wave equations is how to initialize the eikonal and amplitude functions at the source point \cite{avikel63,bab65,karkel59}. Here, inspired by our series works on Hadamard-Babich ansatzes for Helmholtz, Maxwell's, and elastic wave equations \cite{luqiabur18,luqiabur16,luqiabur16b,qiasonlubur21}, our newly proposed Hadamard's ansatz for time-dependent wave equations can be easily initialized as we will show.

Another essential difficulty in applying geometric optics is that it cannot handle caustics easily \cite{lax57,lud66,masfed81,babbul72,kellew95,ben03,engrun03,luoqiabur14a,luqiabur16}, and our Hadamard's ansatz is not an exception as it is also an asymptotic method. Although  caustics occur with high probability for wave propagation in inhomogeneous media \cite{whi84}, we are still able to use the geometrical-optics type method mainly because of the following fact \cite{avikel63,symqia03slowm}: in an isotropic medium the point-source eikonal equation has a locally smooth solution near the source point except the source point itself; this implies that caustics will not develop right away on the expanding wavefront away from the source. Therefore, in a local (spatial) neighborhood of the point source, the eikonal and amplitude functions from solving eikonal and transport equations are smooth except the point source; the resulting asymptotic Green's function is valid locally in that spatial neighborhood except the point source itself and thus is not uniform near the source point.

Then we immediately run into two questions. The first question is how to obtain {\it uniformly} accurate asymptotic Green's functions in that small spatial neighborhood {\it even} at the source point. The proposed Hadamard's ansatz comes to our rescue, where the crucial point is that although the eikonal itself is not differentiable at the source point, the squared eikonal is! This crucial point allows us to absorb the point-source singularity into the Gelfand-Shilov generalized function, so that we can initialize the Hadamard coefficients easily, resulting in a uniformly accurate asymptotic Green's function near and at the source point.

The second question is how to use {\it locally} valid asymptotic Green's functions to solve time-dependent wave equations {\it globally}. The answer is provided by incorporating the locally valid asymptotic Green's function into the time-domain Kirchhoff-Huygens representation formula of the time-dependent wave solution. To appreciate this subtle point, we need to characterize the caustic-free spatial neighborhood of the source point of the eikonal equation in terms of {\it time} and {\it space}. Since a caustic will need some time to develop away from the source point in an isotropic medium, we denote by $\bar{T}(\boldsymbol{x}_0)$ the time when the first caustic occurs for rays issuing from the source point $\boldsymbol{x}_0$, where time actually corresponds to the solution of the point-source eikonal equation.
Therefore, our short-time asymptotic Green's function excited at $(\boldsymbol{x}_0,0)$ is valid in the $(\boldsymbol{x},t)$ space-time domain
\begin{equation}\label{1.4}
\left\{(\boldsymbol{x},t):\; \tau(\boldsymbol{x},\boldsymbol{x}_0)<\bar{T}(\boldsymbol{x}_0), \; 0\leq t< \bar{T}(\boldsymbol{x}_0)\right\}.
\end{equation}
Since the eikonal $\tau(\boldsymbol{x},\boldsymbol{x}_0)$ and Hadamard coefficients are independent of time $t$, we just need to compute these quantities once and use them to construct short-time asymptotic Green's functions for all $(\boldsymbol{x},t)$ in the above space-time domain \eqref{1.4}.
Moreover, since all values of $\boldsymbol{x}_0$ are used by the propagator defined below, we set $\bar{T}$ to be the minimum of $\bar{T}(\boldsymbol{x}_0)$ as $\boldsymbol{x}_0$  varies over some relevant domain, where this domain should not be too large so that $\bar{T}$ might not be too small.
To march forward in time so as to solve the time-dependent wave equation globally in time, we incorporate the short-time asymptotic Green's function into the time-domain Kirchhoff-Huygens representation formula to define a short-time $\Delta t$ propagator, dubbed the Hadamard-Kirchhoff-Huygens (HKH) propagator,
where $\Delta t < \bar{T}$ only depends on the medium and is independent of the initial data.  Recursively applying this propagator in time yields the {\it Hadamard integrator} to solve time-dependent wave equations globally in time, where caustics are treated implicitly. Moreover, by marching forward in time, we are able to treat spatially overturning waves naturally.

The matter at hand now is how to implement the short-time HKH propagator efficiently. To tackle this challenging problem, we must surmount several obstacles.

The first obstacle is how to deal with integrals of Gelfand-Shilov functions, which have singularities near the $t$-wavefront,
\begin{equation}
     \{\boldsymbol{x}: \tau(\boldsymbol{x},\boldsymbol{x}_0)=t\}.
\end{equation}
In a caustic-free local neighborhood of source $\boldsymbol{x}_0$, we introduce the geodesic (ray) polar transformation
$P[\bx_0]: \bx \rightarrow (\tau,\boldsymbol{\omega})$
to simplify the generalized integrals, where $\tau$ is traveltime and $\boldsymbol{\omega}\in \mathbb{S}^{m-1}$ is take-off angle; by assumption, this transformation is well defined.
Since this transformation is from the Cartesian coordinates to the geodesic polar coordinates, we can use the Lagrangian ray-tracing method to trace wavefronts accurately, where locations of the $\tau$-wavefront with take-off angle $\boldsymbol{\omega}$ exactly yield the corresponding Cartesian coordinates $\bx$. Finally, a Gaussian quadrature in geodesic polar coordinates is used to compute the resulting Gelfand-Shilov integrals.

The second obstacle is how to obtain Hadamard coefficients efficiently. Fortunately, over the years we have developed high-order schemes for point-source eikonal and transport equations \cite{qiasym02adapt,zhazhaqia06,qiayualiuluobur16,luqiabur16}, and those schemes can be readily used to compute these coefficients.

The third obstacle is how to accelerate evaluation of Gelfand-Shilov  integrals in the HKH propagator. We first compress computed wavefront locations and Hadamard coefficients into low-rank representations by using multivariate Chebyshev polynomials. On top of low-rank representations, a block-wise matrix based partial summation allows us to evaluate Gelfand-Shilov integrals rapidly.

 \subsection{Related works}
 Hadamard's ansatz that we initiate here is inspired by the Hadamard method \cite{had23} which is outlined in Courant and Hilbert \cite{couhil62}, Chapter VI, Section 15.6. However, since our systematic derivation here is based on Gelfand-Shilov functions as well as their regularization techniques \cite{gelshi64}, it is original.

Fast Huygens sweeping (FHS) methods have been designed to solve Helmholtz equations \cite{luoqiabur14a,luqiabur16}, frequency-domain Maxwell's equations \cite{qialuyualuobur16,luqiabur18}, and frequency-domain elastic wave equations \cite{qiasonlubur23}, and these methods work by incorporating locally valid asymptotic Green's functions into the frequency-domain Kirchhoff-Huygens representations of corresponding wave solutions so that they can treat caustics implicitly in inhomogeneous media at high frequencies. However, since these methods have assumed the sub-horizontal condition \cite{symqia03slowm} for geodesics which is useful in many practical applications, the allowed wave propagation has a certain preferred spatial direction; consequently, the FHS methods are able to propagate wavefields through appropriately partitioned spatial layers by marching in that preferred spatial direction in a layer-by-layer fashion. However, such a spatial preference due to the sub-horizontal condition does come with a cost: the above FHS methods cannot handle overturning waves in that particular spatial direction since marching in a certain spatial direction is {\it unnatural}! Then a question arises immediately: which direction is {\it natural} for marching? It is the time direction. This is exactly what we are achieving in this article!

 \subsection{Plan of the paper}
 We introduce in Section 2 the Kirchhoff-Huygens representation formula which utilizes Green's functions to propagate waves.
We then propose in Section 3 a novel asymptotic series based on the Gelfand-Shilov function, dubbed Hadamard's ansatz, to approximate the Green's function of the time-dependent wave equation, where the governing equations and related initializations for the eikonal and Hadamard coefficients are derived using the properties of the Gelfand-Shilov generalized function. Incorporating the leading term of Hadamard's ansatz into the Kirchhoff-Huygens representation, we develop the Hadamard integrator for the Cauchy problem of the time-dependent wave equation and derive the corresponding Lagrangian formulation in geodesic polar coordinates in Section 4. We develop in Section 5 numerical strategies for implementing the Hadamard integrator. To accelerate evaluations of various Gelfand-Shilov integrals, in Section 6 we construct multivariate Chebyshev polynomial based low-rank representations of wavefront locations and Hadamard ingredients so that block-matrix based fast partial summation can be implemented. Section 7 presents both two-dimensional (2-D) and three-dimensional (3-D) results to demonstrate the performance and accuracy of the new Hadamard integrator. We conclude the paper with some comments in Section 8.

\section{Kirchhoff-Huygens representation formula}
We are interested in solving the following Cauchy problem for the self-adjoint wave equation,
\begin{equation}\label{2.1.0}
  \rho u_{tt}-\nab\cdot(\nu \nab u)=0,\;\; \boldsymbol{x}\in \mathbb{R}^m,\; t>0
\end{equation}
with initial conditions
\begin{equation}\label{2.2.0}
  u(0,\boldsymbol{x})=u^1(\boldsymbol{x}),\,\; u_t(0,\boldsymbol{x})=u^2(\boldsymbol{x}).
\end{equation}
Here $\rho$ and $\nu$ are functions of position $\boldsymbol{x}$.
We look for an integral representation for the wave solution, leading to the Kirchhoff-Huygens representation formula. We give a self-contained derivation here, as the derivation itself sheds some light on how to use it.

Let $u(t,\boldsymbol{x})$ and $v(t,\boldsymbol{x})$ satisfy the non-homogeneous, self-adjoint wave equation:
\begin{equation}\label{2.46}
  \left[\rho(\boldsymbol{x}) u_t(t, \boldsymbol{x})\right]_t-{\nab} \cdot[\nu(\boldsymbol{x}) {\nab} u(t, \boldsymbol{x})]=g_1(t, \boldsymbol{x})
\end{equation}
and
\begin{equation}\label{2.47}
  \left[\rho(\boldsymbol{x}) v_t(t, \boldsymbol{x})\right]_t-{\nab} \cdot[\nu(\boldsymbol{x}) {\nab} v(t, \boldsymbol{x})]=g_2(t, \boldsymbol{x}).
\end{equation}
We will multiply (\ref{2.46}) by $v$ and integrate by parts and then multiply (\ref{2.47}) by $u$ and integrate by parts. On subtraction certain terms will cancel leaving an expression in the form of a divergence to which we will apply the divergence theorem.
So we have
\begin{equation}\label{2.48}
\begin{aligned}
  v g_1&=v\left(\rho u_t\right)_t-v {\nab} \cdot(\nu {\nab} u),\\
  &=\left(\rho v u_t\right)_t-{\nab} \cdot(v \nu {\nab} u)-\rho v_t u_t+\nu {\nab} v \cdot {\nab} u
\end{aligned}
\end{equation}
and
\begin{equation}\label{2.49}
u g_2=\left(\rho u v_t\right)_t-{\nab} \cdot(u \nu {\nab} v)-\rho u_t v_t+\nu {\nab} u \cdot {\nab} v .
\end{equation}
Subtract (\ref{2.49}) from (\ref{2.48}) to get
\begin{equation}\label{2.50}
  v g_1-u g_2=\rho\left(v u_t-u v_t\right)_t-{\nab} \cdot(v \nu {\nab} u-u \nu {\nab} v) .
\end{equation}
Let $V$ be a region of space which does not change in time, and let us integrate (\ref{2.50}) over the cylindrical domain $D$ in space time
\begin{equation}\label{2.51}
  D=[0, T] \times V=\{(t, \boldsymbol{x}) \mid \boldsymbol{x} \in V \;\;\text { and }\;\; t \in[0, T]\}.
\end{equation}
Integrating (\ref{2.50}) over $D$ we get
\begin{equation}\label{2.52}
  \begin{aligned}
\int_D v f-u g \mathrm{~d} D= & \int_0^T \mathrm{~d} t \int_V \rho\left(v u_t-u v_t\right)_t \mathrm{~d} V -\int_0^T \mathrm{~d} t \int_V {\nab} \cdot(v \nu {\nab} u-u \nu {\nab} v) \mathrm{d} V \\
= & {\left[\int_V \rho\left(v u_t-u v_t\right) \mathrm{d} V\right]_0^T } -\int_0^T \mathrm{~d} t \int_S v \boldsymbol{n} \cdot(\nu {\nab} u)-u \boldsymbol{n} \cdot(\nu {\nab} v) \mathrm{d} S\\
=& \int_V \rho(\boldsymbol{x})\left[v(T, \boldsymbol{x}) u_t(T, \boldsymbol{x})-u(T, \boldsymbol{x}) v_t(T, \boldsymbol{x})\right] \mathrm{d} V\\
& -\int_V \rho(\boldsymbol{x})\left[v(0, \boldsymbol{x}) u_t(0, \boldsymbol{x})-u(0, \boldsymbol{x}) v_t(0, \boldsymbol{x})\right] \mathrm{d} V \\
& -\int_0^T \mathrm{~d} t \int_S[v(t, \boldsymbol{x}) \boldsymbol{n} \cdot(\nu {\nab} u(t, \boldsymbol{x}))-u(t, \boldsymbol{x}) \boldsymbol{n} \cdot(\nu {\nab} v(t, \boldsymbol{x}))] \mathrm{d} S,
\end{aligned}
\end{equation}
where $S=\partial V$.

Now we proceed to specialize formula (\ref{2.52}) by choosing $v$ in a special way related to the whole-space Green's function. First consider the Green's function $G(t,\bx_0;\bx)$ which satisfies
\beq{2.53}
[\rho(\bx) G_t(t,\bx_0;\bx)]_t - \nab\bcdot[\nu(\bx)\nab G(t,\bx_0;\bx)] = \dlt(\bx-\bx_0)\,\dlt(t)\,,
\eeq
with initial condition
\beq{2.54}
G(t,\bx_0;\bx) = 0\;\;\mbox{ for }\;\; t<0\,.
\eeq
Alternatively the same $G$ may be specified as the solution of
\beq{2.03}
[\rho(\bx) G_t(t,\bx_0;\bx)]_t - \nab\bcdot[\nu(\bx)\nab G(t,\bx_0;\bx)] = 0\;\;\mbox{ for }\;\; t>0\,,
\eeq
with initial conditions
\beq{2.04}
G(0,\bx_0;\bx) = 0\; \mbox{ and }\; G_t(0,\bx_0;\bx) = \fr{1}{\rho(\bx_0)}\,\dlt(\bx-\bx_0) \,.
\eeq
We will not specify any boundary conditions on $G$.

We then set
\begin{equation}\label{2.55}
  v(t, \boldsymbol{x})=\bar{G}\left(t, \boldsymbol{x}_0; \boldsymbol{x}\right)=G\left(T-t, \boldsymbol{x}_0; \boldsymbol{x}\right).
\end{equation}
We easily verify that
\begin{equation}\label{2.56}
  \left[\rho(\boldsymbol{x}) \bar{G}_t\left(t, \boldsymbol{x}_0; \boldsymbol{x}\right)\right]_t-\nab \cdot\left[\nu(\boldsymbol{x}) \nab \bar{G}\left(t, \boldsymbol{x}_0; \boldsymbol{x}\right)\right]=0 \quad\text { for}\quad t<T,
\end{equation}
with final conditions
\begin{equation}\label{2.57}
  \bar{G}\left(T, \boldsymbol{x}_0; \boldsymbol{x}\right)=0\;,\;\; \bar{G}_t\left(T, \boldsymbol{x}_0; \boldsymbol{x}\right)=-\frac{1}{\rho\left(\boldsymbol{x}_0\right)} \delta\left(\boldsymbol{x}-\boldsymbol{x}_0\right).
\end{equation}
With $v=\bar{G}$ we have already chosen the forcing function $g_2\equiv 0$. Now we take $g_1\equiv 0$ so that the left member of (\ref{2.52}) vanishes. Hence, setting $g_1\equiv g_2\equiv 0$ and $v=\bar{G}$,  (\ref{2.52}) gives
\begin{equation}\label{2.58}
  \begin{aligned}
0= & \int_V \rho(\boldsymbol{x})\left[\bar{G}\left(T, \boldsymbol{x}_0; \boldsymbol{x}\right) u_t(T, \boldsymbol{x})-u(T, \boldsymbol{x}) \bar{G}_t\left(T, \boldsymbol{x}_0; \boldsymbol{x}\right)\right] \mathrm{d} V \\
& -\int_V \rho(\boldsymbol{x})\left[\bar{G}\left(0, \boldsymbol{x}_0; \boldsymbol{x}\right) u_t(0, \boldsymbol{x})-u(0, \boldsymbol{x}) \bar{G}_t\left(0, \boldsymbol{x}_0; \boldsymbol{x}\right)\right] \mathrm{d} V \\
& -\int_0^T \mathrm{~d} t \int_S \bar{G}\left(t, \boldsymbol{x}_0; \boldsymbol{x}\right) \boldsymbol{n} \cdot(\nu {\nab} u(t, \boldsymbol{x}))-u(t, \boldsymbol{x}) \boldsymbol{n} \cdot\left(\nu {\nab} \bar{G}\left(t, \boldsymbol{x}_0; \boldsymbol{x}\right)\right) \mathrm{d} S .
\end{aligned}
\end{equation}
Substituting (\ref{2.55}) and (\ref{2.57}) in (\ref{2.58}) we get
\begin{equation}\label{2.59}
  \begin{aligned}
 u\left(T, \boldsymbol{x}_0\right)&=\int_V \rho(\boldsymbol{x})\left[G\left(T, \boldsymbol{x}_0; \boldsymbol{x}\right) u_t(0, \boldsymbol{x})+u(0, \boldsymbol{x}) G_t\left(T, \boldsymbol{x}_0; \boldsymbol{x}\right)\right] \mathrm{d} V \\
& \quad+\int_0^T \mathrm{~d} t \int_S G\left(T-t, \boldsymbol{x}_0; \boldsymbol{x}\right) \boldsymbol{n} \cdot[\nu {\nab} u(t, \boldsymbol{x})]-u(t, \boldsymbol{x}) \boldsymbol{n} \cdot\left[\nu {\nab} G\left(T-t, \boldsymbol{x}_0; \boldsymbol{x}\right)\right] \mathrm{d} S ,
\end{aligned}
\end{equation}
where the boundary integrals capture information entering the computational domain from outside. Since here we consider the initial value problem with compactly supported initial conditions in a finite time, we can drop the boundary integrals in (\ref{2.59}) by taking a sufficiently large computational domain in space in our formulation without affecting the wave solution, and thus we will do so in the following to obtain
\begin{equation}\label{2.60}
 u\left(T, \boldsymbol{x}_0\right)=\int_V \rho(\boldsymbol{x})\left[G\left(T, \boldsymbol{x}_0; \boldsymbol{x}\right) u_t(0, \boldsymbol{x})+u(0, \boldsymbol{x}) G_t\left(T,\boldsymbol{x}_0; \boldsymbol{x}\right)\right] \mathrm{d} V.
\end{equation}

Differentiating (\ref{2.60}) with respect to time, we get
\begin{equation}\label{2.61}
 u_t\left(T, \boldsymbol{x}_0\right)=\int_V \rho(\boldsymbol{x})\left[G_t\left(T, \boldsymbol{x}_0; \boldsymbol{x}\right) u_t(0, \boldsymbol{x})+u(0, \boldsymbol{x}) G_{tt}\left(T, \boldsymbol{x}_0; \boldsymbol{x}\right)\right] \mathrm{d} V.
\end{equation}

We refer to (\ref{2.60}) and (\ref{2.61}) as the Kirchhoff-Huygens representation formula. Now the question is how to use this formula. When the medium is homogeneous, the Green's function for the wave equation is known so that the formula has been used frequently in practice. However, since the Green's function is usually unknown in an inhomogeneous medium, it is extremely challenging to use this formula efficiently in this case. Therefore, we propose to compute the needed Green's function by developing and implementing a novel Hadamard's asymptotic ansatz. As we will see, using this novel ansatz in the Kirchhoff-Huygens representation formula gives us the Hadamard-Kirchhoff-Huygens (HKH) propagator which is able to propagate highly oscillatory wavefields for a short period of time, but recursively applying this propagator in time yields the {\it Hadamard integrator} to solve time-dependent wave equations globally in time, where caustics are treated implicitly. Moreover, by marching forward in time, we are able to treat spatially overturning waves naturally.

\section{Hadamard's ansatz based local solution}
\subsection{Hadamard's ansatz}
We seek an asymptotic representation of the Green's function of the self-adjoint wave equation which we rewrite here
\begin{equation}\label{2.1}
  \rho u_{tt}-\nab\cdot(\nu \nab u)=0, \boldsymbol{x}\in \mathbb{R}^m, t>0,
\end{equation}
with initial conditions
\begin{equation}\label{2.2}
  u(0,\boldsymbol{x})=0,\, u_t(0,\boldsymbol{x})=\frac{1}{\rho(\boldsymbol{x}_0)}\delta(\boldsymbol{x}-\boldsymbol{x}_0).
\end{equation}

 Essentially, we are looking for the Green's function at the origin
 $\boldsymbol{x}_0 = \boldsymbol{0}$. In what follows we shall regard the dependence upon $\boldsymbol{x}_0$ as understood. Let
\begin{equation}\label{2.3}
  c=\sqrt{\frac{\nu}{\rho}},\, \quad n=\frac{1}{c}=\sqrt{\frac{\rho}{\nu}}.
\end{equation}
We seek an asymptotic representation of $u$ using the Hadamard's ansatz
\begin{equation}\label{2.4}
  u(t, \boldsymbol{x})=\sum_s v_s(\boldsymbol{x}) f_{+}^{s-\frac{m-1}{2}}\left[t^2-\tau^2(\boldsymbol{x})\right],
\end{equation}
where the summation is over all integer values of $s$ using the convention that $ v_s\equiv 0$ for $s < 0$ and that $v_0\neq 0$. $\tau(\boldsymbol{x})$ is the phase function, also known as traveltime, which can be explained as the least travel time at speed $c(\boldsymbol{x})$ from the origin to the point $\boldsymbol{x}$. The generalized function $f_{+}^{\lambda}(\varsigma)$ is called  the Gelfand-Shilov function as described in \cite{gelshi64}. They are defined for $\lambda>-1$ as follows:

\begin{equation}\label{2.5}
  f_{+}^\lambda(\varsigma)=\frac{\varsigma_{+}^\lambda}{\lambda !},
\end{equation}
where
\begin{equation}\label{2.6}
  \varsigma_{+}^\lambda=\left\{\begin{array}{cc}
0, & \text { for } \varsigma<0, \\
\varsigma^\lambda, & \text { otherwise },
\end{array}\right.
\end{equation}
and by analytic continuation for other values of $\lambda$. Thus the support of the function $f_{+}^\lambda\left[t^2-\tau^2(\boldsymbol{x})\right]$ lies within the double cone which is the union of $t \geq \tau$ and $t \leq-\tau$. The poles of $\varsigma_{+}^\lambda$ and of $\lambda$! at the negative integer values of $\lambda$ cancel so that $f_{+}^\lambda(\varsigma)$ is an entire function of $\lambda$. We shall be concerned only with $t>0$ for the wave equation. Also
\beq{SASWE1.06}\ba{ccc}
f_+^{-1}(\varsigma) = \delta(\varsigma),\; & \mbox{and} & \;f_+^{-n}(\varsigma) = \delta^{(n-1)}(\varsigma).
\ea\eeq
We have
\begin{equation}\label{2.7}
  \varsigma f_{+}^{\lambda-1}(\varsigma)=\frac{\varsigma_{+}^\lambda}{(\lambda-1) !}=\lambda \frac{\varsigma_{+}^\lambda}{\lambda !}=\lambda f_{+}^\lambda(\varsigma)
\end{equation}
and the important relationship
\begin{equation}\label{2.8}
  (f_{+}^{\lambda})^{\prime}=f_{+}^{\lambda-1} .
\end{equation}

Writing (\ref{2.1}) in subscript notation we get
\beq{SASWE1.09}
\rho\,\ddot{u}\,-\,(\nu\,u_{,k})_{,k} = 0\;,
\eeq
where $u_{,k}$ indicates the $x_k$-derivative of $u$.
We will calculate successively $u_{,k}$, $\nu\,u_{,k}$, $(\nu u_{,k})_{,k}$, and $\rho\,\ddot{u}$.

Using \rf{2.8},
\beq{SASWE2.02}\ba{rcl}
u_{,k} &=&
\Sm{s}{}(\,v_s\,\tau_{,k}\partial_\tau f_+^{s-\frscr{m-1}{2}} + v_{s,k}\,f_+^{s-\frscr{m-1}{2}}) \\ \\
&=&
\Sm{s}{}(-2 v_s\,\tau\,\tau_{,k}\, f_+^{s-\frscr{m+1}{2}} + v_{s,k}\,f_+^{s-\frscr{m-1}{2}}) .
\ea\eeq
Hence
\beq{SASWE2.03}\ba{rcl}
\nu\,u_{,k} &=&
\Sm{s}{}(-2 v_s\,\nu\,\tau\,\tau_{,k} \,f_+^{s-\frscr{m+1}{2}} +
v_{s,k}\,\nu\,f_+^{s-\frscr{m-1}{2}}) .
\ea\eeq
It follows that
\begin{equation}\label{SASWE2.04}
    \ba{rcl}
(\nu\,u_{,k})_{,k} &=&
\Sm{s}{}[4 v_s\,\nu\,\tau^2\,\tau_{,k} \,\tau_{,k}\,f_+^{s-\frscr{m+3}{2}}
-2 (v_s\,\nu\,\tau\,\tau_{,k})_{,k}\,f_+^{s-\frscr{m+1}{2}} -2\, v_{s,k}\,\nu\,\tau\,\tau_{,k}\,f_+^{s-\frscr{m+1}{2}} +
(v_{s,k}\,\nu)_{,k}\,\,f_+^{s-\frscr{m-1}{2}}] \\ \\
&=&
\Sm{s}{}f_+^{s-\frscr{m+3}{2}}\,[4 v_s\,\nu\,\tau^2\,\tau_{,k} \,\tau_{,k}\,
-2 (v_{s-1}\,\nu\,\tau\,\tau_{,k})_{,k} -2\, v_{s-1,k}\,\nu\,\tau\,\tau_{,k} +
(v_{s-2,k}\,\nu)_{,k}\,]\;.
\ea
\end{equation}

Also, by \rf{2.4},
\beq{SASWE2.05}\ba{rcl}
\rho\,\ddot{u} &=& \Sm{s}{}\rho\,(v_s\,2\,\partial_t(t\,f_+^{s-\frscr{m+1}{2}}) \\
&=& \Sm{s}{}\rho\,\,v_s\,[ 4\,t^2\,f_+^{s-\frscr{m+3}{2}} \,+\,
2\,f_+^{s-\frscr{m+1}{2}}] \\
&=& \Sm{s}{}\rho\,\,v_s\,[ 4\,(t^2-\tau^2)\,f_+^{s-\frscr{m+3}{2}} \,+\,
4\,\tau^2\,f_+^{s-\frscr{m+3}{2}} \,+\,2\,f_+^{s-\frscr{m+1}{2}}] \\
&=& \Sm{s}{}\rho\,\,v_s\,[ 4\,(s-\frtxt{m+1}{2})\,f_+^{s-\frscr{m+1}{2}} \,+\,
4\,\tau^2\,f_+^{s-\frscr{m+3}{2}} \,+\,2\,f_+^{s-\frscr{m+1}{2}}] \\
&=& \Sm{s}{}\rho\,\,v_s\,[4\,\tau^2\,f_+^{s-\frscr{m+3}{2}} \,+\, 4\,(s-\frtxt{m}{2})\,f_+^{s-\frscr{m+1}{2}}] \\
&=&4\,\Sm{s}{}\,f_+^{s-\frscr{m+3}{2}}\,[\rho\,\tau^2\,v_{s}\, +
(s-1-\frtxt{m}{2})\,\rho\,v_{s-1}\,] ,
 \ea\eeq
 where we have also used \rf{2.7}.
Using \rf{SASWE2.04} and \rf{SASWE2.05} in \rf{SASWE1.09} we see that
\beq{SASWE2.06}\ba{rcl}
0\,&=&
\Sm{s}{}f_+^{s-\frscr{m+3}{2}}\,[\,4 v_{s}\,\nu\,\tau^2\,\tau_{,k} \,\tau_{,k}\,
-\,4\,v_{s}\,\rho\,\tau^2\,-2 (v_{s-1}\,\nu\,\tau\,\tau_{,k})_{,k}\,\\
&&\qquad -2\, v_{s-1,k}\,\nu\,\tau\,\tau_{,k}\,-\,4\,(s-1-\frtxt{m}{2})\,\rho\,v_{s-1}   +
(v_{s-2,k}\,\nu)_{,k}\,] .
\ea\eeq
Assuming that \rf{2.4} is an asymptotic power series in powers of $t^2-\tau^2$ near $t=\tau$, we may equate to zero the coefficient of each of the $f_+^{s-\frscr{m+3}{2}}$ for $s=0, 1, 2, \dots$

Thus
\beq{SASWE2.07}\ba{rcl}
0\,&=&
4 v_{s}\,\tau^2\,(\nu\,\tau_{,k} \,\tau_{,k}\, \,-\,\rho)\,-2 (v_{s-1}\,\nu\,\tau\,\tau_{,k})_{,k}\,
-2\, v_{s-1,k}\,\nu\,\tau\,\tau_{,k}\,-\,4\,(s-1-\frtxt{m}{2})\,\rho\,v_{s-1} +
(v_{s-2,k}\,\nu)_{,k} .
\ea\eeq
By setting $s=0$ in \rf{SASWE2.07} and remembering that $v_{-1} = v_{-2} = 0$ we get
\beq{SASWE2.08}
4\,v_0\,\tau^2\,(\nu\,\tau_{,k} \,\tau_{,k}\, \,-\,\rho)\,=\,0 .
\eeq

Since we are assuming that $v_0\neq0$ and $\tau(\bx)\neq0$ except at $\bx=\bzero$, we have
\beq{SASWE2.09}
\nu\,\tau_{,k} \,\tau_{,k}\,-\,\rho\,=\,0 ,
\eeq
which is the eikonal equation
\beq{SASWE2.10}
|\bnab\,\tau|^2 = \fr{\rho}{\nu} = n^2 .
\eeq
Then, equation \rf{SASWE2.06} with $s$ replaced by $s+1$ reduces to the transport equations for $v_s$,
\beq{SASWE2.11}\ba{rcl}
2 (\nu\,\tau\,\tau_{,k}\,v_{s})_{,k}\,
+2\, \nu\,\tau\,\tau_{,k}\,v_{s,k}\,
+\,4\,(s-\frtxt{m}{2})\,\rho\,v_{s}  &=&(v_{s-1,k}\,\nu)_{,k} .
\ea\eeq
Or
\beq{SASWE2.12}\ba{rcl}
4\, \nu\,\tau\,\tau_{,k}\,v_{s,k}\,
+[\nu\,(\tau^2)_{,k}]_{,k}\,v_{s}\,
+\,2\,(2s-m)\,\rho\,v_{s}  &=&(v_{s-1,k}\,\nu)_{,k} .
\ea\eeq

\subsection{Interpretation as energy conservation when $s=0$}
We now study \rf{SASWE2.12} for the leading amplitude, i.e. for $s=0$.  Notice that the leading amplitude in \rf{2.4} is not $v_0$ but $v_0$ multiplied by the amplitude of the leading singularity near $t=\tau$ of $f_+^{-\frscr{m-1}{2}}(t^2-\tau^2)=(t+\tau)^{-\frscr{m-1}{2}}\,f_+^{-\frscr{m-1}{2}}(t-\tau)$.\footnote{This is true even for odd $m$, in which case $f_+^{-\frscr{m-1}{2}} = \dlt^{(\frscr{m-3}{2})}$.  (See Gelfand and Shilov (1964), III, 1.7.)}  But for $t$ near $\tau$ this amplitude is proportional to $\tau^{-\frscr{m-1}{2}}$.  So the true leading amplitude is (proportional to) $u_0 = v_0/\tau^{\frscr{m-1}{2}}$, and we will set $v_0 = \tau^{\frscr{m-1}{2}}\,u_0$ in \rf{SASWE2.11} to get
\beq{SASWE2.13}\ba{rcl}
0\, &=&(u_{0}\,\nu\,\tau^\frscr{m+1}{2}\,\tau_{,k})_{,k}\, +
\, u_{0,k}\,\nu\,\tau^\frscr{m+1}{2}\,\tau_{,k}\, +
\,\frtxt{m-1}{2}\,u_{0}\,\nu\,\tau^\frscr{m-1}{2}\,\tau_{,k}\, \tau_{,k}\,-\,m\,\rho\,\tau^\frscr{m-1}{2}\,u_{0} \\
&=&(u_{0}\,\nu\,\tau^\frscr{m+1}{2}\,\tau_{,k})_{,k}\, +
\, u_{0,k}\,\nu\,\tau^\frscr{m+1}{2}\,\tau_{,k}\,-\,\frscr{m+1}{2}\,\,u_{0}\,\nu\,\tau^\frscr{m-1}{2}\,\tau_{,k}\, \tau_{,k} \\
&=&\tau^\frscr{m+1}{2}\,[(u_{0}\,\nu\,\tau_{,k})_{,k}\, +
\, u_{0,k}\,\nu\,\tau_{,k}] .
\ea\eeq
So
\beq{SASWE2.14}
(u_{0}\,\nu\,\tau_{,k})_{,k}\, +
\, u_{0,k}\,\nu\,\tau_{,k} = 0 .
\eeq
On multiplying by $u_0$ we get
\beq{SASWE2.15}
u_0\,(u_{0}\,\nu\,\tau_{,k})_{,k}\, +
\, u_{0,k}\,u_0\,\nu\,\tau_{,k} = 0 .
\eeq
But this is
\beq{SASWE2.16}
(\nu\,u_{0}^2\,\tau_{,k})_{,k} = 0 .
\eeq
I.e.
\beq{SASWE2.17}
(\rho\,u_{0}^2\,c^2\,\tau_{,k})_{,k} = \bnab\bcdot(\rho\,u_{0}^2\,c^2\,\bnab\tau) = 0.
\eeq
But $\rho\,u_{0}^2$ can be thought of as the energy density\footnote{The energy density is $\hftxt\,\rho\,\dot{u_0}^2\,+\,\hftxt\,\nu\,|\bnab u_0|^2$, but these two terms, which are analogous to kinetic and potential energy densities, are equal in this leading asymptotic term.}
and $c^2\,\bnab\tau$ is the group velocity vector with magnitude $c$ in the direction of the ray, so that $\rho\,u_{0}^2\,c^2\bnab\tau$ is the energy flux vector.  Thus the divergence of the energy flux is zero, and so energy is conserved for the leading-order singularity, and energy flux is conserved along tubes of rays as it should be.  This also verifies that \rf{SASWE2.12} leads to the conventional transport equation for the leading term, which always has this interpretation as energy conservation.

\subsection{Solution of the transport equations}
In this subsection, we assume that the traveltime $\tau$ has been found in the source neighborhood by the method of characteristics.
\subsubsection{Analytic form of $v_0$}
By method of characteristics, along a ray traced out from the source $\boldsymbol{0}$ to $\boldsymbol{x}$, the directional derivative operator along the traveltime $\tau$ satisfies
\begin{equation}\label{2.22}
  n^2 \frac{\mathrm{d}}{\mathrm{d} \tau}=\nab \tau \cdot \nab=\tau_{,k}\partial_k .
\end{equation}

Using \rf{SASWE2.12} with $s=0$ and \rf{2.22} we may write the equation for $v_0$ in the form
\beq{SASWE2.19}
4\,\rho\,\tau\,\fr{\drm v_0}{\drm\tau}\,+\,v_0\,[\bnab\bcdot(\nu\,\bnab\,\tau^2)\,-\,2\,m\,\rho] = 0 ,
\eeq
from which we see that
\beq{SASWE2.20}
\fr{\drm\log v_0}{\drm\tau} = \fr{1}{v_0}\,\fr{\drm v_0}{\drm \tau}
= -\,\fr{\bnab\bcdot(\nu\,\bnab\,\tau^2) - 2\,m\,\rho}{4\,\rho\,\tau} .
\eeq

So,
\begin{equation}\label{2.25}
  \log v_0=-\int_{\mathcal{R}(\boldsymbol{x})} \frac{{\nab} \cdot\left(\nu {\nab} \tau^2\right)-2m \rho}{4 \tau \rho} \mathrm{d} \tau+\log \left[v_0(\mathbf{0})\right],
\end{equation}
and
\begin{equation}\label{2.26}
  \begin{aligned}
v_0(\boldsymbol{x})&=v_0(\boldsymbol{0}) \exp \left\{-\int_{\mathcal{R}(\boldsymbol{x})} \frac{{\nab} \cdot\left(\nu {\nab} \tau^2\right)-2m \rho}{4 \rho \tau} \mathrm{d} \tau\right\} \\
&=v_0(\boldsymbol{0}) h_0(\bx).
\end{aligned}
\end{equation}
Here $\mathcal{R}(\boldsymbol{x})$ is the segment of ray joining $\boldsymbol{ 0}$ to the point $\boldsymbol{x}$. The dependency upon $\boldsymbol{x}$ arises from the dependency of $\mathcal{R}(\boldsymbol{x})$, which is parameterized by $\tau$. It turns out that $v_0(\boldsymbol{0})$ does not depend upon the initial direction of the ray at $\boldsymbol{0}$, and in fact $v_0(\boldsymbol{x})$ is analytic in $\boldsymbol{x}$ at $\boldsymbol{0}$ if $n(\bx)$ is also.
We have defined
\begin{equation}\label{2.26a}
  \begin{aligned}
h_0(\boldsymbol{x})&=\exp \left\{-\int_{\mathcal{R}(\boldsymbol{x})} \frac{{\nab} \cdot\left(\nu {\nab} \tau^2\right)-2m \rho}{4 \rho \tau} \mathrm{d} \tau\right\},
\end{aligned}
\end{equation}
which satisfies \eqref{SASWE2.19} for $v_0$ with $h_0(\boldsymbol{0})=1$.

\subsubsection{Analytic forms of $v_s$ for $s\geq 1$}
We begin by restating \rf{SASWE2.12}
\beq{SASWE3.01}
4\,\rho\,\tau\,\fr{\drm v_{s}}{\drm \tau}\,+\,v_{s}\,[\nab\bcdot(\nu\,\nab\,\tau^2)\,+\,2\,(2s-m)\,\rho] = \nab\bcdot(\nu\nab v_{s-1}) .
\eeq
Let us first consider the solution $h_{s}$ to the homogeneous form of \rf{SASWE3.01}
\beq{SASWE3.02}
4\,\rho\,\tau\,\fr{\drm h_s}{\drm \tau}\,+\,h_s\,[\nab\bcdot(\nu\,\nab\,\tau^2)\,+\,2\,(2s-m)\,\rho] = 0,
\eeq
which leads to
\beq{SASWE3.03}\ba{rcl}
\fr{\drm\log h_s}{\drm \tau}
&=& -\,\fr{\nab\bcdot(\nu\,\nab\,\tau^2)\, - 2\,m\,\rho}{4\,\rho\,\tau}
\,-\,\fr{s}{\tau}\\
&=& \fr{\drm\log h_0}{\drm \tau}\,-\,\fr{s}{\tau} \\
&=& \fr{\drm\log \tau^{-s}\,h_0}{\drm \tau}\;.
\ea\eeq
So we may take
\beq{SASWE3.04}
h_s = \tau^{-s}\,h_0 .
\eeq
To solve \rf{SASWE3.01} we set
\beq{SASWE3.05}
v_{s} = h_s\,w_s =  \tau^{-s}\,h_0\,w_s.
\eeq
Substituting this into \rf{SASWE3.01} and taking \rf{SASWE3.02} into account we have
\beq{SASWE3.06}
4\,\rho\,\tau^{1-s}\,h_0\,\fr{\drm w_s}{\drm \tau}
=  \nab\bcdot(\nu\nab v_{s-1}),
\eeq
leading to
\beq{SASWE3.07}
w_s = \Int{\calR(\scriptsize\bx)}{}\fr{\tau^{s-1}\,\nab\bcdot(\nu\nab v_{s-1})}
{4\,\rho\,v_0}\,\drm\tau + c_{s},
\eeq
and finally we have
\beq{SASWE3.08}
v_{s}(\bx) = \tau^{-s}\,h_0(\bx)\,\Biggl[\Int{\calR(\scriptsize\bx)}{}\fr{\tau^{s-1}\,\nab\bcdot(\nu\nab v_{s-1})}{4\,\rho\,h_0}\,\drm\tau\,+ \,c_{s}\Biggr] .
\eeq
Let us consider the behavior of $v_{s}$ as $|\bx|\to 0$ assuming that $v_{s-1}$ and its derivatives are finite.  We easily see that for $s\geq 1$
\beq{SASWE3.09}
v_{s} \quad\simeq\quad \left\{\ba{cc}c_{s}\,\tau^{-s}&\mbox{for $c_{s}\neq 0$} , \\ \\
\fr{\nab\bcdot(\nu\nab v_{s-1})}{4\,s\,\rho}\Biggr|_{{\scriptstyle \bx}=\bzero}
 & \mbox{for $c_{s}= 0$} .
\ea\right.
\eeq
We choose $c_s=0$ for all $s\geq1$ in order for $v_s$ to be finite.  It will also be analytic in $\bx$ according to Babich\cite{bab65} just below his equation (8), to which \rf{SASWE3.08} should reduce when $\nu=1$ and $\rho = n^2$, but there appears to be an error in Babich's equation (8).

Thus we see that the $v_s$ are determined by, and depend linearly upon, $v_0(\bzero)$, whose value we obtain from the initial conditions \rf{2.2}.

\subsubsection{Initialization of $v_0(\boldsymbol{0})$ for $m$ even }
To obtain the initialization of $v_0(\boldsymbol{0})$, we will follow closely the presentation in the Example in Chapter III, Section
1.6 of Gelfand and Shilov \cite{gelshi64}. We will consider only the first term of the asymptotic series (\ref{2.4}) since it will subsequently be clear that later terms contribute zero to the initial values of $u$ and $u_t$.

We write the leading term in the series (\ref{2.4}) as
\begin{equation}\label{2.27}
  G_0(t, \bx_0;\boldsymbol{x})=v_0(\bx_0;\boldsymbol{x}) f_{+}^{-\frac{m-1}{2}}\left[t^2-\tau^2(\bx_0;\boldsymbol{x})\right].
\end{equation}
When there is no ambiguity, we suppress source $\bx_0 = \bzero$ in the display to simplify notation. Let
\begin{equation*}
    c_0=c(\boldsymbol{0}),\quad n_0=n(\boldsymbol{0}),\quad \nu_0=\nu(\boldsymbol{0}),\quad\mbox{and}\quad \rho_0=\rho(\boldsymbol{0}).
\end{equation*}
We first introduce the geodesic (ray) polar transformation for a given source $\bx_0$,
\begin{equation}\label{gec}
    P[\bx_0]: \bx \rightarrow (\tau,\boldsymbol{\omega}),
\end{equation}
where $\tau$ is the traveltime, $\boldsymbol{\omega}\in \mathbb{S}^{m-1}$ is the take-off angle of the ray, and $(\tau,\boldsymbol{\omega})$ is the geodesic polar coordinates. Within any neighborhood of $\bx_0$ not containing any caustics other than $\bx_0$, there is one and only one ray connecting $\bx$ and $\bx_0$, which means $P[\bx_0]$ is well-defined and one-to-one, and $\boldsymbol{x}$ is a smooth function of the point $\boldsymbol{y}=\bx_0+\tau \boldsymbol{\omega}$. To facilitate our following discussions, we also write down the volume element as the following,
\begin{equation}
    \mathrm{d} \bx=\mathrm{d} s \mathrm{d} S=c \mathrm{d} \tau \mathrm{d} S=c |\frac{\partial S}{\partial \boldsymbol{\omega}}|\mathrm{d} \tau \mathrm{d}\boldsymbol{\omega},
\end{equation}
where $s$ is the arc length along the ray, $ \mathrm{d} S=\left|\frac{\partial S}{\partial \boldsymbol{\omega}}\right| \mathrm{d} \boldsymbol{\omega}$ is the element of area cut out on the wave front $\tau=$ const. by rays emanating from the solid angle element $\mathrm{d} \omega$ at the source. We will provide a specific expression for the Jacobian $\left|\frac{\partial S}{\partial \boldsymbol{\omega}}\right|$ later. However, for now, we will use the following approximation when $\tau\rightarrow 0$,
\begin{equation}\label{2.29}
  \left|\frac{\partial S}{\partial \boldsymbol{\omega}}\right|=c^{m-1}_0\tau^{m-1}(1+O(\tau)),
\end{equation}
which can be obtained by comparing with the polar coordinates in a homogeneous medium near the source.

We now apply \rf{2.27} to a test function $\phi(\boldsymbol{x})\in C_c^{\infty}(\mathbb{R}^m)$ as follows,
\begin{equation}\label{2.28}
\begin{aligned}
\int G_0(t, \bx_0;\boldsymbol{x}) \phi(\boldsymbol{x}) \mathrm{d} \boldsymbol{x} & =\int v_0(\bx_0;\boldsymbol{x}) f_{+}^{-\frac{m-1}{2}}\left[t^2-\tau^2(\bx_0;\boldsymbol{x})\right] \phi(\boldsymbol{x}) \mathrm{d} \boldsymbol{x} \\
& =\int v_0(\boldsymbol{x}) f_{+}^{-\frac{m-2}{2}-\frac{1}{2}}\left[t^2-\tau^2(\boldsymbol{x})\right] \phi(\boldsymbol{x}) \mathrm{d} \boldsymbol{x} \\
& =\int\left(-\frac{\partial}{2 \tau \partial \tau}\right)^{\frac{m-2}{2}} f_{+}^{-\frac{1}{2}}\left[t^2-\tau^2(\boldsymbol{x})\right] v_0[\boldsymbol{x}(\tau, \boldsymbol{\omega})] \phi[\boldsymbol{x}(\tau, \boldsymbol{\omega})]\left|\frac{\partial S}{\partial \boldsymbol{\omega}}\right| c \mathrm{~d} \tau \mathrm{d} \boldsymbol{\omega},
\end{aligned}
\end{equation}
where
\beq{2.30}\ba{cc}
  k=\fr{m-2}{2}\; \text { and }\; \psi_1(\tau, \boldsymbol{\omega})=v_0[\boldsymbol{x}(\tau, \boldsymbol{\omega})] \phi[\boldsymbol{x}(\tau, \boldsymbol{\omega})]\left|\frac{\partial S}{\partial \boldsymbol{\omega}}\right| .
\ea\eeq
Then as $|\tau|\rightarrow 0$, using \rf{2.29} we have
\begin{equation}\label{2.31}
  \psi_1(\boldsymbol{x})=c^{m-1}_0 v_0(\mathbf{0}) \phi(\mathbf{0}) \tau^{m-1}[1+O(\tau)],
\end{equation}
leading to
\begin{equation}\label{2.32}
  \begin{aligned}
 \int G_0(t, \bx) \phi(\bx) \mathrm{d} \boldsymbol{x}&= c^m_0 v_0(\mathbf{0}) \phi(\mathbf{0}) \int_{\Omega_m} \int_0^{\infty}\left(-\frac{\partial}{2 \tau \partial \tau}\right)^k f_{+}^{-\frac{1}{2}}\left[t^2-\tau^2(\boldsymbol{x})\right]\tau^{m-1}[1+O(\tau)] \mathrm{d} \tau \mathrm{d} \boldsymbol{\omega} \\
&= \frac{1}{2} c^m_0 v_0(\mathbf{0}) \phi(\mathbf{0}) \int_{\Omega_m} \int_0^{\infty} f_{+}^{-\frac{1}{2}}\left[t^2-\tau^2(\boldsymbol{x})\right]
 \left(\frac{\partial}{2 \tau \partial \tau}\right)^k \tau^{m-2}[1+O(\tau)] \mathrm{d} \tau^2 \mathrm{~d} \boldsymbol{\omega}\\
&= \frac{1}{2} c^m_0 v_0(\mathbf{0}) \phi(\mathbf{0}) \int_{\Omega_m} \int_0^{\infty} f_{+}^{-\frac{1}{2}}\left[t^2-\tau^2(\boldsymbol{x})\right] \left(\frac{\partial}{2 \tau \partial \tau}\right)^k \tau^{2k}[1+O(\tau)] \mathrm{d} \tau^2 \mathrm{~d} \boldsymbol{\omega}.
\end{aligned}
\end{equation}
Following Gelfand and Shilov \cite{gelshi64} (Chapter III, Section 1.6) we remark that each application of $\fr{\partial}{2\,\tau\,\partial \tau}$ reduces the lowest power of $\tau$ by two and comes from applying all the derivatives to the factor $\tau^{2 k}$, i.e. performing $\left(\fr{\partial}{2\,\tau\,\partial \tau}\right)^k \tau^{2 k}$. But after $k$ applications of $\fr{\partial}{2\,\tau\,\partial \tau}$ the lowest power of $\tau$ is just $\tau^0=1$, so that, as $t \rightarrow 0$, we are left with
\begin{equation}\label{2.33}
  \begin{aligned}
& \quad \frac{1}{2} c^m_0 v_0(\mathbf{0}) \phi(\mathbf{0}) \int_{\Omega_m} \int_0^{\infty} f_{+}^{-\frac{1}{2}}\left[t^2-\tau^2(\boldsymbol{x})\right]\left(\frac{\partial}{2 \tau \partial \tau}\right)^k \tau^{2 k}[1+O(\tau)] \mathrm{d} \tau^2 \mathrm{~d} \boldsymbol{\omega} \\
& =\frac{(2 k) ! !}{2^k\left(-\frac{1}{2}\right) !} c^m_0 \omega_m v_0(\mathbf{0}) \phi(\mathbf{0}) \int_0^t\left[t^2-\tau^2(\boldsymbol{x})\right]^{-\frac{1}{2}} \tau \mathrm{d} \tau[1+O(\tau)] \\
& =\frac{k !}{\left(-\frac{1}{2}\right) !} c^m_0\omega_m v_0(\boldsymbol{0}) \phi(\boldsymbol{0}) t[1+O(t)],
\end{aligned}
\end{equation}
where $\omega_m$ is the area of the unit sphere in $\mathbb{R}^m$.
So as $t\rightarrow 0$
\begin{equation}\label{2.34}
  \int G_0(t, \boldsymbol{x}) \phi(\boldsymbol{x}) \mathrm{d} \boldsymbol{x} \rightarrow 0,
\end{equation}
and the initial condition of $u$ is satisfied.

We now consider the initial condition of $u_t$. Differentiating (\ref{2.32}) with respect to $t$ and using $\dot{G}_0$ to represent the derivative of the leading term, we obtain
\begin{equation}\label{2.35}
\int \dot{G}_0(t, \bx) \phi(\boldsymbol{x}) \mathrm{d} \boldsymbol{x}=v_0(\boldsymbol{0}) \phi(\boldsymbol{0}) t c_0^m \omega_m \int_0^{\infty}\left(-\frac{\partial}{2 \tau \partial \tau}\right) f_{+}^{-\frac{1}{2}}\left[t^2-\tau^2(\boldsymbol{x})\right]\left(\frac{\partial}{2 \tau \partial \tau}\right)^k \tau^{2 k}[1+O(\tau)] \mathrm{d} \tau^2 .
\end{equation}
The lowest power of $\tau$ is again generated from
\begin{equation}\label{2.36}
  \omega_m\left(\frac{\partial}{2 \tau \partial \tau}\right)^k \tau^{2 k}=\frac{(2 k) ! ! \omega_m}{2^k}=k ! \omega_m.
\end{equation}
But it is known \cite{gelshi64} that
\begin{equation}\label{2.37}
  \omega_m=\frac{2 \pi^{\frac{1}{2} m}}{\left(\frac{m}{2}-1\right) !}=\frac{2 \pi^{k+1}}{k !} .
\end{equation}
Hence
\begin{equation}\label{2.38}
  \omega_m\left(\frac{\partial}{2 \tau \partial \tau}\right)^k \tau^{2 k}=2 \pi^{k+1}.
\end{equation}

Using (\ref{2.38}) in (\ref{2.35}) we have
\begin{equation}\label{2.39}
  \begin{aligned}
\int \dot{G}_0(t, \boldsymbol{x}) \phi(\boldsymbol{x}) \mathrm{d} \boldsymbol{x} & =2 \pi^{k+1} v_0(\mathbf{0}) \phi(\mathbf{0}) t c_0^m \int_0^{\infty}\left(-\frac{\partial}{\partial \tau^2}\right) f_{+}^{-\frac{1}{2}}\left[t^2-\tau^2\right] \mathrm{d} \tau^2[1+O(t)] \\
& =2 \pi^{k+1} v_0(\mathbf{0}) \phi(\mathbf{0}) t c_0^m \int_0^{\infty}\left(-\frac{\partial}{\partial \tau^2}\right) \frac{\left[t^2-\tau^2\right]_{+}^{-\frac{1}{2}}}{\left(-\frac{1}{2}\right) !} \mathrm{d} \tau^2[1+O(t)] \\
& =-\left.2 \pi^{k+\frac{1}{2}} v_0(\mathbf{0}) \phi(\mathbf{0}) t c_0^m\left[t^2-\tau^2\right]_{+}^{-\frac{1}{2}}\right|_{\tau=\infty} ^0[1+O(t)] \\
& =2 \pi^{k+\frac{1}{2}} v_0(\mathbf{0}) \phi(\mathbf{0}) t c_0^m \frac{1}{t}[1+O(t)] \\
& \rightarrow 2 \pi^{k+\frac{1}{2}} v_0(\mathbf{0}) \phi(\mathbf{0}) c_0^m,
\end{aligned}
\end{equation}
where we have used $(-\frac{1}{2})!=\pi^{\frac{1}{2}}$. Hence
\begin{equation}\label{2.40}
  \lim _{t \rightarrow 0} \dot{G}_0(t, \bx)=2 c^m_0 \pi^{k+\frac{1}{2}} v_0(\mathbf{0}) \delta(\boldsymbol{x})=2 c^m_0 \pi^{\frac{m-1}{2}} v_0(\boldsymbol{0}) \delta(\boldsymbol{x}).
\end{equation}

Then if we take
\begin{equation}\label{2.41}
  v_0(\mathbf{0})=\frac{n_0^m}{2 \rho_0 \pi^{\frac{m-1}{2}}}=\frac{n_0^{m-2}}{2 \nu_0 \pi^{\frac{m-1}{2}}}\, ,
\end{equation}
we obtain
\begin{equation}\label{2.42}
  \lim _{t \rightarrow 0} u_t(t, \boldsymbol{x})=\lim _{t \rightarrow 0} \dot{G}_0(t, \boldsymbol{x})=\frac{1}{\rho} \delta(\boldsymbol{x}),
\end{equation}
as required by the initial condition of $u_t$.

\subsubsection{Initialization of $v_0(\boldsymbol{0})$ for $m$ odd }
Similar to even $m$, we only consider the leading term in the series (\ref{2.4})
\beq{SASWE4.01}
G_0(t,\bx) = v_0(\bx)\,f_+^{-\frscr{m-1}{2}}[t^2-\tau^2(\bx)] \;.
\eeq
Since $G_0(t,\bx)$ is zero for $\tau(\bx)>t$ and the region of $\bx$-space  for which $\tau(\bx)\leq t$ shrinks to the origin as $t\rightarrow0$, the limiting value of $G_0(t,\bx)$ is either zero or some combination of derivatives of $\dlt(\bx)$.
In order to study the initial value of $G_0(t,\bx),$ we apply it to a test function $\phi(\bx)\in C_c^{\infty}(\mathbb{R}^m)$ for fixed $t$ as follows,
\beq{SASWE4.02}\ba{rcl}
\Int{}{} G_0(t,\bx)\,\phi(\bx)\,\drm\bx
&=& \Int{}{} v_0(\bx)\,f_+^{-\frscr{m-1}{2}}[t^2-\tau^2(\bx)]
\,\phi(\bx)\,\drm\bx
\\ \\
&=& \Int{}{} v_0(\bx)\,\dlt^{(\frscr{m-3}{2})}[t^2-\tau^2(\bx)]
\,\phi(\bx)\,\drm\bx
\\ \\
&=& \Int{}{} \dlt^{(\frscr{m-3}{2})}[t^2-\tau^2]
v_0[\bx(\tau,\bom)]\,\phi[\bx(\tau,\bom)]\,\Big|\fr{\partial S}{\partial\bld\omega}\Big|\,c\,\drm \tau\,\drm\bld{\omega}
\\ \\
&=& \Int{\Omega_m}{}\,\Int{0}{\infty}\dlt^{(k)}[t^2-\tau^2]
\,\psi_2[\tau,\bom]\,c\,\drm \tau\,\drm\bld{\omega}\;,
\ea\eeq
where
\beq{SASWE4.03}\ba{cc}
k\,=\,\fr{m-3}{2}\,, \;\;\mbox{and}\; & \psi_1(\tau,\boldsymbol{\omega})
\,=v_0[\boldsymbol{x}(\tau, \boldsymbol{\omega})] \phi[\boldsymbol{x}(\tau, \boldsymbol{\omega})]\left|\frac{\partial S}{\partial \boldsymbol{\omega}}\right|.
\ea\eeq
Using (\ref{2.29}), we have
\beq{SASWE4.04}
 \psi_1(\bx)
\,=\,c_0^{m-1}\,v_0(\bzero)\,\phi(\bzero)\,\tau^{m-1}[1+O(\tau)]\;,
\eeq
leading to
\beq{SASWE4.05}\ba{rcl}
\Int{}{} G_0(t,\bx)\,\phi(\bx)\,\drm\bx
&=&\,c_0^{m}\, \Int{\Omega_m}{}\,\Int{0}{\infty}\dlt^{(k)}[t^2-\tau^2]
\,v_0(\bzero)\,\phi(\bzero)\,\tau^{m-1}[1+O(\tau)]\,\drm \tau\,\drm\bld{\omega}
\\ \\
&=&\,\hf\,c_0^{m}\, \Int{\Omega_m}{}\,\Int{0}{\infty}\dlt^{(k)}[t^2-\tau^2]
\,v_0(\bzero)\,\phi(\bzero)\,\tau^{m-2}[1+O(\tau)]\,\drm \tau^2\,\drm\bld{\omega}
\\ \\
&=&\frac{(-1)^{k}}{2}\,c_0^{m}\Int{\Omega_m}{}\,\Int{0}{\infty}\dlt^{(k)}[\tau^2-t^2]
\,v_0(\bzero)\,\phi(\bzero)\,\tau^{m-2}[1+O(\tau)]\,\drm \tau^2\,\drm\bld{\omega}
\\ \\
&=&\frac{(-1)^{k}(-1)^k}{2}\,c_0^{m}\, \Int{\Omega_m}{}
\bigg(\fr{\partial}{2\,\tau\,\partial \tau}\bigg)^k
\,\bigg[v_0(\bzero)\,\phi(\bzero)\,\tau^{2k+1}[1+O(\tau)]\bigg]_{\tau=t}\,\drm\bld{\omega}\;.
\ea\eeq
 Each application of $\fr{\partial}{2\,\tau\,\partial \tau}$ reduces the lowest power of $\tau$ by two and comes from performing $(\fr{\partial}{2\,\tau\,\partial \tau})^k\tau^{2k+1}$.  But after $k$ applications of $\fr{\partial}{2\,\tau\,\partial \tau}$ the lowest power of $\tau$ is just $\tau$, and so as $t\rightarrow0$,
\beq{SASWE4.06}
\Int{}{} G_0(t,\bx)\,\phi(\bx)\,\drm\bx \,\to\,0\;.
\eeq
The other terms in the series \rf{2.4} yield higher powers of $\tau$ and therefore also tend to zero.  Thus, if $v(\bzero)$ is finite, the first initial condition of \rf{2.2}  is satisfied.

We next consider the limit of $\dot{u}(t,\bx)$ as $t\rightarrow0$.
  Differentiating \rf{SASWE4.05} with respect to $t$ we obtain
\beq{SASWE4.07}\ba{rcl}
\Int{}{} \dot{G}_0(t,\bx)\,\phi(\bx)\,\drm\bx
&=&c_0^{m}\,t\, \Int{\Omega_m}{}
\bigg(\fr{\partial}{2\,\tau\,\partial \tau}\bigg)^{k+1}
\,\bigg[v_0(\bzero)\,\phi(\bzero)\,\tau^{2k+1}[1+O(\tau)]\bigg]_{\tau=t}\,\drm\bld{\omega}\;.
\ea\eeq
This time the lowest power of $\tau$ is generated from
\beq{SASWE4.08}
\Int{\Omega_m}{}\,\bigg(\!\fr{\partial}{2\,\tau\,\partial \tau}\!\bigg)^{k+1}\,\,\tau^{2k+1}\,\drm\bom
\,=\,\omega_m\,\bigg(\!\fr{\partial}{2\,\tau\,\partial \tau}\!\bigg)^{k+1}\,\tau^{2k+1}
\,=\,\fr{(2k+1)!!\,\omega_m}{2^{k+1}\,\tau}\;,
\eeq
where $\omega_m$ is the area of the unit sphere on $\mathbb{R}^m$ and
\beq{SASWE4.09}
\omega_m\,=\,\fr{2\pi^{\hfscr\,m}}{(\frtxt{m}{2}-1)!}
\,=\,\fr{2^{k+2}}{(2k+1)!!}\,\pi^{k+1}\;.
\eeq
Hence
\beq{SASWE4.10}
\omega_m\,\bigg(\!\fr{\partial}{2\,\tau\,\partial \tau}\!\bigg)^{k+1}\,\,\tau^{2k+1}
\,=\,\fr{2\,\pi^{k+1}}{\tau}\;.
\eeq
Using this in \rf{SASWE4.07} we see that
\beq{SASWE4.11}
\lim_{t\rightarrow 0}\,\Int{}{} \dot{G}_0(t,\bx)\,\phi(\bx)\,\drm\bx
\,=\,2\,c_0^m\,\pi^{k+1}\,v_0(\bzero)\,\phi(\bzero)
\,=\,2\,c_0^m\,\pi^{k+1}\,v_0(\bzero)\,\phi(\bzero)\;,
\eeq
so that
\beq{SASWE4.12}
\lim_{t\rightarrow0}\dot{G}_0(t,\bx)
\,=\,2\,c_0^m\,\pi^{k+1}\,v_0(\bzero)\,\dlt(\bx)
\,=\,2\,c_0^m\,\pi^{\frscr{m-1}{2}}\,v_0(\bzero)\,\dlt(\bx)\;.
\eeq
Then if we take
\beq{SASWE4.13}
v_0(\bzero)
\,=\,\fr{1}{2\rho_0\,c_0^m\,\pi^{\frscr{m-1}{2}}}
\,=\,\fr{1}{2\nu_0\,c_0^{m-2}\,\pi^{\frscr{m-1}{2}}}
\,=\,\fr{\,n_0^{m-2}}{2\nu_0\,\pi^{\frscr{m-1}{2}}}
\,=\,\fr{\,n_0^m}{2\rho_0\,\pi^{\frscr{m-1}{2}}}\,,
\eeq
we obtain
\beq{SASWE4.14}
\lim_{t\rightarrow0}\dot{u}(t,\bx)
\,=\,\lim_{t\rightarrow0}\dot{G}_0(t,\bx)
\,=\,\fr{1}{\rho_0}\,\dlt(\bx),
\eeq
as required by the second initial condition of \rf{2.2}. Let us define the constant
\beq{SASWE4.141}
C_m
\,\equiv\,\fr{\,n_0^m}{2\rho_0\,\pi^{\frscr{m-1}{2}}}\;.
\eeq
Compare \rf{2.41} with \rf{SASWE4.13} and with equation (5) of Babich \cite{bab65}.  So this general formula holds for all $m$, even and odd.  Gelfand and Shilov recommend using the method of descent from odd to even dimensions \cite{gelshi64}, but this is not so straightforward in a nonuniform medium.

\subsection{Local, short-time validity of Hadamard's ansatz}
In an inhomogeneous, isotropic medium, rays issuing from a point source will refocus away from the source with high probability \cite{whi84} so that multiple rays pass through some physical locations, where the traveltime (a.k.a eikonal or phase) and amplitude functions consist of multiple branches and thus are multivalued; moreover, caustics occur at ray envelopes, where the traveltime function transits from one branch to another, and the amplitude function becomes infinite and is thus burning\footnote{The word `caustic' comes from Greek via Latin, meaning `combustible'.} at such transition points (caustics); for examples of caustics, see \cite{ben99,benlafsensol03} for caustic folds and \cite{qialeu04,leuqiaosh04,qialeu06} for caustic cusps. However, the validity of geometric optics including the Hadamard's ansatz is based on the tacit, implicit assumption that the traveltime function is a well-defined single-valued function. Once the traveltime function becomes multivalued, the geometric-optics ansatz has to be modified to produce correct asymptotic solutions according to the structure of caustics  \cite{lud66,masfed81,babbul72}; since the caustic structure of wave phenomenon is unknown a prior, such modifications are nontrivial and infeasible in general. Therefore, since it cannot easily accommodate possibly multi-branched traveltime functions directly, the Hadamard's ansatz can be used only locally around the source point before a caustic occurs. Then the question that we ask is: how to handle caustics in practice? To go beyond caustics at the level of wave propagation, we will use the Huygens secondary-source principle in the form of the Kirchhoff-Huygens representation formula of wave solutions so that we can treat caustics implicitly rather than `stare' at them directly.

The crux of the matter is the following fact: in isotropic media such as considered here, there is a neighborhood of the source in which no caustic occurs except the source itself \cite{luqiabur16}. In terms of traveltime $\tau$, such a caustic-free neighborhood of the source can be characterized as
$\{\boldsymbol{x}:\tau(\bx_0;\boldsymbol{x})< \bar{T}(\boldsymbol{x}_0)\}$,
where $\bar{T}(\boldsymbol{x}_0)$ is the time when the first caustic transpires.
In this caustic-free neighborhood, the Hadamard asymptotic expansion (\ref{2.4}) of the Green's function is valid for $t<\bar{T}(\boldsymbol{x}_0)$. We do not need to find $\bar{T}(\boldsymbol{x}_0)$ exactly. Since all values of $\boldsymbol{x}_0$ are used by the propagator defined below, we set $\bar{T}$ to be the minimum of $\bar{T}(\boldsymbol{x}_0)$ as $\boldsymbol{x}_0$ varies over some relevant domain, where this domain should not be too large so that $\bar{T}$ might not be too small.
Based on this observation, we incorporate the short-time asymptotic Green's function into the time-domain Kirchhoff-Huygens representation formula to define a short-time $\Delta t$ HKH propagator, where $\Delta t < \bar{T}$ only depends on the medium and is independent of the initial data.

%%%%%%%%%%%%%%%%%%%%%%%%%%%%%%%%%%%%%%%%%%%
\section{Hadamard-Kirchhoff-Huygens (HKH) propagator}
By incorporating the leading term of Hadamard's ansatz \rf{2.4} into Kirchhoff-Huygens representation formulas, we develop the Hadamard-Kirchhoff-Huygens (HKH) propagator to propagate the highly oscillatory wavefields for a short period of time.

To implement this propagator efficiently, we need to integrate Gelfand-Shilov generalized functions over curved or surface wavefronts by using geodesic polar coordinates; thus, we start by providing an exact expression for the Jacobian of the geodesic polar transformation,
\begin{equation}\label{2.68.1}
  \left|\frac{\partial S}{\partial \omega}\right|=\frac{\tau^{m-1}}{4 \rho_0 c_0^m \pi^{m-1} \rho c v_0^2}.
\end{equation}
Its derivation is detailed in Appendix. With the exact Jacobian (\ref{2.68.1}) at our disposal, we further derive and simplify Gelfand-Shilov integrals for odd and even dimensions, respectively.

\subsection{The integral of $G_0(t,\boldsymbol{x}_0;\bx)$ and $\ddot{G}_{0}(t,\boldsymbol{x}_0;\bx)$ for $m$ even}
In order to compute the integral of $G_0(t,\bx;\bx_0)$, we apply it to a test function $\phi\in C_0^{\infty}(\mathbb{R}^m)$ for fixed $t$ as follows,
\begin{equation}\label{2.69}
\begin{aligned}
\int G_0(t,\boldsymbol{x}_0;\boldsymbol{x}) \phi(\boldsymbol{x}) \mathrm{d} \boldsymbol{x} &=\Int{}{} v_0(\bx)\,f_+^{-\frscr{m-1}{2}}[t^2-\tau^2(\bx)]
\,\phi(\bx)\,\drm\bx\\
& =\int\left(-\frac{\partial}{2 \tau \partial \tau}\right)^{\frac{m-2}{2}} f_{+}^{-\frac{1}{2}}\left[t^2-\tau^2(\boldsymbol{x})\right] v_0[\boldsymbol{x}(\tau, \boldsymbol{\omega})] \phi[\boldsymbol{x}(\tau, \boldsymbol{\omega})]\left|\frac{\partial S}{\partial \boldsymbol{\omega}}\right| c \mathrm{~d} \tau \mathrm{d} \boldsymbol{\omega}\\
&=\int\left(-\frac{\partial}{2 \tau \partial \tau}\right)^{\frac{m-2}{2}} f_{+}^{-\frac{1}{2}}\left[t^2-\tau^2(\boldsymbol{x})\right] \frac{\tau^{m-1}\phi}{4 \rho_0 c_0^m \pi^{m-1} \rho  v_0}  \mathrm{~d} \tau \mathrm{d} \boldsymbol{\omega}.
\end{aligned}
\end{equation}

When $m=2$, we have the following integral:
\begin{equation}\label{2.70}
  \int G_0(t,\bx_0;\boldsymbol{x}) \phi(\boldsymbol{x}) \mathrm{d} \boldsymbol{x}=\frac{1}{\left(-\frac{1}{2}\right) !} \int_{0}^{2\pi} \int_0^t\left[t^2-\tau^2\right]^{-\frac{1}{2}} \frac{\tau\phi}{4 \rho_0 c_0^2 \pi \rho v_0} \mathrm{d} \tau \mathrm{d} \theta,
\end{equation}
where the take-off angle $\boldsymbol{\omega}=\theta\in [0,2\pi]$.

To handle the integral of $\ddot{G}_0(t,\boldsymbol{x}_0;\bx)$, we use the fact that the Green's function satisfies the wave equation (\ref{2.1}). Therefore, for the given test function $\phi$ and $t>0$, we carry out integration by parts:
\begin{equation}\label{2.74}
  \begin{aligned}
\int \rho(\boldsymbol{x}) \ddot{G}_{0}\left(t, \boldsymbol{x}_0; \boldsymbol{x}\right) \phi(\boldsymbol{x}) d \boldsymbol{x} & =\int \nab \cdot\left(\nu \nab G_0\left(t, \boldsymbol{x}_0;\bx\right)\right) \phi(\boldsymbol{x}) d \boldsymbol{x} \\
& =-\int\left(\nu \nab G_0\left(t, \boldsymbol{x}_0; \boldsymbol{x}\right)\right) \cdot \nab(\phi(\boldsymbol{x})) d \boldsymbol{x} \\
& =\int G_0\left(t, \boldsymbol{x}_0; \boldsymbol{x}\right) \nab \cdot(\nu \nab \phi(\boldsymbol{x})) d \boldsymbol{x}.
\end{aligned}
\end{equation}

Further, when $m=2$, combining (\ref{2.74}) with (\ref{2.70}) yields:
\begin{equation}\label{2.75}
  \int \ddot{G}_0(t,\bx_0;\boldsymbol{x}) \phi(\boldsymbol{x}) d\boldsymbol{x}=\frac{1}{\left(-\frac{1}{2}\right) !} \int_{0}^{2\pi} \int_0^t\left[t^2-\tau^2\right]^{-\frac{1}{2}} \frac{\tau \nab \cdot\left(\nu\nab (\phi/\rho)\right)}{4 \rho_0 c_0^2 \pi \rho v_0} \mathrm{d} \tau \mathrm{d} \theta\; .
\end{equation}

\subsection{The integral of $\dot{G}_0(t,\boldsymbol{x}_0;\bx)$ for $m$ even}
Differentiating (\ref{2.69}) with respect to $t$ we obtain
\begin{equation}\label{2.71}
  \begin{aligned}
 \int \dot{G}_0(t,\bx_0; \boldsymbol{x}) \phi(\boldsymbol{x}) \mathrm{d} \boldsymbol{x}&=2 t \int v_0(\boldsymbol{x}) f_{+}^{-\frac{m-1}{2}-1}\left[t^2-\tau^2(\boldsymbol{x})\right] \phi(\boldsymbol{x}) \mathrm{d} \boldsymbol{x} \\
 &=2 t \int v_0(\boldsymbol{x}) f_{+}^{-\frac{m-2}{2}-\frac{1}{2}-1}\left[t^2-\tau^2(\boldsymbol{x})\right] \phi(\boldsymbol{x}) \mathrm{d} \boldsymbol{x} \\
&=2 t \int\left(-\frac{\partial}{2 \tau \partial \tau}\right)^{(k+1)} f_{+}^{-\frac{1}{2}}\left[t^2-\tau^2\right] v_0 \phi \left|\frac{\partial S}{\partial \boldsymbol{\omega}}\right| c \mathrm{d} \tau \mathrm{d} \boldsymbol{\omega}\\
&=2 t \int\left(-\frac{\partial}{2 \tau \partial \tau}\right)^{(k+1)} f_{+}^{-\frac{1}{2}}\left[t^2-\tau^2\right]  \frac{\tau^{m-1}\phi}{4 \rho_0 c_0^m \pi^{m-1} \rho  v_0}  \mathrm{d} \tau \mathrm{d} \boldsymbol{\omega},\\
\end{aligned}
\end{equation}
where $k=\frac{m-2}{2}$.

When $m=2$ such that $k=0$, we carry out integration by parts, yielding
\begin{equation}\label{2.72}
  \begin{aligned}
 \int \dot{G}_0(t,\bx_0;\boldsymbol{x}) \phi(\boldsymbol{x}) \mathrm{d} \boldsymbol{x}
&=2 t \int_{0}^{2\pi}\int_0^t \left(-\frac{\partial}{2 \tau \partial \tau}\right) f_{+}^{-\frac{1}{2}}\left[t^2-\tau^2\right]  \frac{\tau\phi}{4 \rho_0 c_0^2 \pi \rho  v_0}  \mathrm{d} \tau \mathrm{d} \theta\\
&=\frac{t}{\left(-\frac{1}{2}\right) !} \int_{0}^{2\pi}\left[\left[t^2-\tau^2\right]_+^{-\frac{1}{2}} \frac{\phi}{4 \rho_0 c_0^2 \pi \rho  v_0}\right]_{\tau=t}^0 \mathrm{~d} \theta \\
&\quad +\frac{t}{\left(-\frac{1}{2}\right) !} \int_{0}^{2\pi} \int_0^t\left[t^2-\tau^2\right]^{-\frac{1}{2}}\left(\frac{\partial}{\partial \tau}\right)\left[\frac{\phi}{4 \rho_0 c_0^2 \pi \rho  v_0}\right] \mathrm{d} \tau \mathrm{d} \theta\\
&=\frac{\phi(\bx_0)}{\rho(\bx_0)} +\frac{t}{\left(-\frac{1}{2}\right) !} \int_{0}^{2\pi} \int_0^t\left[t^2-\tau^2\right]^{-\frac{1}{2}}c^2\nab\tau\cdot \nab\left[\frac{\phi}{4 \rho_0 c_0^2 \pi \rho  v_0}\right] \mathrm{d} \tau \mathrm{d} \theta,\\
\end{aligned}
\end{equation}
where we have used the definition (\ref{2.6}), the initialization (\ref{2.41}) and the fact that $\boldsymbol{x}(\tau,\theta)|_{\tau=0}=\boldsymbol{x}_0$.
%%%%%%%%%%%%%%%%%%%%%%%%%%%%%%%%%%%%%%%%%%%%%%%%%%%%%%%%%%%%%%%%%%%%%%%%%%%%%%%%
\subsection{The integral of $G_0(t,\boldsymbol{x}_0;\bx)$ and $\ddot{G}_{0}(t,\boldsymbol{x}_0;\bx)$ for $m$ odd}

As before, we apply $G_0(t,\boldsymbol{x}_0;\bx)$ to a test function $\phi(\bx)$ for fixed $t$ as follows,
\beq{add4.8}\ba{rcl}
\Int{}{} G_0(t,\bx)\,\phi(\bx)\,\drm\bx
&=& \Int{}{} v_0(\bx)\,f_+^{-\frscr{m-1}{2}}[t^2-\tau^2(\bx)]
\,\phi(\bx)\,\drm\bx
\\ \\
&=& \Int{}{} v_0(\bx)\,\dlt^{(\frscr{m-3}{2})}[t^2-\tau^2(\bx)]
\,\phi(\bx)\,\drm\bx
\\ \\
&=& \Int{}{} \dlt^{(\frscr{m-3}{2})}[t^2-\tau^2]
v_0[\bx(\tau,\bom)]\,\phi[\bx(\tau,\bom)]\,\Big|\fr{\partial S}{\partial\bld\omega}\Big|\,c\,\drm \tau\,\drm\bld{\omega}
\\ \\
&=& \Int{}{} \dlt^{(\frscr{m-3}{2})}[t^2-\tau^2]\,\fr{\tau^{m-1}\phi\,}{4\rho_0c_0^m\pi^{m-1}\rho v_0}\,\drm \tau\,\drm\bld{\omega}.
\ea\eeq When $m=3$, we have the following integral,
\beq{SASWE4.02r}\ba{rcl}
\Int{}{} G_0(t,\bx_0;\bx)\,\phi(\bx)\,\drm\bx
&=& \Int{}{} \dlt[t^2-\tau^2]\,\fr{\tau^{2}\phi\,}{4\rho_0c_0^3\pi^{2}\rho v_0}\,\drm \tau\,\drm\bld{\omega}
\\ \\
&=& \Int{}{}\fr{\dlt[t-\tau] }{2t} \fr{\tau^{2}\phi}{4  \rho_0 c_0^3 \pi^{2} \rho  v_0}d \tau d\bld{\omega}\\ \\
&=& t\Int{}{} \left[\fr{\phi }{8  \rho_0 c_0^3 \pi^{2} \rho  v_0}\right]_{\tau=t} d\boldsymbol{\omega}.
\ea\eeq

Combining (\ref{2.74}) with \rf{SASWE4.02r}, we obtain the integral of $\ddot{G}_{0}(t,\boldsymbol{x}_0;\bx)$ for $m=3$:
\beq{SASWE4.06r3}\ba{rcl}
\Int{}{} \ddot{G}_0(t,\bx_0;\bx)\,\phi(\bx)\,\drm\bx
&=&t\Int{}{} \left[\fr{\nab \cdot(\nu\nab\left(\phi/\rho)\right)}{8\rho_0c_0^3\pi^2\rho v_0 }\right]_{\tau=t}d\boldsymbol{\omega}.
\ea\eeq
%%%%%%%%%%%%%%%%%%%%%%%%%%%%%%%%%%%%%%%%%%%%%%%%%%%%%%%%%%%%%%%%

\subsection{The integral of $\dot{G}_0(t,\boldsymbol{x}_0;\boldsymbol{x})$ for $m$ odd}
We first introduce the following properties of the $\delta$-function: For $t, \tau\geq 0$,
\begin{equation}
     \label{eq:prop:delta:1}
        \delta^{(1)}(t^2-\tau^2)=
        \frac{1}{4\tau^3}(\delta(t-\tau)+\tau\delta^{(1)}(t-\tau)).\\
\end{equation}
Differentiating \eqref{SASWE4.02r} with respect to $t$ and letting $m=3$, we obtain
\beq{SASWE4.06r}\ba{rcl}
\Int{}{} \dot{G}_0(t,\bx_0;\bx)\,\phi(\bx)\,\drm\bx
&=& \Int{}{} 2t \dlt^{(1)}[t^2-\tau^2]
v_0[\bx(\tau,\bld{\omega})]\,\phi[\bx(\tau,\bld{\omega})]\,\Big|\fr{\partial S}{\partial\bld\omega}\Big|c[\bx(\tau,\bld{\omega})]\,\drm \tau\,\drm\bld{\omega}
\\ \\
&=& \Int{}{} \dlt[t-\tau]
v_0[\bx(\tau,\bld{\omega})]\,\phi[\bx(\tau,\bld{\omega})]\,\Big|\fr{\partial S}{\partial\bld\omega}\Big|\fr{t\,c[\bx(\tau,\bld{\omega})]}{2\,\tau^3}\,\drm \tau\,\drm\bld{\omega}
\\ \\
& & + \Int{}{} \dlt^{(1)}[t-\tau]
v_0[\bx(\tau,\bld{\omega})]\,\phi[\bx(\tau,\bld{\omega})]\,\Big|\fr{\partial S}{\partial\bld\omega}\Big|\,\fr{t\,c[\bx(\tau,\bld{\omega})]}{2\,\tau^2}\,\drm \tau\,\drm\bld{\omega}
\\ \\
&=& \Int{}{} \dlt[t-\tau]
\,\psi_2(\tau,\bld{\omega})\,\drm \tau\drm\bld{\omega}+ \Int{}{} \fr{d}{d\tau}\left[
v_0[\bx(\tau,\bld{\omega})]\,\phi[\bx(\tau,\bld{\omega})]\,\Big|\fr{\partial S}{\partial\bld\omega}\Big|\,\fr{t\,c[\bx(\tau,\bld{\omega})]}{2\,\tau^2}\right]_{\tau=t}\,\drm\bld{\omega}
\\ \\
&=& \,\Int{}{}\psi_2(t,\bld{\omega})\,\drm\bld{\omega} + \Int{}{} \fr{d}{d\tau}\left[
v_0[\bx(\tau,\bld{\omega})]\,\phi[\bx(\tau,\bld{\omega})]\,\Big|\fr{\partial S}{\partial\bld\omega}\Big|\,\fr{t\,c[\bx(\tau,\bld{\omega})]}{2\,\tau^2}\right]_{\tau=t}\,\drm\bld{\omega} \;
\\ \\
&=& \,\Int{}{}\psi_2(t,\bld{\omega})\,\drm\bld{\omega} + \Int{}{} \fr{d}{d\tau}\left[\tau \psi_2(\tau,\bld{\omega})\right]_{\tau=t}\, \drm\bld{\omega} \;,
\ea\eeq
where
\beq{SASWE4.07r}
\psi_2(\tau,\bld{\omega})
\,=\,v_0[\bx(\tau,\bld{\omega})]\,\phi[\bx(\tau,\bld{\omega})]\Big|\fr{\partial S}{\partial\bld\omega}\Big|\,\fr{t\,c[\bx(\tau,\bld{\omega})]}{2\,\tau^3}= \frac{t\phi}{8 \rho_0 c_0^3 \pi^{2} \rho  v_0\tau}  \;.
\eeq
That is,
\beq{SASWE4.06r2}\ba{rcl}
\Int{}{} \dot{G}_0(t,\bx_0;\bx)\,\phi(\bx)\,\drm\bx
&=&\Int{}{} \left[\fr{\phi}{8\rho_0c_0^3\pi^2\rho v_0 }\right]_{\tau=t}d\boldsymbol{\omega}+t\Int{}{} \left[c^2\nab \tau \cdot \nab \left(\fr{\phi}{8\rho_0c_0^3\pi^2\rho v_0} \right)\right]_{\tau=t}d\bld{\omega}.
\ea\eeq
%%%%%%%%%%%%%%%%%%%%%%%%%%%%%%%%%%%%%%%%%%%%%%%%%%%%%%%%%%%%%%%%%%%%%%%%%%%%%%%%%
\subsection{HKH propagator}
By choosing the initial conditions as the test function $\phi$, we obtain the HKH propagator.

When $m=2$, by taking the test function $\phi(\bx)$ as $\rho(\bx) u_t(0,\bx)$ and $\rho(\bx) u(0,\bx) $ in \rf{2.70} and \rf{2.72}, respectively, we obtain
\begin{equation}\label{hkh2.1}
\begin{aligned}
       u(t,\boldsymbol{x}_0)&=u(0,\boldsymbol{x}_0)+\frac{1}{\left(-\frac{1}{2}\right) !} \int_{0}^{2\pi} \int_0^t\left[t^2-\tau^2\right]^{-\frac{1}{2}} \frac{\tau u_t(0,\boldsymbol{x}(\tau,\theta))}{4 \rho_0 c_0^2 \pi v_0} \mathrm{d} \tau \mathrm{d} \theta \\
       &\quad+\frac{t}{\left(-\frac{1}{2}\right) !} \int_{0}^{2\pi} \int_0^t\left[t^2-\tau^2\right]^{-\frac{1}{2}}c^2\nab \tau \cdot \nab \left[\frac{ u(0,\boldsymbol{x}(\tau,\theta))}{4 \rho_0 c_0^2 \pi  v_0} \right] \mathrm{d} \tau\mathrm{d} \theta.\\
       \end{aligned}
\end{equation}
Further taking $\phi(\bx)$ as $\rho(\bx) u_t(0,\bx)$ and $\rho(\bx) u(0,\bx)$ in \rf{2.72} and \rf{2.75}, respectively, yields
\begin{equation}\label{hkh2.2}
    \begin{aligned}
      u_t(t,\boldsymbol{x}_0)
       &=u_t(0,\boldsymbol{x}_0)+\frac{1}{\left(-\frac{1}{2}\right) !} \int_{0}^{2\pi} \int_0^t\left[t^2-\tau^2\right]^{-\frac{1}{2}} \frac{\tau \nab\cdot(\nu\nab u(0,\boldsymbol{x}(\tau,\theta)))}{4 \rho_0 \rho c_0^2 \pi v_0} \mathrm{d} \tau \mathrm{d} \theta\\
       &\quad+\frac{t}{\left(-\frac{1}{2}\right) !} \int_{0}^{2\pi} \int_0^t\left[t^2-\tau^2\right]^{-\frac{1}{2}}c^2\nab \tau \cdot \nab \left[\frac{ u_t(0,\boldsymbol{x}(\tau,\theta))}{4 \rho_0 c_0^2 \pi  v_0} \right] \mathrm{d} \tau \mathrm{d} \theta.\\
    \end{aligned}
\end{equation}

When $m=3$, by taking $\phi(\bx)$ as $\rho(\bx) u_t(0,\bx)$ and
$\rho(\bx) u(0,\bx) $ in \rf{SASWE4.02r} and \rf{SASWE4.06r2}, respectively, we obtain
\begin{equation}\label{hkh3.1}
\begin{aligned}
    u(t,\bx_0)=&\Int{}{} t\left[\fr{u_t(0,\bx(\tau,\boldsymbol{\omega})) }{8  \rho_0 c_0^3 \pi^{2}  v_0}\right]_{\tau=t} +\left[\fr{u(0,\bx(\tau,\boldsymbol{\omega}))}{8\rho_0c_0^3\pi^2 v_0 }\right]_{\tau=t}+t\left[c^2\nab \tau \cdot \nab \left(\fr{u(0,\bx(\tau,\boldsymbol{\omega}))}{8\rho_0c_0^3\pi^2 v_0} \right)\right]_{\tau=t}d\bld{\omega}.
\end{aligned}
\end{equation}
Further taking $\phi(\bx)$ as $\rho(\bx) u_t(0,\bx)$ and $\rho(\bx) u(0,\bx) $ in \rf{SASWE4.06r2} and \rf{SASWE4.06r3}, respectively, yields
\begin{equation}\label{hkh3.2}
\begin{aligned}
    u_t(t,\bx_0)=&\Int{}{} t\left[\fr{\nab \cdot (\nu\nab u(0,\bx(\tau,\boldsymbol{\omega}))) }{8  \rho_0 \rho c_0^3 \pi^{2}  v_0}\right]_{\tau=t} +\left[\fr{u_t(0,\bx(\tau,\boldsymbol{\omega}))}{8\rho_0c_0^3\pi^2 v_0 }\right]_{\tau=t}+t\left[c^2\nab \tau \cdot \nab \left(\fr{u_t(0,\bx(\tau,\boldsymbol{\omega}))}{8\rho_0c_0^3\pi^2 v_0} \right)\right]_{\tau=t}d\bld{\omega}.
\end{aligned}
\end{equation}

Formulas \rf{hkh2.1}-\rf{hkh2.2} and \rf{hkh3.1}-\rf{hkh3.2} are two-dimensional and three-dimensional HKH propagators, respectively, which are used to propagate the wavefield from $\tau=0$ to $\tau=t$, where $0<t<\bar{T}$. Although HKH propagators are only valid for a
short-time period in a caustic-free neighborhood, recursively applying this propagator in time yields the {\it Hadamard integrator} to solve time-dependent wave equations globally in time, where caustics are treated implicitly and spatially overturning waves are handled naturally.
%%%%%%%%%%%%%%%%%%%%%%%%%%%%%%%%%%%%%%%%%%%%%%%%%%%%%%%%%%%%%%%%%%%%%%%%%%%%%%%%

\section{Numerics for Hadamard integrators}
Here, we present numerics for the Hadamard integrator. Essentially, we will numerically discretize the HKH propagator and obtain the Hadamard integrator by recursively applying the propagator in time. For highly  oscillatory wavefields, we maintain a fixed number of points per wavelength (PPW) to uniformly discretize the computational domain into regular grid points.

To begin with, we briefly discuss high-order numerical schemes for the eikonal and transport equations, yielding the squared-phase function  $\tau^2$ and the Hadamard coefficient $v_0$, respectively. Subsequently, we present Gaussian quadrature formulas for the HKH propagator. To ensure accurate computation of the numerical quadrature, we utilize the ray tracing method to obtain locations of the equal-time wavefront. We then construct cubic spline interpolations to approximate the integrands at the wavefront locations, using the oscillatory initial data and computed Hadamard ingredients on the uniform grid, respectively. Finally, we present a preliminary algorithm of the Hadamard integrator for time-dependent wave equations.

\subsection{Numerical schemes for eikonal and transport equations}
The leading term of the Hadamard's ansatz is defined by two functions, the eikonal $\tau$ satisfying the eikonal equation (\ref{SASWE2.10}) and the Hadamard coefficient satisfying the transport equation (\ref{SASWE2.12}). Since we have assumed that the Hadamard's ansatz is valid locally around the source point, we need access to these two functions in order to construct the ansatz. Since the eikonal equation as a first-order nonlinear partial differential equation does not have analytical solutions in general, we have to use a robust, high-order numerical scheme to compute it; moreover, the eikonal equation equipped with a point-source condition is even more tricky to deal with due to the upwind singularity at the source point \cite{qiasym02adapt}. To make the situation even more complicated, the transport equation (\ref{SASWE2.12}) for the Hadamard coefficient $v_0$ is weakly coupled with the eikonal equation (\ref{SASWE2.10}) in that the coefficients of the former equation depend on the solution of the latter.

Fortunately, this set of weakly coupled equations with point-source conditions has been solved to high-order accuracy by using Lax-Friedrichs weighted essentially non-oscillatory (LxF-WENO) sweeping schemes as demonstrated in \cite{qiayualiuluobur16}. The high-order schemes in \cite{qiayualiuluobur16} have adopted essential ideas from many sources including  \cite{oshshu91,jiapen00,qiasym02adapt,kaooshqia04,zha05,zharechov05,fomluozha09,luoqiabur14b,luqiabur16} and have been used in many applications. Consequently, we will adopt these schemes to our setting as well and we omit details here.

\subsection{Numerical quadrature for the HKH propagator}
Starting from the Gelfand-Shilov integrals in the Lagrangian formulation, we apply the Gaussian quadrature to numerically implement the HKH propagator in geodesic polar coordinates.
\subsubsection{When $m=2$}
We set $\tau=t\cos\zeta$ in \rf{hkh2.1} and \rf{hkh2.2} to obtain
\begin{equation}\label{ray1.1}
\begin{aligned}
    u(t,\boldsymbol{x}_0)
      &= u(0,\boldsymbol{x}_0)+\frac{t}{\sqrt{\pi} } \int_{0}^{2\pi} \int_0^{\frac{\pi}{2}} \frac{\cos\zeta u_t(0,\boldsymbol{x}(t\cos\zeta,\theta))}{4 \rho_0 c_0^2 \pi v_0}+c^2\nab \tau \cdot \nab \left[\frac{ u(0,\boldsymbol{x}(t\cos\zeta,\theta))}{4 \rho_0 c_0^2 \pi  v_0} \right] \mathrm{d} \zeta \mathrm{d} \theta\\
        &\equiv u(0,\boldsymbol{x}_0)+\frac{t}{\sqrt{\pi} } \int_{0}^{2\pi} \int_0^{\frac{\pi}{2}} F_1(\zeta,\theta)+F_2(\zeta,\theta) \mathrm{d} \zeta \mathrm{d} \theta,
\end{aligned}
\end{equation}
where
\begin{equation}
    F_1(\zeta,\theta)=\frac{\cos\zeta u_t(0,\boldsymbol{x}(t\cos\zeta,\theta))}{4 \rho_0 c_0^2 \pi v_0},\quad F_2(\zeta,\theta)=c^2\nab \tau \cdot \nab \left[\frac{ u(0,\boldsymbol{x}(t\cos\zeta,\theta))}{4 \rho_0 c_0^2 \pi  v_0} \right],
\end{equation}
and
\begin{equation}\label{ray1.2}
\begin{aligned}
    u_t(t,\boldsymbol{x}_0)&=u_t(0,\boldsymbol{x}_0)+\frac{t}{\sqrt{\pi}} \int_{0}^{2\pi} \int_0^{\frac{\pi}{2}} \frac{\cos\zeta \nab\cdot(\nu\nab u(0,\boldsymbol{x}(t\cos\zeta,\theta)))}{4 \rho_0 \rho c_0^2 \pi v_0} +c^2\nab \tau \cdot \nab \left[\frac{ u_t(0,\boldsymbol{x}(t\cos\zeta,\theta))}{4 \rho_0 c_0^2 \pi  v_0} \right]\mathrm{d}\zeta \mathrm{d} \theta\\
        &\equiv u_t(0,\boldsymbol{x}_0)+\frac{t}{\sqrt{\pi} } \int_{0}^{2\pi} \int_0^{\frac{\pi}{2}} F_3(\zeta,\theta)+F_4(\zeta,\theta) \mathrm{d} \zeta \mathrm{d} \theta,
\end{aligned}
\end{equation}
where
\begin{equation}
    F_3(\zeta,\theta)=\frac{\cos\zeta \nab\cdot(\nu\nab u(0,\boldsymbol{x}(t\cos\zeta,\theta)))}{4 \rho_0 \rho c_0^2 \pi v_0},\quad F_4(\zeta,\theta)=c^2\nab \tau \cdot \nab \left[\frac{ u_t(0,\boldsymbol{x}(t\cos\zeta,\theta))}{4 \rho_0 c_0^2 \pi  v_0} \right].
\end{equation}

We uniformly discretize $[0,\frac{\pi}{2}]\times [0,2\pi]$ into a grid of size $M_1\times M_2,$ which corresponds to using the tensorized  Gaussian-Chebyshev quadrature nodes with respect to $(\tau,\theta)$. Then, we obtain
\begin{equation}\label{ray4.13}
    u(t,\boldsymbol{x}_0)\approx u(0,\bx_0)+\frac{t}{\sqrt{\pi}}\sum_{i=1}^{M_1}\sum_{j=1}^{M_2} (F_1(\zeta_i,\theta_j)+ F_2(\zeta_i,\theta_j))\frac{\pi}{2M_1}\frac{2\pi}{M_2},
\end{equation}
and
\begin{equation}\label{ray4.14}
    u_t(t,\boldsymbol{x}_0)\approx u_t(0,\bx_0)+\frac{t}{\sqrt{\pi}}\sum_{i=1}^{M_1}\sum_{j=1}^{M_2} (F_3(\zeta_i,\theta_j)+ F_4(\zeta_i,\theta_j))\frac{\pi}{2M_1}\frac{2\pi}{M_2}.
\end{equation}
%For each source $\boldsymbol{x}_0$, we trace the rays via the ODE system (\ref{ode2}) to obtain the wavefront locations:
%\begin{equation}
%    \bx(i,j)=P^{-1}\left[(t\cos\zeta_i,\theta_j)\right], i=1,2,\cdots,M_1,j=1,2,\cdots,M_2,
%\end{equation}
%and the arrival angles $\nab \tau(i,j)=\nab \tau(\bx(i,j)).$ Next, we solve eikonal equaiton and transport equation with source $\bx_0$ then interpolate $F_1$, $F_2$, $F_3$, $F_4$ onto ${\bx(i,j)}$ and perform the summations\rf{ray4.13} and \rf{ray4.14} to compute the HKH propagator.

\subsubsection{When $m=3$}
We rewrite \rf{hkh3.1} and \rf{hkh3.2} as follows,
\begin{equation}
\begin{aligned}
    u(t,\bx_0)
    \equiv &\Int{}{}F_5(\boldsymbol{\omega})+F_6(\boldsymbol{\omega})\ d\boldsymbol{\omega},
\end{aligned}
\end{equation}
where
\begin{equation}
    F_5(\boldsymbol{\omega})=\left[\fr{t u_t(0,\bx(\tau,\boldsymbol{\omega}))+ u(0,\bx(\tau,\boldsymbol{\omega}))}{8  \rho_0 c_0^3\pi^{2} v_0}\right]_{\tau=t},\quad  F_6(\boldsymbol{\omega})=t\left[c^2\nab \tau \cdot \nab \left(\fr{u(0,\bx(\tau,\boldsymbol{\omega}))}{8\rho_0c_0^3\pi^2 v_0} \right)\right]_{\tau=t},
\end{equation}
and
\begin{equation}
\begin{aligned}
    u_t(t,\bx_0)\equiv &\Int{}{}F_7(\boldsymbol{\omega})+F_8(\boldsymbol{\omega})\ d\boldsymbol{\omega},
\end{aligned}
\end{equation}
where
\begin{equation}
    F_7(\boldsymbol{\omega})=\left[\fr{t \nab\cdot(\nu \nab u(0,\bx(\tau,\boldsymbol{\omega})))+ \rho u_t(0,\bx(\tau,\boldsymbol{\omega}))}{8 \rho \rho_0 c_0^3\pi^{2} v_0}\right]_{\tau=t},\quad  F_8(\boldsymbol{\omega})=t\left[c^2\nab \tau \cdot \nab \left(\fr{u_t(0,\bx(\tau,\boldsymbol{\omega}))}{8\rho_0c_0^3\pi^2 v_0} \right)\right]_{\tau=t}.
\end{equation}

Now we utilize the tensorized Gauss-Legendre quadrature with respect to $\boldsymbol{\omega}=(\theta,\xi)\in[0,2\pi]\times [0,\pi]$  to evaluate the spherical integrals:
\begin{equation}\label{ray4.22}
    u(t,\boldsymbol{x}_0)\approx \frac{\pi}{M_3}\sum_{i=1}^{2M_3}\sum_{j=1}^{M_3} w_j\left(F_5(\theta_i,\xi_j)+ F_6(\theta_i,\xi_j)\right),
\end{equation}
\begin{equation}\label{ray4.23}
    u_t(t,\boldsymbol{x}_0) \approx \frac{\pi}{M_3}\sum_{i=1}^{2M_3}\sum_{j=1}^{M_3} w_j\left(F_7(\theta_i,\xi_j)+ F_8(\theta_i,\xi_j)\right),
\end{equation}
where we set $\{\xi_j\}$ so that $\{\cos\xi_j\}$ and $w_j$ are nodes and weights of the $M_3$-order Gauss-Legendre quadrature in $[-1,1]$, and $\{\theta_i\}$ are the equidistant nodes on $[0,2\pi]$.

%For each source $\boldsymbol{x}_0$, we trace the rays via the ODE system (\ref{ode2}) to obtain the wavefront locations:
%\begin{equation}
%    \bx(i,j)=P^{-1}\left[(\theta_i,\xi_j)\right], i=1,2,\cdots,2M_3, j=1,2,\cdots,M_3,
%\end{equation}
%and the arrival angles $\nab\tau(i,j)=\nab \tau(\bx(i,j))$. Next, we solve eikonal equaiton and transport equation with source $\bx_0$ then interpolate $F_5$, $F_6$, $F_7$, $F_8$ onto ${\bx(i,j)}$ and perform the summation \rf{ray4.22} and \rf{ray4.23} to compute the HKH propagator.

%For highly oscillatory initial data, we use local spline interpolation to calculate the values on the wavefront. For the wavefront locations and Hadamard ingredients, the low-rank representations are introduced in the next section.

\subsection{Hadamard ingredients}
In order to construct the HKH propagator, we need the following ingredients,
\begin{equation}\label{2.76}
  \tau^2,\; v_0,\; \nab\tau^2,\;  \nab v_0,\; \nab \tau^2 \cdot \nab v_0,
\end{equation}
which will be referred to as the Hadamard ingredients in the following. Here we use $\tau^2$ rather than $\tau$ since $\tau^2$ is differentiable at the point source while $\tau$ itself is not.

In principle, all these ingredients can be obtained from the eikonal $\tau$ and $v_0$ by solving the eikonal and transport equations with the LxF-WENO schemes \cite{qiayualiuluobur16}, respectively; however, since the LxF-WENO schemes yield solutions on uniform regular meshes, they are not completely fit with the geodesic polar coordinates used for evaluating Gelfand-Shilov integrals. Therefore, in order to compute these ingredients, we will develop a hybrid computational framework by combining the results from both the Eulerian LxF-WENO solver and the Lagrangian ray-tracing method.

The Lagrangian ray-tracing method, introduced in the next section, provides us with the arrival-angle related slowness vector $\nabla \tau$ as a by-product, and we have $\tau$ known in the geodesic polar coordinates. Consequently, we only need to solve the eikonal and transport equations to obtain the following ingredients,
\begin{equation}\label{2.76.1}
    v_0,\; \nabla v_0, \; \nab \tau^2\cdot \nab v_0.
\end{equation}
To obtain a third-order accurate $v_0$, we use the third-order WENO approximations \cite{jiapen00} to compute $\nabla v_0$. However, to avoid evaluating $\nabla \tau^2 \cdot \nabla v_0$ through numerical differentiation of $v_0$, we employ equation (\ref{SASWE2.12}) to obtain the following expression,
\begin{equation}\label{5.23.1}
    \nab \tau^2\cdot \nab v_0=-\frac{v_0[\nab \cdot (\nu \nab \tau^2)-2m\rho]}{2\nu}.
\end{equation}

It is worth noting that both $\nab \tau$ and $\nab\tau^2\cdot \nab v_0$  are used to compute $\nab \tau \cdot \nab v_0$ when needed. When $m$ is even, we combine $\nab v_0$ with $\nab \tau$ to form $\nab \tau \cdot \nab v_0$ directly, since using the expression
\begin{equation}\label{5.18.1}
    \nab \tau \cdot \nab v_0=\frac{1}{2\tau}\nab \tau^2\cdot \nab v_0
\end{equation}
involves dividing by $\tau$, which is unstable near the source, where $\tau=0$. When $m$ is odd, the Gelfand-Shilov integrals concentrate on the wavefront $\tau=\Delta t$ which is away from the source, and we can safely utilize $\nab \tau^2 \cdot \nab v_0$ in (\ref{5.18.1}) to evaluate $\nab \tau \cdot \nab v_0$ indirectly so as to reduce the number of needed ingredients and achieve better accuracy as well.

\subsection{Ray tracing methods}
We have derived the numerical quadrature formulas for the HKH propagator in geodesic polar coordinates. The crucial task now is to accurately compute these formulas. To achieve this, it is imperative to determine the wavefront locations corresponding to the Gaussian-quadrature nodes, as both the wavefields and Hadamard ingredients are given in Cartesian coordinates. Hence, to facilitate the evaluation of the propagator, we introduce the Lagrangian ray-tracing method to compute wavefront locations accurately and obtain geodesic polar coordinates accordingly.

As mentioned earlier, in a caustic-free neighborhood of source $\boldsymbol{x}_0$, the geodesic polar transformation $P[\bx_0]$ is well-defined and one-to-one. So we define the inverse geodesic polar transformation for a fixed $\bx_0$ as follows,
\begin{equation}
    P^{-1}[\bx_0]: (\tau,\boldsymbol{\omega})\rightarrow \boldsymbol{x},
\end{equation}
which maps the traveltime and take-off angle to the corresponding physical location of the ray. Essentially, it traces the ray starting from $\bx_0$ with the take-off angle $\boldsymbol{\omega}$ until time $\tau$, which aligns very well with the method of characteristics for a short period of time $\tau$.

Solving the eikonal equation via the method of characteristics with the Hamiltonian $H(\boldsymbol{x,p})=c(\boldsymbol{x})|\boldsymbol{p}|$, where $\boldsymbol{p}=\nab\tau$, we obtain the following Hamiltonian system,
\begin{equation}\label{ode1}
    \begin{aligned}
        &\frac{\mathrm{d} \boldsymbol{x}}{\mathrm{d} t}=\nabla_{\boldsymbol{p}} H(\boldsymbol{x}, \boldsymbol{p})=c(\boldsymbol{x})\frac{ \boldsymbol{p}}{|\boldsymbol{p}|}\;,\\
        &\frac{\mathrm{d} \boldsymbol{p}}{\mathrm{d} t}=-\nabla_{\boldsymbol{x}} H(\boldsymbol{x}, \boldsymbol{p})=\nab c(\boldsymbol{x}) |\boldsymbol{p}|\;,\\
        &\frac{\mathrm{d} \tau}{\mathrm{d} t}=\nabla_{\boldsymbol{p}} H(\boldsymbol{x}, \boldsymbol{p})\cdot \boldsymbol{p}=c(\boldsymbol{x})\frac{\bld{p}}{|\bld{p}|}\cdot \boldsymbol{p}=c(\boldsymbol{x})|\boldsymbol{p}|=1\;,
    \end{aligned}
\end{equation}
with initial conditions
\begin{equation}\label{ode1ini}
    \boldsymbol{x}(0)=\boldsymbol{x}_0,\; \boldsymbol{p}(0)=n(\bx_0)\boldsymbol{\omega},\; \tau(0)=0.
\end{equation}
By the initial conditions (\ref{ode1ini}) for the system \eqref{ode1}, we immediately get $\tau=t$ so that we can rewrite the Hamiltonian  system as
\begin{equation}\label{ode2}
    \begin{aligned}
        &\frac{\mathrm{d} \boldsymbol{x}}{\mathrm{d} \tau}=\nabla_{\boldsymbol{p}} H(\boldsymbol{x}, \boldsymbol{p})=c(\boldsymbol{x})\frac{ \boldsymbol{p}}{|\boldsymbol{p}|}\;,\\
        &\frac{\mathrm{d} \boldsymbol{p}}{\mathrm{d} \tau}=-\nabla_{\boldsymbol{x}} H(\boldsymbol{x}, \boldsymbol{p})=\nab c(\boldsymbol{x})|\boldsymbol{p}|
    \end{aligned}
\end{equation}
with initial conditions
\begin{equation}\label{ode2ini}
    \boldsymbol{x}(0)=\boldsymbol{x}_0,\; \boldsymbol{p}(0)=n(\bx_0)\boldsymbol{\omega}.
\end{equation}

Given the source location $\bx_0$ and the take-off angle $\boldsymbol{\omega}$, we solve the Hamiltonian system (\ref{ode2}) with the initial conditions (\ref{ode2ini}) using Runge-Kutta methods until $\tau=\Delta t$. This allows us to obtain the physical location of the wavefront (ray) $\boldsymbol{x}|_{\tau=\Delta t}$, which yields the inverse geodesic polar transformation $\boldsymbol{x}=P^{-1}[\boldsymbol{x}_0](\Delta t, \boldsymbol{\omega})$. Additionally, as a by-product, we obtain the arrival angle $\boldsymbol{p}|_{\tau=\Delta t}=\nabla \tau |_{\tau=\Delta t}$, which serves as one of the Hadamard ingredients.

After obtaining the wavefront locations through ray tracing, the next step is to compute the integrands $F_i$ at these irregular locations, which consist of Hadamard ingredients \rf{2.76} and highly oscillatory initial data, specifically:
\begin{equation}\label{inidata}
    u(0,\cdot),\; u_t(0,\cdot),\; \nab u(0,\cdot),\; \nab u_t(0,\cdot),\; \nab\cdot(\nu\nab u(0,\cdot)).
\end{equation}

To obtain Hadamard ingredients on the curved wavefront, we use two steps: in the first step, we solve the eikonal and transport equations in a caustic-free neighborhood of the source $\bx_0$ on a regular grid, where the grid size $h_1$ is independent of $h$ used for interpolating the initial data; in the second step, we interpolate these quantities defined on the regular grid onto the curved wavefront by using cubic splines, where the accuracy of the interpolation is ensured by high-order numerical schemes for Hadamard ingredients \cite{liusonburqia23}.

On the other hand, to obtain the oscillatory initial data \rf{inidata} on the curved wavefront, we use a two-step strategy: in the first step, we compute the required numerical derivatives on a regular grid to guarantee high-order accuracy of these quantities; in the second step, we use cubic splines to interpolate these initial data defined on the regular grid onto the curved wavefront. Since $u(0,\cdot)=u^1$ and $u_t(0,\cdot)=u^2$ are highly oscillatory $L^2$ functions, accurately interpolating these initial data in a {\it weak} sense requires a sufficiently fine grid. Numerically, according to the oscillation frequency of $u^1$ and $u^2$, we will choose an appropriate grid size $h$ to discretize the computational domain so as to ensure that the number of points per wavelength (PPW) is approximately $10$.

To numerically implement the cubic spline interpolations, we utilize the MATLAB function \textit{griddedInterpolant} while setting the 'option' to 'spline'. This provides an interpolant that performs cubic spline interpolation at the input locations based on gridded data.

Finally, we can update the wavefield at $\bx_0$ by computing the summations \rf{ray4.13} and \rf{ray4.14} for the 2-D case or \rf{ray4.22} and \rf{ray4.23} for the 3-D case, leading to the desired HKH propagator.

\subsection{Algorithm for time-dependent wave equations}
Recursively using the current data of the wavefield as `initial' data to apply the HKH propagator yields the Hadamard integrator which solves the Cauchy problem of  time-dependent wave equations. We have the following preliminary algorithm.\\
\textbf{Algorithm 1}
\begin{enumerate}
    \item  Uniformly discretize the computational domain into a wave-resolution-satisfying regular grid; choose an appropriate time step size $\Delta t<\bar{T}$ according to the medium and set an ending time $T_{end}$ to ensure that the waves do not reach the computational boundary; initialize $u(0,\cdot)=u^1$ and $u_t(0,\cdot)=u^2$ on the regular grid; set time $T=0$ and the loop variable $k=0$.
    \item For the current time step with $T=k \Delta t$:
    \begin{enumerate}
        \item compute the following numerical differentiations on the regular grid
    \begin{equation}\label{al1.1}
         \nab u(k\Delta t,\cdot),\quad \nab u_t(k \Delta t,\cdot),\quad  \mbox{and}\quad\nab\cdot(\nu\nabla u(k\Delta t,\cdot)),
    \end{equation}
    which, along with $u(k\Delta t,\cdot)$ and $u_t(k\Delta t,\cdot)$, constitute the current data at the current time step;
     \item construct the cubic spline interpolations for the gridded current data using \textit{griddedInterpolant};
     %interpolate the current data \rf{al1.1} onto the predetermined regular grid via cubic splines;
     \item determine a subregion $\Omega_k$ of $\Omega$ that contains the region of influence of current data by extending each direction of the non-zero region of the current data outward by $\max c \Delta t$;
     \item  for $\bx_0\in \Omega_k$:     \begin{enumerate}
    \item solve the Hamiltonian system (\ref{ode2}) with the initial conditions (\ref{ode2ini}) to obtain wavefront locations and arrival angles; solve the eikonal and transport equations in a caustic-free neighborhood of the source $\bx_0$ and interpolate the Hadamard ingredients onto the wavefront locations  via cubic splines; interpolate the current data onto wavefront locations and generate the integrands $F_i$;
    \item use formulas \eqref{ray4.13} and \eqref{ray4.14} for 2-D or \eqref{ray4.22} and \eqref{ray4.23} for 3-D to update wavefields:   $$u=\left((k+1)\Delta t,\bx_0\right), \quad u_t\left((k+1)\Delta t,\bx_0\right);$$
 %   \item Repeat Step 3 to Step 4 for all $\Omega_k^{\ell}$.
    \end{enumerate}
%    \item Update $T=(k+1)\Delta t.$
    \item  set $u\left((k+1)\Delta t,\Omega\backslash \Omega_k\right)=u_t\left((k+1)\Delta t,\Omega\backslash \Omega_k\right)=0$.
    \end{enumerate}
    \item Update $T=(k+1)\Delta t.$ If $T<T_{end}$, then $k\leftarrow k+1$ and go to Step 2; else, stop.
\end{enumerate}

\iffalse
\textbf{Algorithm 1.}
\begin{enumerate}
    \item For the current time step $T=k \Delta t$, compute the following numerical derivatives
    \begin{equation}\label{al1.1}
         \nab u(k\Delta t,\cdot),\quad \nab u_t(k \Delta t,\cdot),\quad  \mbox{and}\quad\nab\cdot(\nu\nabla u(k\Delta t,\cdot)),
    \end{equation}
    and determine the region $\Omega_k$ to be updated in this step using the region of influence by extending each direction of the non-zero region of the wavefield outward by $\max c \Delta t$.  Interpolate the initial data \rf{al1.1} onto a predetermined wave-resolution-satisfying regular grid via cubic splines.

    \item For $\bx_0\in\Omega_k$: solve the Hamiltonian system (\ref{ode2}) with the initial conditions (\ref{ode2ini}) to obtain
    wavefront locations, locations of the Gaussian-quadrature nodes on the wavefront, and arrival angles.
    \item Solve the eikonal and transport equations
    in a caustic-free neighborhood of the source $\bx_0$ and interpolate the Hadamard ingredients onto the wavefront locations  via cubic splines.
    Then preform interpolation of the initial data on the wavefront locations to generate the integrands $F_i$.
    \item Use formulas \eqref{ray4.13} and \eqref{ray4.14} for 2-D or \eqref{ray4.22} and \eqref{ray4.23} for 3-D to update wavefields  $u(T+\Delta t,\bx_0)$ and $u_t(T+\Delta t,\bx_0)$.
    \item Repeat Step 2 to Step 4 for all $\bx_0\in \Omega_k$. Update $T=(k+1)\Delta t$.
    \item Repeat Step 1 to Step 5 until $T=T_{end}$.
    \end{enumerate}
    \fi
In this preliminary algorithm, for every $\bx_0$ we need to trace rays, solve the eikonal and transport equations, and perform corresponding interpolations, and these operations are expensive since the set of $\{\bx_0\}$ occupies a volume. Therefore, to accelerate this algorithm, we will construct low-rank representations of wavefront locations and Hadamard ingredients so that those Gelfand-Shilov integrals can be evaluated rapidly, where such low-rank representations result in an algorithm which is amenable to fast block-matrix operations.

\section{Fast computation of Hadamard integrators}
In the preliminary Algorithm 1, for each given source $\bx_0$ we trace rays, solve eikonal and transport equations, and interpolate the current data and Hadamard ingredients onto current wavefronts at every time step. Notably, except for the interpolation of the oscillatory current data, tracing rays and solving eikonal and transport equations only depend on the given medium so that wavefront locations and Hadamard ingredients can be precomputed and reused for different initial conditions. Moreover, assuming that the medium parameters $\rho$ and $\nu$ are analytic, we can construct low-rank representations of wavefront locations and Hadamard ingredients. To achieve this, we first introduce generic  multivariate Chebyshev interpolations.

\subsection{Multivariate Chebyshev interpolation}
Let us consider a function
\begin{equation}
    f(\boldsymbol{\eta}),\quad \boldsymbol{\eta}=[\eta_1,\eta_2,\cdots,\eta_M] \in [-1,1]^{M},
\end{equation}
which permits a low-rank representation. We can then expand $f$ in terms of Chebyshev polynomials of the first kind
\begin{equation}\label{4.3}
    f(\boldsymbol{\eta})\approx \sum_{i_1=1}^{n_1}\cdots \sum_{i_M=1}^{n_M}C(i_1,\cdots,i_M)T_{i_1}(\eta_1)\cdots T_{i_M}(\eta_M),
\end{equation}
where $n_k$ is the order of Chebyshev interpolation with respect to $\eta_k$, $T_{i_k}(\eta_k)=\cos(i_k\arccos(\eta_k))$ is the Chebyshev polynomial of the first kind of order $i_k$, and the tensor $C$ contains the spectral coefficients to be determined.

To construct the multivariate low-rank representation (\ref{4.3}), we create an $M$-dimensional tensor $F$ that contains the function values of $f$ at the tensor-product Chebyshev nodes. These nodes are defined as the following,
\begin{equation}
    [\boldsymbol{\eta}^c(i_1,i_2,\cdots,i_M)]=\left[\cos\left(\frac{2i_1-1}{2n_1}\right),\cos\left(\frac{2i_2-1}{2n_2}\right),\cdots,\cos\left(\frac{2i_M-1}{2n_M}\right) \right]
\end{equation}
where $1\leq i_k\leq n_{k}$. Subsequently,
\begin{equation}
  F\equiv f(\boldsymbol{\eta}^c).
\end{equation}
After applying the fast cosine transform to each dimension of $F$, we generate the spectral-coefficient tensor $C$ \cite{boy01}. When we aim to interpolate $f$ onto a mesh of size $N_1\times N_2\times \cdots\times N_{M}$, the computational cost of the direct summation of (\ref{4.3}) is $O(\prod_{i=1}^{M}n_iN_i)$. However, the Orszag partial summation method \cite{boy01} can be introduced to greatly reduce the cost. In \cite{luoqiabur14a}, we have given an efficient Chebyshev summation method for $M=3$, and this summation can be directly generalized to the current generic M-dimensional case. The computational cost is then reduced to
\begin{equation}\label{cost}
    O\left(\sum_{i=1}^{M}\left[ \prod_{k=1}^{i}n_k\prod_{\ell=i}^{M}N_{\ell} \right]\right).
\end{equation}
Since $n_k$ is the order of Chebyshev interpolation which is much smaller than $N_k$, where the mesh size $N_k$ depends on the initial conditions, the partial summation does reduce the cost of the interpolation onto a regular mesh. However, when we consider to interpolate onto  irregular locations, the above partial summation trick no longer works. This motivates us to construct low-rank representations in a principled manner so that we can perform interpolation onto irregular locations rapidly.

\subsection{Low-rank representations of wavefront locations and Hadamard ingredients}
Given a source $\bx_0$ and Gaussian-quadrature nodes in geodesic polar coordinates, we need wavefront locations $\bx$ and Hadamard ingredients (\ref{2.76}) so as to compute Gelfand-Shilov integrals numerically. To facilitate these evaluations, we now construct the Chebyshev-polynomial based low-rank representations with respect to traveltime $\tau$, source point $\bx_0$, and take-off angle $\boldsymbol{\omega}$.

The first set of low-rank representations, $\Lambda_1$, is constructed  to provide representations for wavefront locations and arrival angles, effectively avoiding repeatedly ray-tracing at each source,
\begin{equation}
    \Lambda_1: (\tau,\bx_0, \boldsymbol{\omega}) \rightarrow \bx,\nab \tau(\bx_0;\bx),
\end{equation}
which is actually the low-rank representation of the short-time solver for the autonomous Hamiltonian system \rf{ode2} with analytic coefficients $c$ and $\nab c$.

The second set of low-rank representations, $\Lambda_2$, is constructed to provide representations for Hadamard ingredients on the wavefront locations $\bx$, effectively avoiding repeatedly solving the eikonal and transport equations at each point source,
\begin{equation}
    \Lambda_2: (\tau,\bx_0, \boldsymbol{\omega})\rightarrow v_0(\bx_0;\bx), \nab v_0(\bx_0;\bx), \nab \tau^2(\bx_0;\bx)\cdot \nab v_0(\bx_0;\bx),
\end{equation}
which additionally avoids interpolations on curved wavefronts. In fact, we have constructed Chebyshev interpolants for Hadamard ingredients with respect to source $\bx_0$ and wavefront location $\bx$ as is done in \cite{liusonburqia23}. Now since both $\Lambda_1$ and $\Lambda_2$ are low-rank representations with respect to $(\tau,\bx_0, \boldsymbol{\omega})$, they allow us to use block-matrix based partial summation to accelerate evaluations of numerical integrals.

When the underlying medium changes rapidly, we can divide the computational region into several sub-regions and construct the low-rank representations in each sub-region separately. Although this may slightly increase the cost of pre-computation, we may use lower  order Chebyshev interpolants in each sub-region so that the overall accuracy for the entire region can be improved and the construction time of interpolants may be reduced.

Different low-rank representations $\Lambda_1$ and $\Lambda_2$ are constructed for the 2-D and 3-D case, respectively. The main difference between these representations lies in whether the traveltime $\tau$ is treated as an interpolation variable (in the 2-D case) or a fixed parameter (in the 3-D case). This difference arises due to the distinct properties of Green's functions for wave equations in odd and even spatial dimensions.

\subsubsection{When $m=2$}
Given a source $\bx_0$ and time step $\Delta t$, the spatial support of the two-dimensional Green's function at $\bx_0$ is
$\{\bx:\tau(\bx_0;\bx)\leq \Delta t\}.$ Thus we construct the low-rank representations for wavefront locations $\bx=[x,y]$, Hadamard ingredients $v_0$, $\nab \tau=[n(\bx)\cos(\theta),n(\bx)\sin(\theta)]$ and $\nab v_0=[\frac{\partial v_0}{\partial x}, \frac{\partial v_0}{\partial y}]$ with respect to traveltime $\tau$, source $\bx_0=[x_0,y_0]$ and take-off angle $\theta_0$.

To construct $\Lambda_1$, we take
\begin{equation}
    f(\boldsymbol{\eta})=x,\;y,\;n(\bx)\cos\theta,\; n(\bx)\sin\theta,
\end{equation}
respectively, in (\ref{4.3}), where
\begin{equation}\label{4.8}
    \boldsymbol{\eta}=[\tau,x_0,y_0,\theta_0]\in[0,\Delta t]\times \Omega_c \times [0,2\pi]\equiv D_2,
\end{equation}
with $\Omega_c \subset \Omega$, and we can further map $D_2$ to $[-1,1]^4$ by translation and scaling.

To obtain $F=\{f(\tau_{i_1},x_{0,i_2},y_{0,i_3},\theta_{0,i_4})\}$,
we first use the arrival angle $\theta$ to parameterize the slowness vector $\boldsymbol{p}$ so as to reduce the Hamiltonian system (\ref{ode2}) into the following system \cite{qialeu04},
\begin{equation}\label{4.9}
    \frac{d x}{d\tau}=c\cos \theta,\quad \frac{d y}{d \tau}=c\sin\theta,\quad \frac{d \theta}{d \tau}=c_x\sin\theta-c_y\cos\theta\,;
\end{equation}
we then use the Runge-Kutta method(RK4) to solve the Hamiltonian system (\ref{4.9}) with initial conditions
\begin{equation}
    \bx(0)=\bx_0=[x_{0,i_2},y_{0,i_3}],\quad \boldsymbol{p}(0)=[n(\bx_0)\cos(\theta_{0,i_4}),n(\bx_0)\sin(\theta_{0,i_4})]
\end{equation}
until $\tau=\tau_{i_1}$; the resulting solutions yield $F$.

To construct $\Lambda_2$, we take
\begin{equation}
    f(\boldsymbol{\eta})=v_0,
    \;\frac{\partial v_0}{\partial x},\;\frac{\partial v_0}{\partial y},
\end{equation}
respectively, where we set $\boldsymbol{\eta}$ as in (\ref{4.8}). We follow the same procedure as before to solve the Hamiltonian system and obtain wavefront locations  $\bx(\tau_{i_1},x_{0,i_2},y_{0,i_3},\theta_{0,i_4})$. Additionally, we solve the eikonal and transport equations with the source located in the neighborhood defined by
\begin{equation}\label{4.11.1}
    [x_{0,i_2}-H,x_{0,i_2}+H]\times [y_{0,i_3}-H,y_{0,i_3}+H],
\end{equation}
where we carefully choose the value of $H$ to ensure that the neighborhood is caustic-free and includes the wavefronts required for our HKH propagator. After interpolating $v_0$, $\frac{\partial v_0}{\partial x}$, $\frac{\partial v_0}{\partial y}$ onto wavefront locations $\bx(\tau_{i_1},x_{0,i_2},y_{0,i_3},\theta_{0,i_4})$ using cubic splines, we obtain $F.$

Finally, we compute the spectral coefficients $C$ by applying the fast cosine transform to $F$.

\subsubsection{When $m=3$}
Given a source $\bx_0$ and time step $\Delta t$, the spatial support of the three-dimensional Green's function is $\{\bx:\tau(\bx_0;\bx)=\Delta t\}.$ Thus we take $\tau=\Delta t$ and construct the low-rank representations for wavefront locations $\bx=[x,y,z]$, Hadamard ingredients $\nab \tau=[n(\bx)\cos\theta\sin\xi,n(\bx)\sin\theta\sin\xi,n(\bx)\cos\xi]$, $v_0$ and $\nab \tau^2\cdot \nab v_0$ with respect to source $\bx_0=[x_0,y_0,z_0]$ and take-off angle $[\theta_0,\xi_0]$.

To construct $\Lambda_1$, we take
\begin{equation}
f(\boldsymbol{\eta})=x,\;y,\;z,\;n(\bx)\cos\theta\sin\xi,\;n(\bx)\sin\theta\sin\xi,\;n(\bx)\cos\xi,
\end{equation}
respectively, where
\begin{equation}\label{6.13.1}
    \boldsymbol{\eta}=[x_0,y_0,z_0,\theta_0,\xi_0]\in \Omega_c\times [0,2\pi]\times [0,\pi]\doteq D_3,
\end{equation}
with $\Omega_c \subset \Omega$, and we can further map $D_3$ to $[-1,1]^5$ by translation and scaling.

To obtain
$F=\{f(x_{0,i_1},y_{0,i_2},z_{0,i_3},\theta_{0,i_4},\xi_{0,i_5})\}$, we solve the Hamiltonian system (\ref{ode2}) by the Runge-Kutta method (RK4) equipped with the following initial conditions
\begin{eqnarray}
    \bx(0)&=&\bx_0=[x_{0,i_1},y_{0,i_2},z_{0,i_3}],\\
    \boldsymbol{p}(0)&=&[n(\bx_0)\cos\theta_{0,i_4}\sin\xi_{0,i_5},n(\bx_0)\sin\theta_{0,i_4}\sin\xi_{0,i_5},n(\bx_0)\cos\xi_{0,i_5}]
\end{eqnarray}
until $\tau=\Delta t$.

To construct $\Lambda_2$, we take
\begin{equation}
    f(\boldsymbol{\eta})=v_0,\; \nab \tau^2 \cdot \nab v_0,
\end{equation}
respectively, where we set $\boldsymbol{\eta}$ as (\ref{6.13.1}). Then we solve the eikonal and transport equations with the source $[x_{0,i_1},y_{0,i_2},z_{0,i_3}]$ located in the neighborhood defined by
\begin{equation}\label{6.16}
    [x_{0,i_1}-H,x_{0,i_1}+H]\times[y_{0,i_2}-H,y_{0,i_2}+H]\times[z_{0,i_3}-H,z_{0,i_3}+H],
\end{equation}
 Subsequently, we interpolate $v_0$ and $\nab \tau^2\cdot \nab v_0$ onto $\bx(x_{0,i_1},y_{0,i_2},z_{0,i_3},\theta_{0,i_4},\xi_{0,i_5})$ using cubic splines to obtain $F$. With $F$ at our disposal, we calculate the spectral coefficients $C$ via the fast cosine transform.

\subsection{Low-rank representation based fast algorithm.}
Based on the low-rank representations, we upgrade the preliminary Algorithm 1 into a blockwise fast algorithm.\\
\textbf{Algorithm 2}
\begin{enumerate}
    \item Construct the low-rank representations for wavefront locations and Hadamard ingredients according to Section 6.2.
    \item Uniformly discretize the computational domain into a wave-resolution-satisfying regular grid; choose an appropriate time step size $\Delta t<\bar{T}$ according to the medium and set an ending time $T_{end}$ to ensure that the waves do not reach the computational boundary; initialize $u(0,\cdot)=u^1$ and $u_t(0,\cdot)=u^2$ on the regular grid; set time $T=0$ and the loop variable $k=0$.
    \item For the current time step $T=k \Delta t$:
    \begin{enumerate}
     \item compute the following numerical differentiations on the regular grid
    \begin{equation}\label{al1.2}
         \nab u(k\Delta t,\cdot),\quad \nab u_t(k \Delta t,\cdot),\quad  \mbox{and}\quad\nab\cdot(\nu\nabla u(k\Delta t,\cdot)),
    \end{equation}
    which, along with $u(k\Delta t,\cdot)$ and $u_t(k\Delta t,\cdot)$, constitute the current data at the current time step;
     \item  construct the cubic spline interpolations for the gridded current data using \textit{griddedInterpolant};
    % interpolate the current data \rf{al1.2} onto the predetermined  regular grid via cubic splines;
     \item determine a block subregion $\Omega_k$ of $\Omega$ which contains the region of influence of current data by first extending each direction of the non-zero region of the current data outward by $\max c \Delta t$ and then finding a minimal block region that contains the current region, where the minimal block region is assigned to be $\Omega_k$;
     \item divide $\Omega_k$ into $N_k$ sub-domains $\Omega_k^{\ell}$ such that the size of each sub-domain is smaller than a predefined constant. For each sub-domain $\Omega_k^{\ell}$:     \begin{enumerate}
    \item use low-rank approximations to obtain the wavefront locations and the Hadamard ingredients on these locations; interpolate the current data onto wavefront locations and generate the integrands $F_i$;
    \item use formulas \eqref{ray4.13} and \eqref{ray4.14} for 2-D or \eqref{ray4.22} and \eqref{ray4.23} for 3-D to update wavefields:   $$u\left((k+1)\Delta t,\Omega_k^{\ell}\right),\quad u_t\left((k+1)\Delta t,\Omega_k^{\ell}\right);$$
 %   \item Repeat Step 3 to Step 4 for all $\Omega_k^{\ell}$.
    \end{enumerate}
    \item set $u\left((k+1)\Delta t,\Omega\backslash \Omega_k\right)=u_t\left((k+1)\Delta t,\Omega\backslash \Omega_k\right)=0$.
%    \item Update $T=(k+1)\Delta t.$
    \end{enumerate}
    \item Update $T=(k+1)\Delta t.$  If $T<T_{end}$, then $k\leftarrow k+1$ and go to Step 3; else, stop.
\end{enumerate}

We have in Step 3 carried out an additional domain decomposition, which is intended to control the memory usage for updating the oscillatory wavefields block by block. The rational is the following. Since the estimate of the computational cost (\ref{cost}) suggests that, the larger the size of the grid for multivariate Chebyshev interpolation, the more the computational cost is saved by using partial summation. Therefore, when computing resources are sufficient, we can skip this additional partition and update the wavefields on the entire $\Omega_k$ all at once to achieve the lowest computational cost. Numerically, we aim to choose larger subdomains $\Omega_k^{\ell}$ whenever possible, subject to the limitation of computing resources.
But when the available computing resources are limited, this additional partition might help us manage resources more effectively so that we can carry out large-scale computations.
In the 2-D case, we alternate directions and perform successive bisections to obtain sub-domains $\Omega_k^{\ell}$ of the same size. In the 3-D case, we partition $\Omega_k$ into layered $\Omega_k^{\ell}$ along a certain direction, such as the $z$-direction.

\section{Numerical examples}
This section provides numerical examples to demonstrate accuracy and performance of the proposed Hadamard integrator. Because an exact solution for the wave equation is not available in general, we numerically solve the wave equation with a pseudospectral method to obtain highly accurate numerical solutions, and these numerical solutions will serve as exact solutions to calibrate our Hadamard integrator. In addition, we use `RT' to indicate solutions computed by the Hadamard integrator.

\subsection{Two-dimensional examples}
\textbf{Example 1.} We set up the problem as the following.
\begin{itemize}
    \item[$\bullet$] $\rho=\frac{1}{(1+0.1\sin(2\pi x)\cos(2\pi y))^2},$\; $\nu=1,$\; and $c=1+0.1\sin(2\pi x)\cos(2\pi y)$.
    \item[$\bullet$] $u(0,x,y)=\sin \left( \pi \beta(x+y-1)\right) \exp \left(-600\left((x-0.5)^{2}+(y-0.5)^{2}\right)\right)$, and $u_t(0,x,y)=0$, where $\beta$ is a positive frequency parameter.
    \item[$\bullet$]The computational domain is $\Omega=[0,1]^2$ and the grid size used to discretize $\Omega$ is $h=\frac{1}{5\beta}.$
    %{\bf what is $h$ here?}
    \item[$\bullet$]The orders of the tensorized Chebyshev polynomials with respect to $[x_0,y_0,\tau,\theta_0]$ are $[15,15,15,15].$
    \item[$\bullet$]The numbers of Gaussian-quadrature nodes are $M_1=M_2=2\beta.$
\end{itemize}

Figure \ref{figure10} shows the velocity model, some rays and wavefronts, where no caustic transpires in the computational domain. Accordingly, we set $\Delta t=0.1$ for the HKH propagator. To construct the low-rank representations, we solve the Hamiltonian system (\ref{4.9}) using the Runge–Kutta method (RK4) for $1000$ time steps to obtain accurate wavefront locations, and we solve the eikonal and transport equations with grid size $h_1=0.004$ in the squared neighborhood (\ref{4.11.1}) with $H=0.2$.

     \begin{figure}[htbp]
     \centering
     \subfigure[]{
     \includegraphics[scale=0.385]{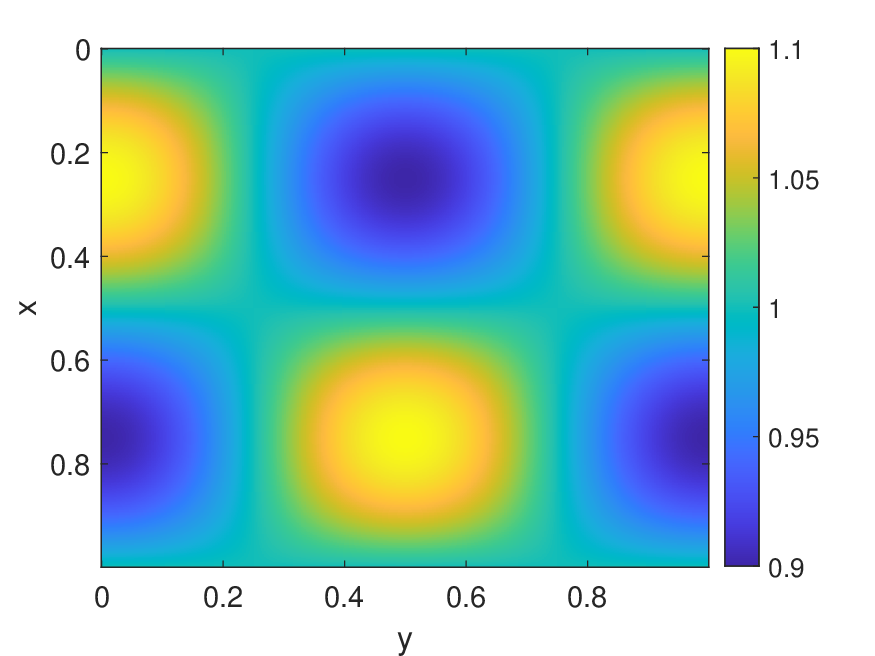}}
     \hspace{8mm}
     \subfigure[]{
     \includegraphics[scale=0.385]{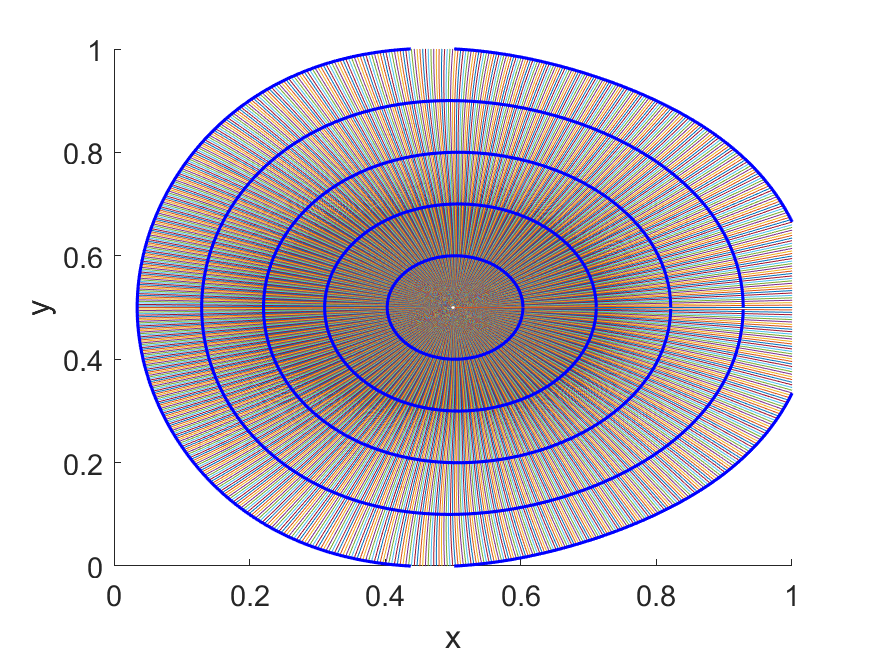}}
     \caption{Example 1. (a) The velocity; (b) Rays and wavefronts with source $\bx_0=[0.5,0.5]$. The thick blue lines represent equal-time wavefronts (traveltime contours) with the contour interval equal to $0.1$, and thin colored lines represent rays with different take-off angles.}
     \label{figure10}
     \vspace{-3.75mm}
     \end{figure}

We first illustrate how the sampling rate of waves in terms of points per wavelength (PPW) affects the accuracy of wave propagation. Since, essentially, a  wave-resolution-satisfying discretization of the computational domain affects the accuracy of numerical differentiation and interpolation of oscillatory wave data, which in turn will influence the accuracy of the Hadamard integrator, we will  appreciate such discretization effect by using numerical experiments. Of course, at the outset of such experiments, we can avoid numerical errors from differentiation and interpolation by using the given {\it exact} initial conditions in their analytic forms. Therefore, we will apply the HKH propagator for one step to solve the wave equation in the setting of Example 1 with the following parameters: $\beta=64$, $h=\frac{1}{320}$, and $M_1=M_2=128$, where we will use different initialization methods to fulfill the initial data on the underlying wavefronts.

Specifically, we consider the following four fulfilling cases:
\begin{enumerate}
    \item (PPW$\sim$5) Cubic-spline interpolation created under a mesh with step size $\bar{h}=\frac{1}{160}$;
    \item (PPW$\sim$10) Cubic-spline interpolation created under a mesh with step size $\bar{h}=\frac{1}{320}$;
    \item (PPW$\sim$20) Cubic-spline interpolation created under a mesh with step size $\bar{h}=\frac{1}{640}$;
    \item (Analytical expressions) Analytical expressions of the initial conditions and their derivatives are used directly.
\end{enumerate}

The numerical errors compared with the exact solution for the wavefield at $T=0.1$ are shown in Figure \ref{figure9}, and the relative $L^2$ and $L^{\infty}$ errors  are shown in Table \ref{tbl1}. Since we use the same low-rank representations and Gaussian-quadrature nodes for the four cases, we can conclude that different error behaviors are due to different initializations. We also observe that both the initialization by PPW=20 and the analytical initialization yield the similar level of accuracy, indicating that an over-resolved numerical initialization will achieve the same accuracy as an analytical initialization. However, since we are recursively applying the HKH propagator to solve the wave equation, at intermediate steps an analytical initialization is not available and an over-resolved numerical initialization is expensive; consequently, as a trade-off between accuracy and computational complexity, we choose to uniformly discretize the computational domain by using PPW$\sim 10$ which suffices to resolve wave motion.

          \begin{figure}[htbp]
     \centering
          \subfigure[]{
     \includegraphics[scale=0.35]{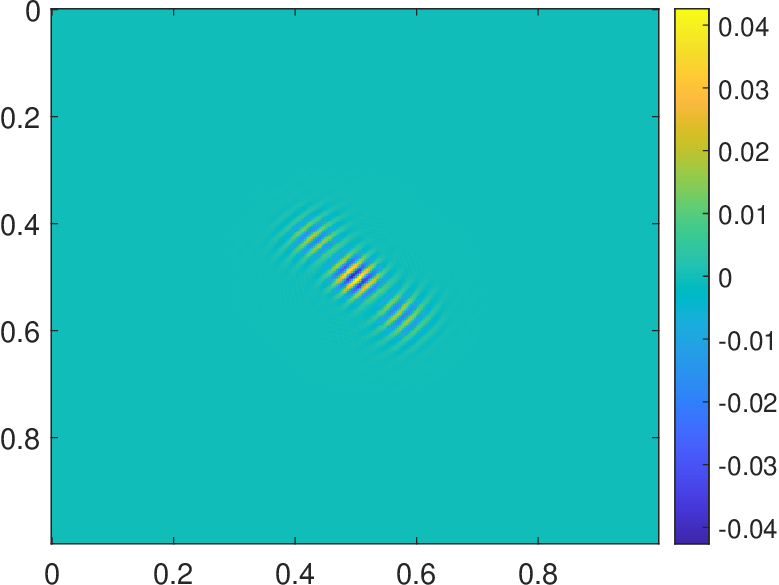}}
     \hspace{8mm}
          \subfigure[]{
          \includegraphics[scale=0.35]{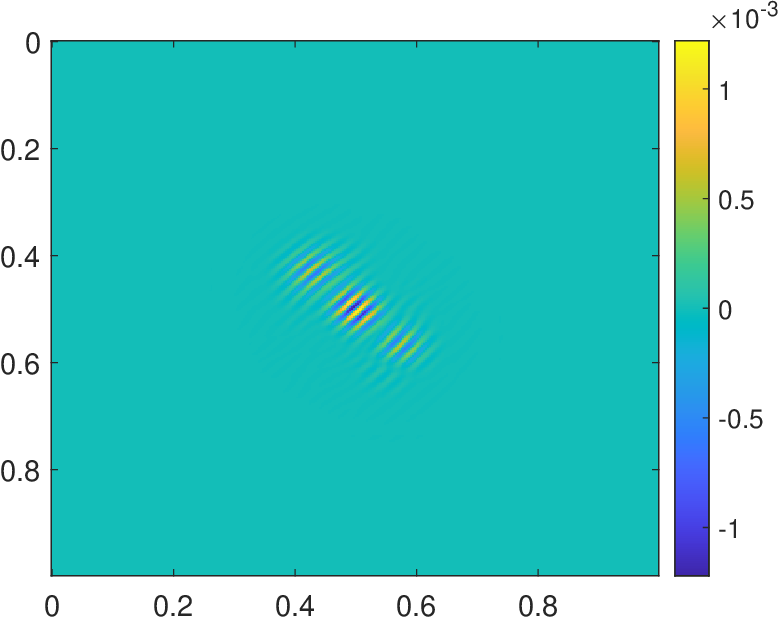}}\\
          \subfigure[]{
     \includegraphics[scale=0.35]{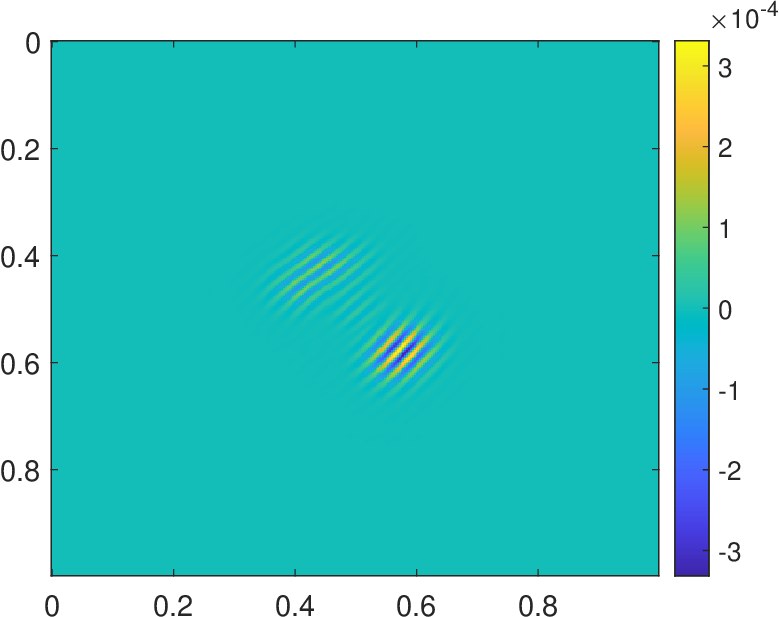}}
     \hspace{8mm}
      \subfigure[]{
          \includegraphics[scale=0.35]{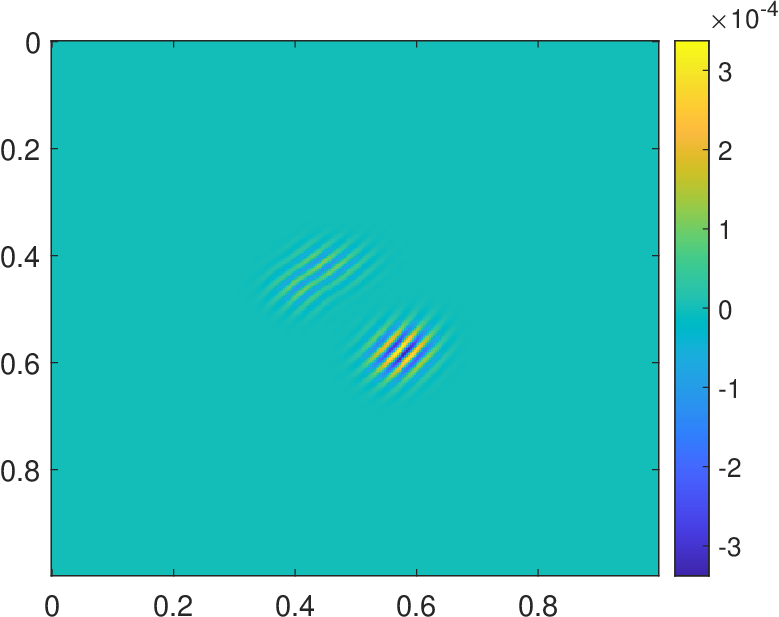}}\\
     \caption{Example 1. $T=0.1$. Errors with different initializations. (a) PPW$\sim5$; (b) PPW$\sim10$; (c) PPW$\sim 20$; (d) Analytical expressions }\label{figure9}
     \vspace{-3.75mm}
     \end{figure}

\begin{table}[width=1\linewidth,cols=4,pos=h]
\caption{The relative $L^2$ and $L^{\infty}$ errors at $T=0.1$ with different initializations.}\label{tbl1}
\begin{tabular*}{\tblwidth}{@{} LLLLL@{} }
\toprule
Initialization methods & PPW $\sim$ 5 & PPW $\sim$ 10 & PPW $\sim$ 20 &  Analytical expressions\\
\midrule
Relative $L^2$ error & $6.8e-2$ & $2.0e-3$ & $6.0e-4$ & $6.0e-4$ \\
Relative $L^{\infty}$ error & $8.9e-2$ & $2.5e-3$ & $7.0e-4$ & $7.0e-4$\\
\bottomrule
\end{tabular*}
\end{table}

          \begin{figure}[htbp]
     \centering
               \subfigure[]{
          \includegraphics[scale=0.35]{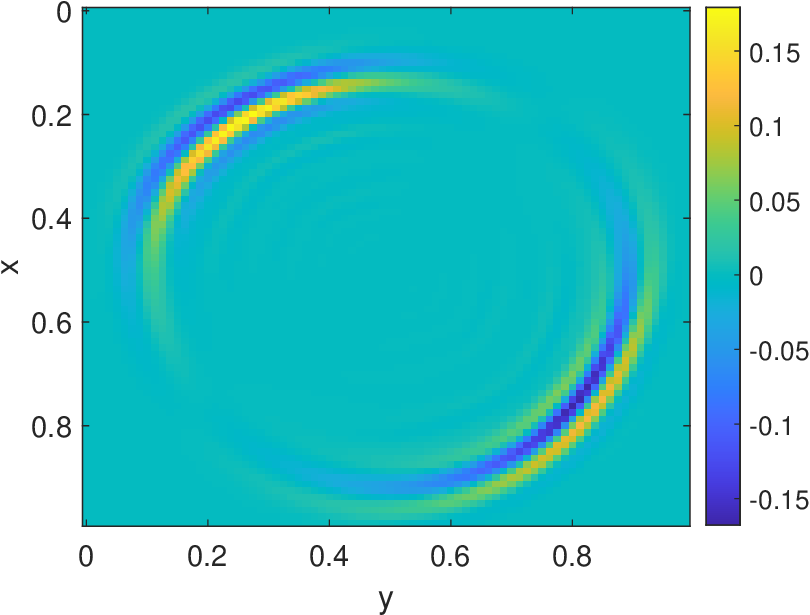}}
          \hspace{8mm}
          \subfigure[]{
     \includegraphics[scale=0.35]{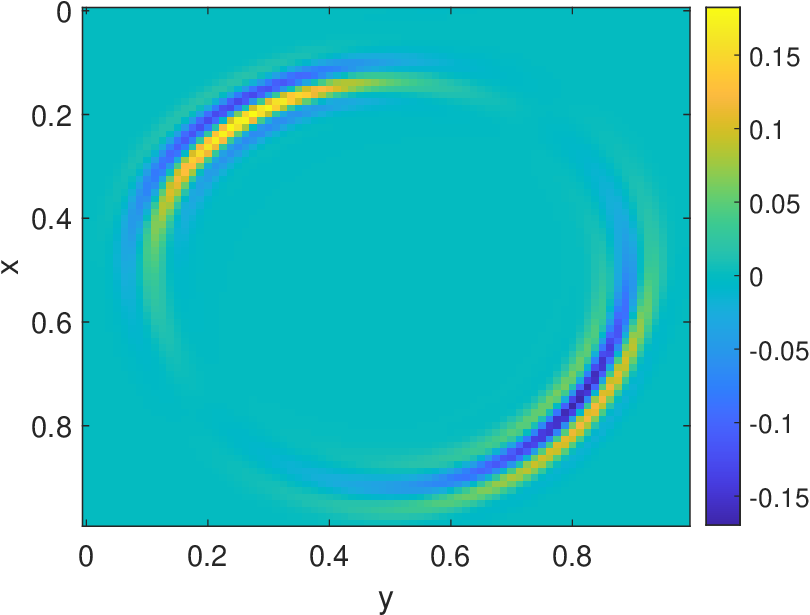}}\\
          \subfigure[]{
          \includegraphics[scale=0.35]{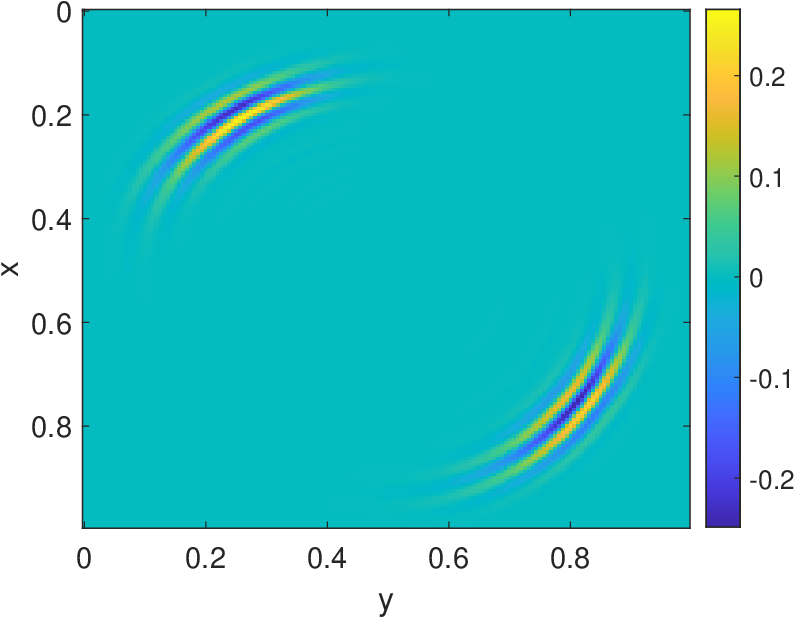}}
          \hspace{8mm}
          \subfigure[]{
     \includegraphics[scale=0.35]{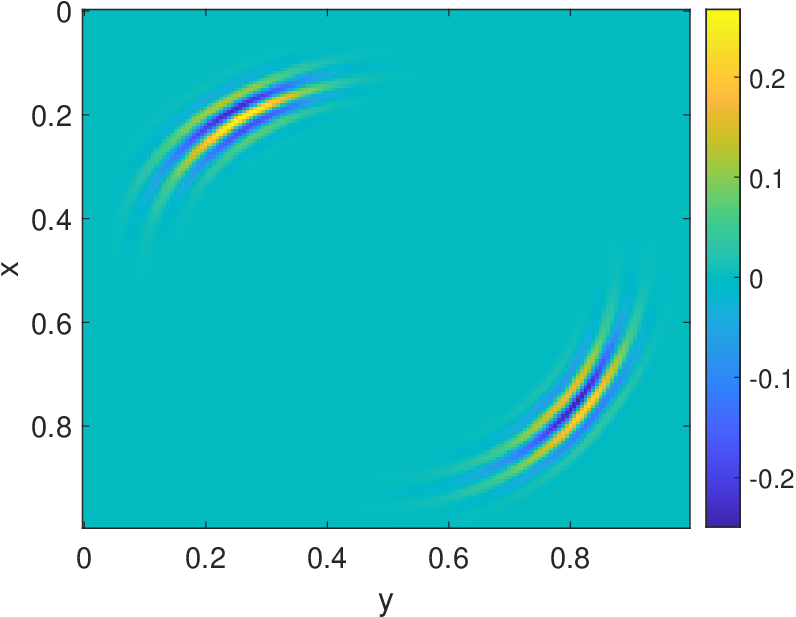}}\\
          \subfigure[]{
          \includegraphics[scale=0.35]{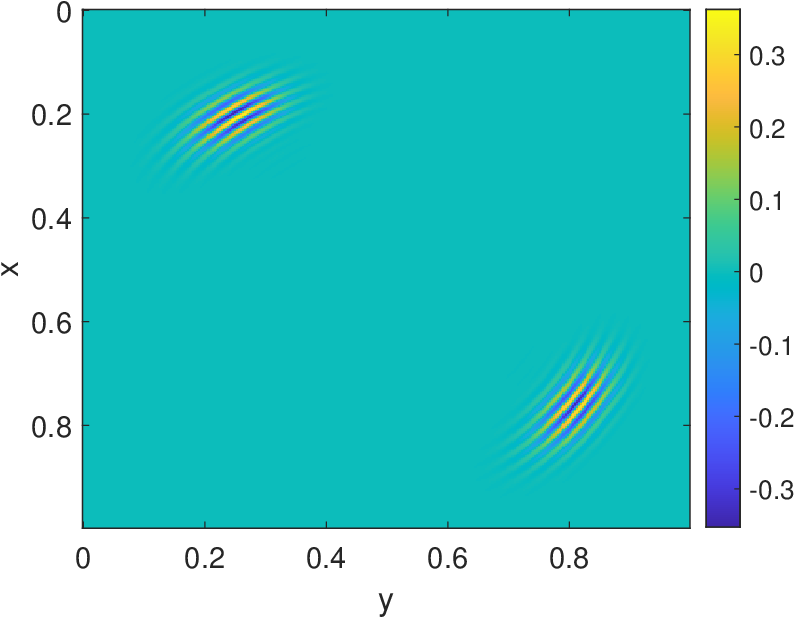}}
          \hspace{8mm}
          \subfigure[]{
     \includegraphics[scale=0.35]{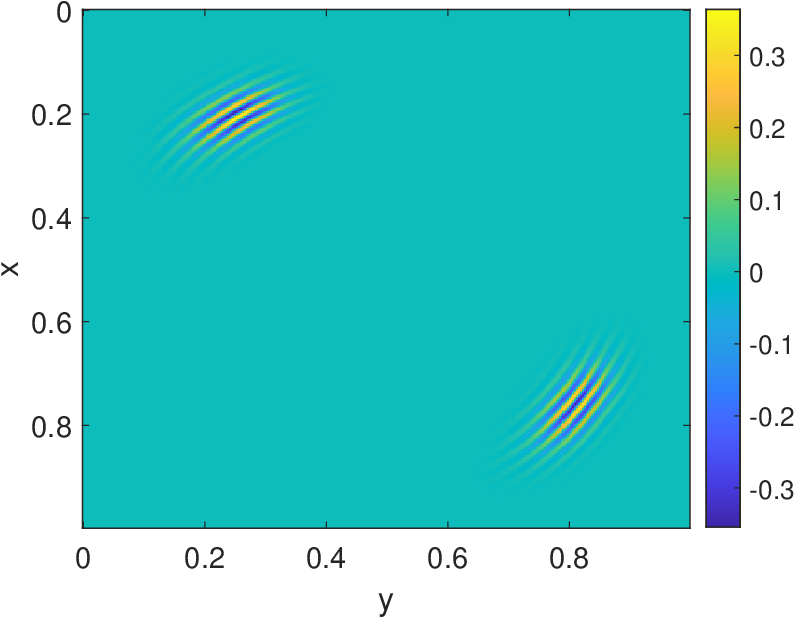}}\\
     \subfigure[]{
     \includegraphics[scale=0.35]{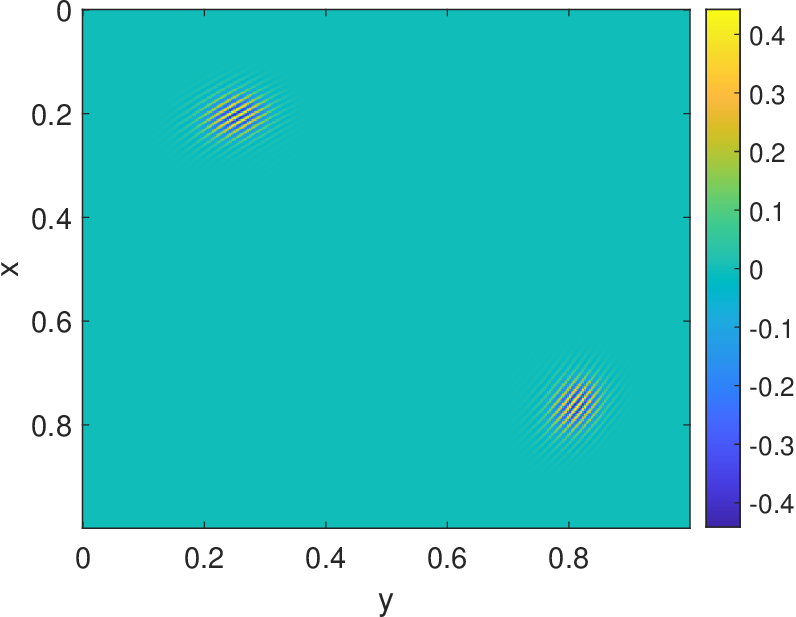}}
     \hspace{8mm}
     \subfigure[]{
     \includegraphics[scale=0.35]{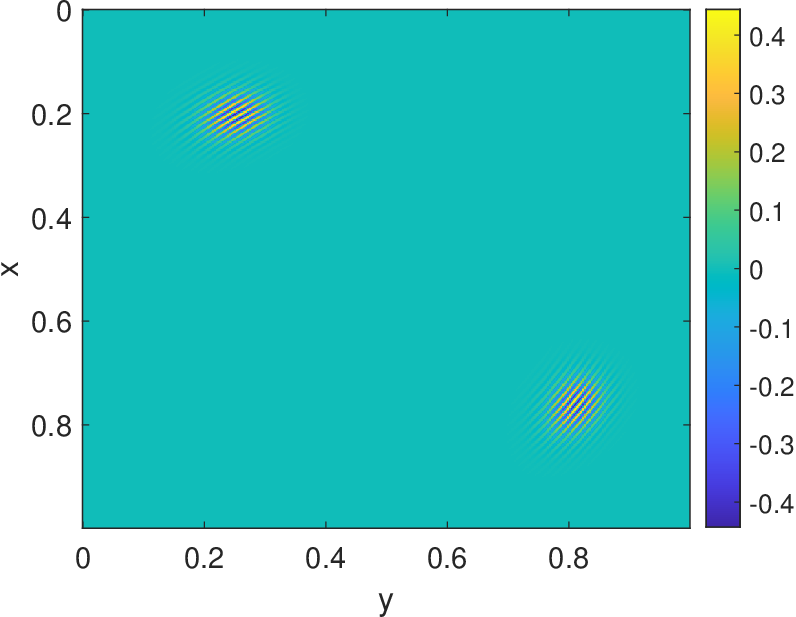}}\\
     \caption{Example 1. $T=0.4$. (a) RT solution with $\beta=16$; (b) Exact solution with $\beta=16$; (c) RT solution with $\beta=32$; (d) Exact solution with $\beta=32$; (e) RT solution with $\beta=64$; (f) Exact solution with $\beta=64$; (g) RT solution with $\beta=128$; (h) Exact solution with $\beta=128$.}\label{figure11}
     \vspace{-3.75mm}
     \end{figure}

\begin{figure}[htbp]
     \centering
               \subfigure[]{
          \includegraphics[scale=0.35]{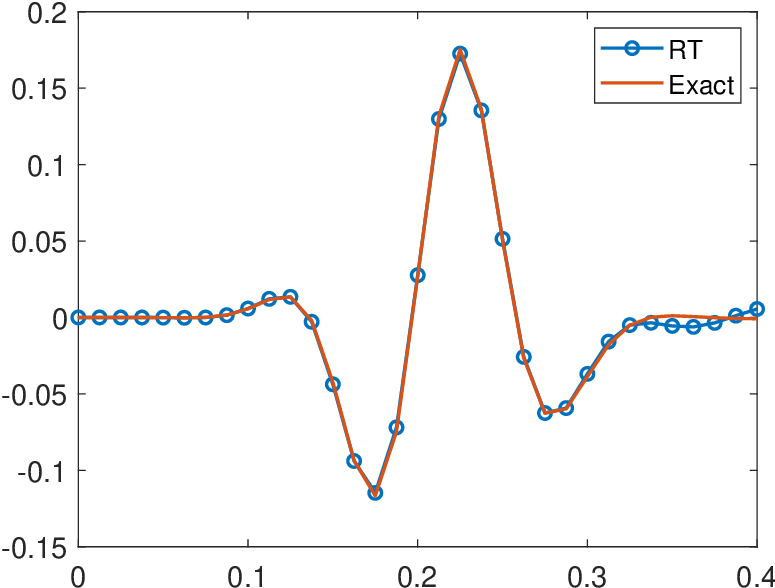}}
          \hspace{8mm}
          \subfigure[]{
     \includegraphics[scale=0.35]{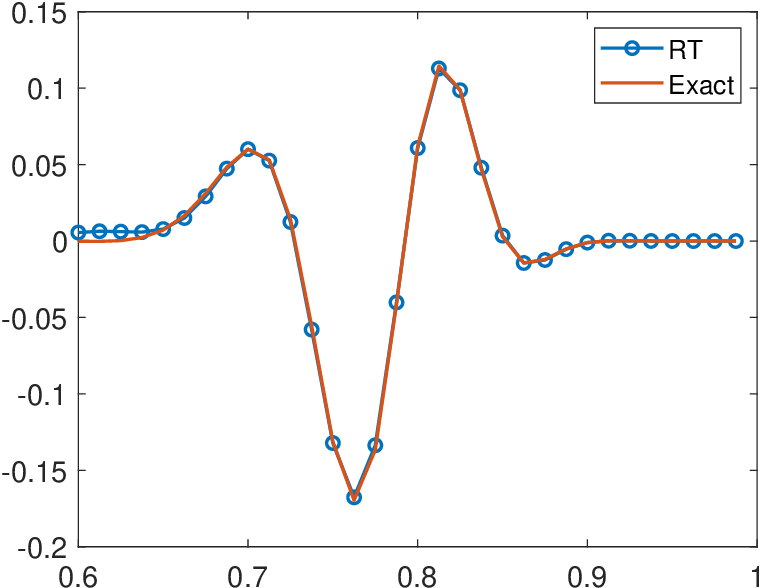}}\\
          \subfigure[]{
          \includegraphics[scale=0.35]{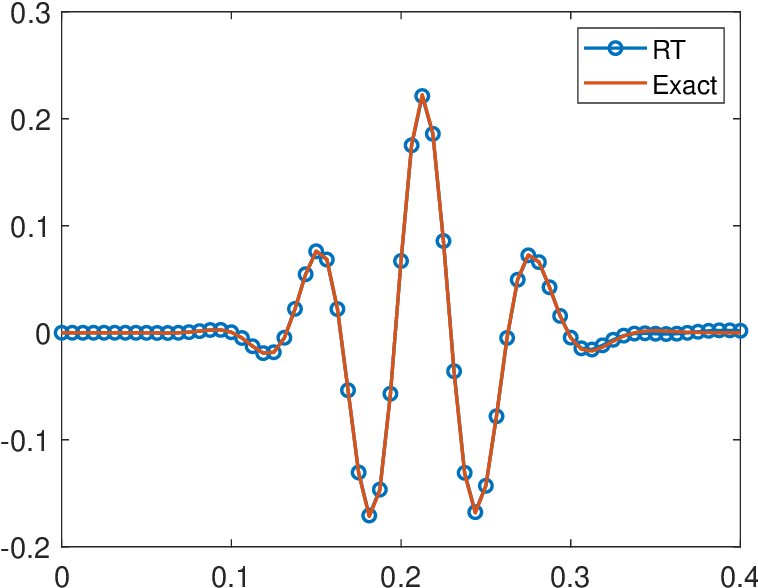}}
          \hspace{8mm}
          \subfigure[]{
     \includegraphics[scale=0.35]{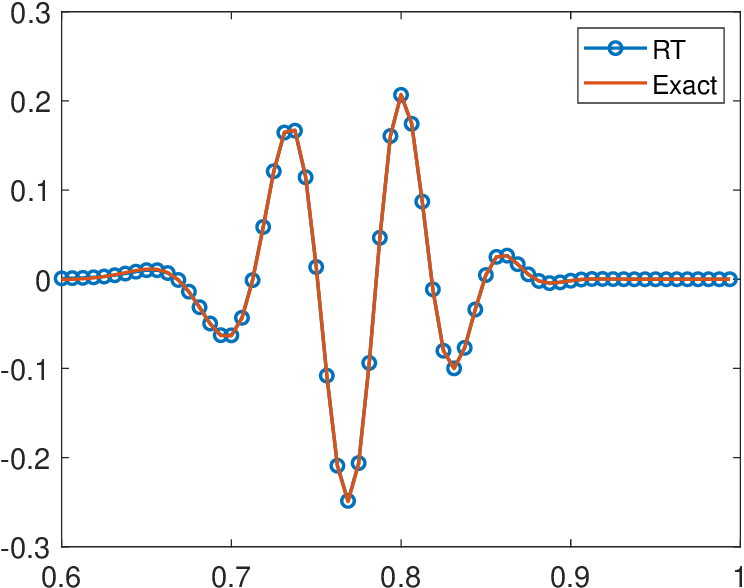}}\\
          \subfigure[]{
          \includegraphics[scale=0.35]{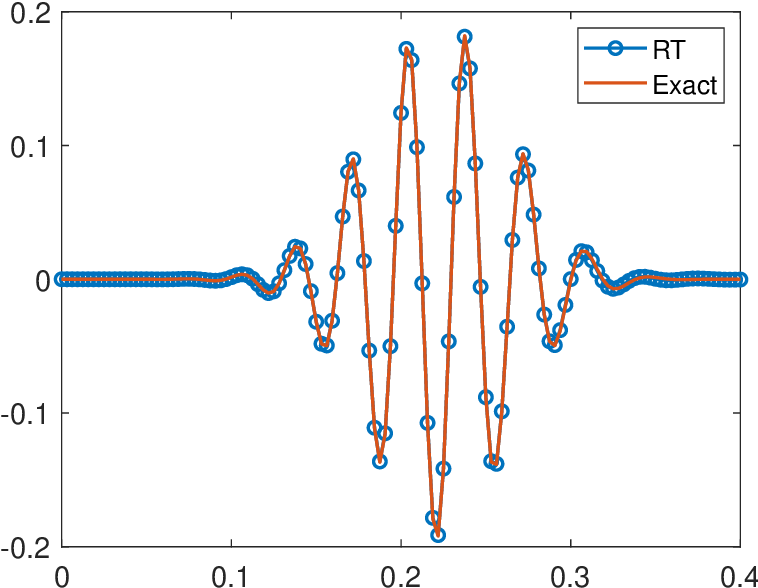}}
          \hspace{8mm}
          \subfigure[]{
     \includegraphics[scale=0.35]{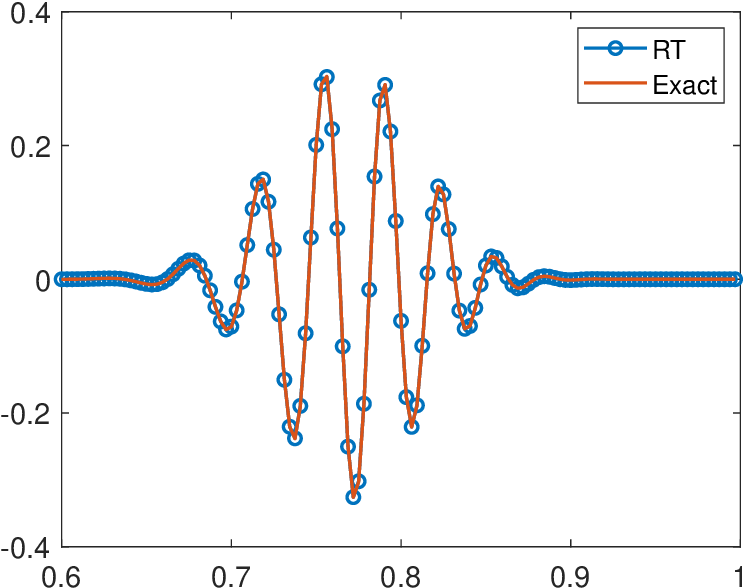}}\\
     \subfigure[]{
     \includegraphics[scale=0.35]{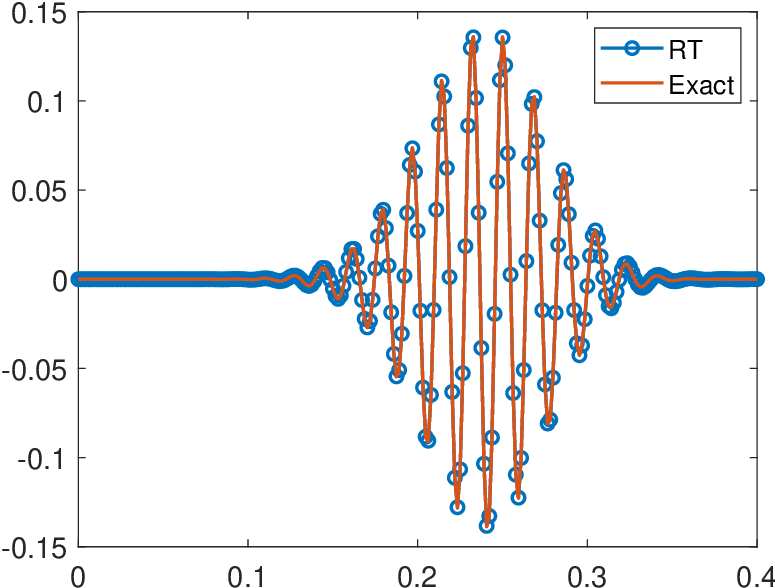}}
     \hspace{8mm}
     \subfigure[]{
     \includegraphics[scale=0.35]{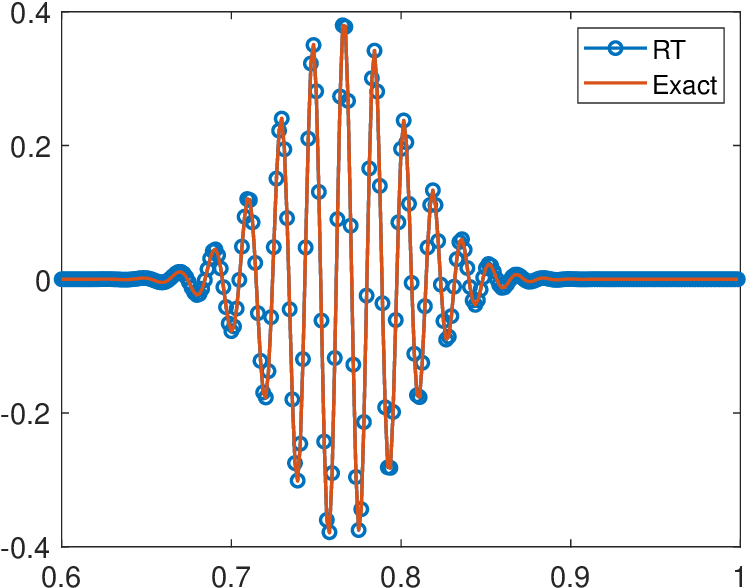}}\\
     \caption{Example 1. $T=0.4$. (a) a windowed slice at $x=0.25$ with $\beta=16$;
 (b) a windowed slice at $y=0.8$ with $\beta=16$; (c) a windowed slice at $x=0.25$ with $\beta=32$; (d) a windowed slice at $y=0.8$ with $\beta=32$; (e) a windowed slice at $x=0.25$ with $\beta=64$; (f) a windowed slice at $y=0.8$ with $\beta=64$; (e) a windowed slice at $x=0.25$ with $\beta=128$; (f) a windowed slice at $y=0.8$ with $\beta=128$.}\label{figure12}
 \vspace{-3.75mm}
     \end{figure}

Next, we consider initial conditions with variable $\beta$ to illustrate the effectiveness of the Hadamard integrator for highly oscillatory wavefields. We set $\beta=16,32,64,128$ and $h=\frac{1}{80}, \frac{1}{160}, \frac{1}{320}, \frac{1}{640}$, respectively, and we present the RT solutions and exact solutions at $T=0.4$ in Figure \ref{figure11}. Figure \ref{figure12} shows line comparisons of the two solutions, where the exact solution in `-' overlays the RT solution in `o'.
The RT solutions match the exact solution well, especially when the $\beta$ is large.

The relative $L^\infty$ and $L^2$ errors of RT solutions with different $\beta$ are shown in Figure \ref{figure13}. The relative errors drop significantly and stay below $1\%$ as $\beta$ increases, and such error behaviors are typical for microlocal analysis based numerical methods for high-frequency wave propagation \cite{liusonburqia23}. While a detailed numerical analysis of our new method is an ongoing work, we provide a brief explanation for such a phenomenon here. The overall error of our Hadamard integrator consists of three parts: the leading-term truncation error of the asymptotic series in terms of $\beta$, the  numerical errors due to approximations to oscillatory data, Hadamard ingredients, and low-rank representations, and the interaction of the asymptotic error and the numerical errors; when $\beta$ is small, the asymptotic error dominates over the other errors so that we can observe the obvious error reduction when $\beta$ increases; but when $\beta$ is large enough, the asymptotic error is no longer dominant over the other errors so that the error reduction saturates when $\beta$ increases.

\begin{figure}[htbp]
     \centering
     \includegraphics[scale=0.385]{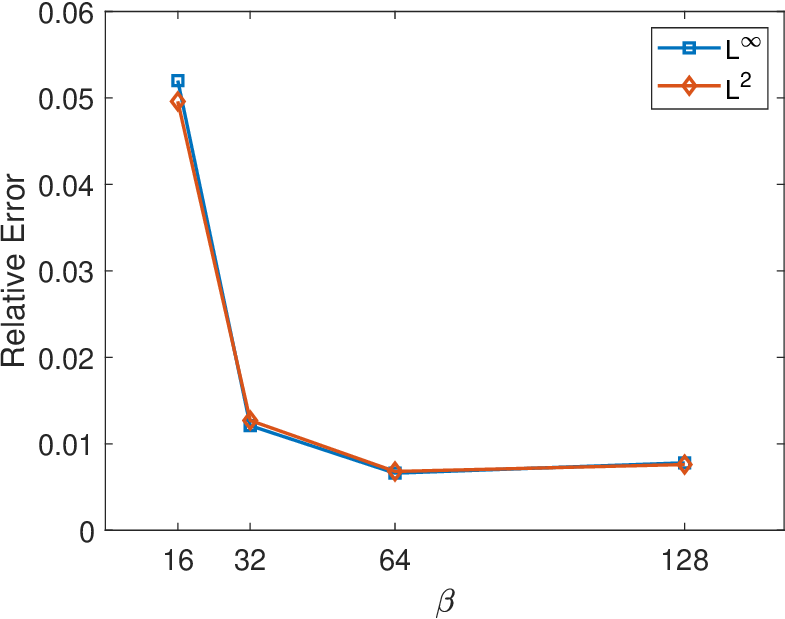}
     \caption{The relative $L^\infty$ and $L^2$ errors of RT solutions with different $\beta$}
     \label{figure13}
     \vspace{-3.75mm}
     \end{figure}

\noindent \textbf{Example 2.} We set up the problem as the following.
    \begin{itemize}
    \item[$\bullet$] $\rho=\frac{4}{(1+0.5\sin(5\pi x)\cos(3\pi y))^2},
    \;$ $\nu=1,\;$ $c=0.5(1+0.5\sin(5\pi x)\cos(3\pi y)).$
    \item[$\bullet$] $u(0,x,y)=\sin \left( 64\pi(x+y-1)\right) \exp \left(-600\left((x-0.5)^{2}+(y-0.5)^{2}\right)\right)$, and $\;u_t(0,x,y)=0.$
    \item[$\bullet$]The computational domain is $\Omega=[0,1]^2$ and the grid size used to discretize $\Omega$ is $h=\frac{1}{640}.$
    \item[$\bullet$]The orders of the tensorized Chebyshev polynomials with respect to $[x_0,y_0,\tau,\theta_0]$ are $[31,31,21,21].$
    \item[$\bullet$]The numbers of Gaussian-quadrature nodes are $M_1=M_2=192.$
\end{itemize}

This example serves the purpose of illustrating that the Hadamard integrator maintains accuracy in a medium that induces caustics and is able to treat spatially overturning waves smoothly. In Figure \ref{figure17}, we show the velocity and the ray diagram overlaid with several wavefronts (traveltime isocontours) emitted from the source point at $[0.5,0.5]$, where caustics occur  when the rays form envelopes and the wavefronts self-intersect.

To evolve the wavefield in this example, we take $\Delta t=0.1$ to construct the HKH propagator. To handle the rapid variation in the velocity model, we increase the orders of Chebyshev polynomials but still keep PPW$\sim$10 by considering the smallest  wavelength. To construct the low-rank representations, we solve the system (\ref{4.9}) using the Runge-Kutta method(RK4) for $1000$ steps, and we solve the eikonal and transport equations with mesh size $h_1=0.002$ in the square neighborhood (\ref{4.11.1}) with $H=0.1$.

\begin{figure}[htbp]
     \centering
     \subfigure[]{
     \includegraphics[scale=0.385]{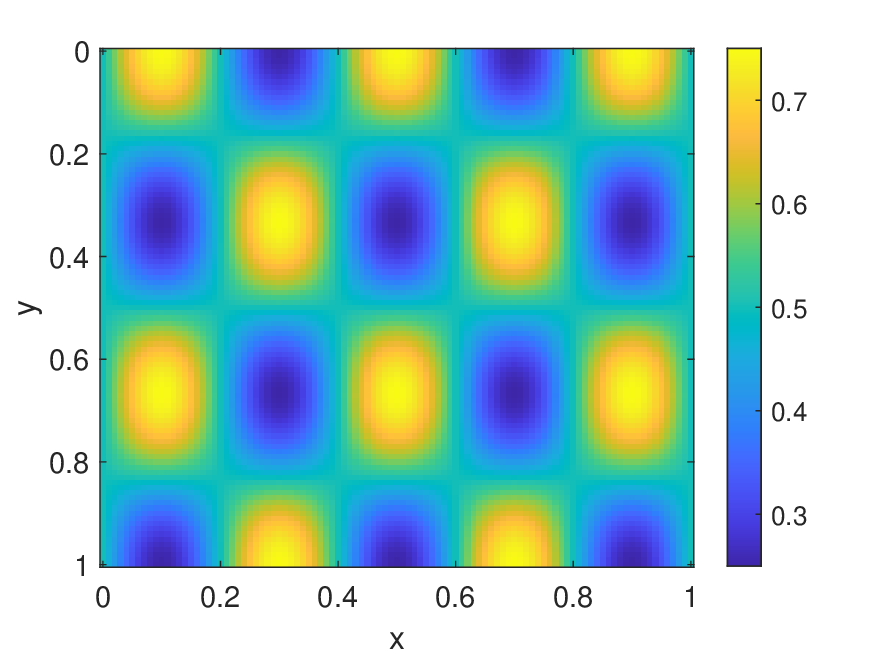}}
     \hspace{8mm}
     \subfigure[]{
     \includegraphics[scale=0.385]{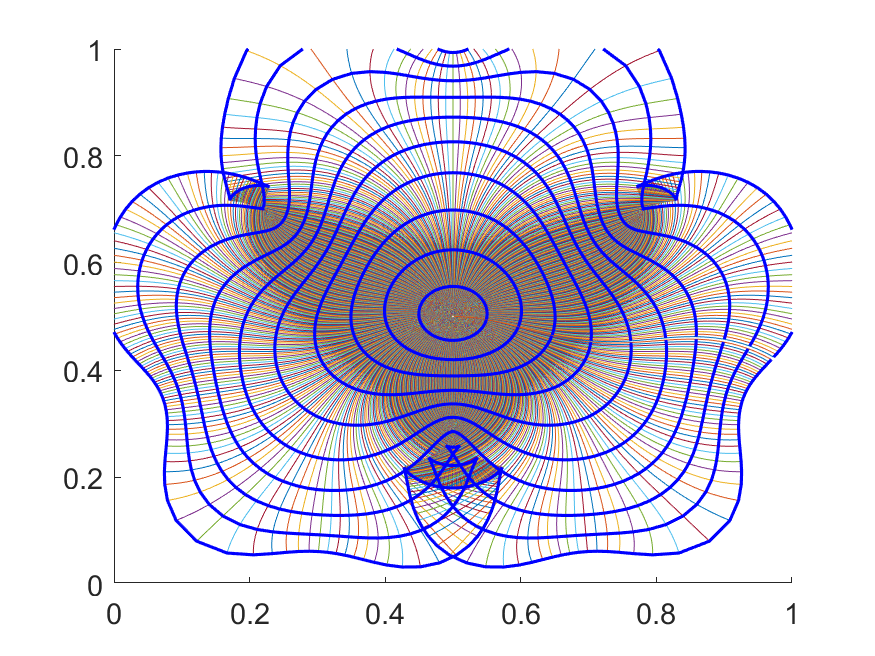}}
     \caption{Example 2. (a) The velocity; (b) The rays and wavefronts with source $\bx_0=[0.5,0.5]$. The thick blue lines represent the equal-time wavefronts with the contour interval $0.1$, and thin colored lines represent rays from different take-off angles.}
     \label{figure17}
     \vspace{-3.75mm}
     \end{figure}

 We show in Figure \ref{figure18} the wavefields at $T=0,0.1,\cdots,0.8$. We can perceive the effect of caustics from the variations of the wave amplitude which imply that the rays are gradually focusing at caustics.
\begin{figure}[htbp]
     \centering
     \subfigure[]{
     \includegraphics[scale=0.35]{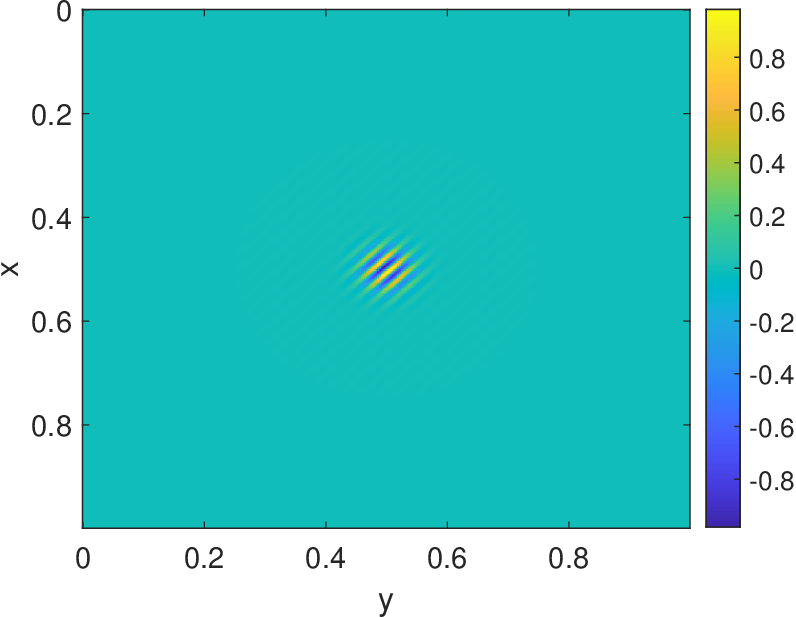}}
     \subfigure[]{
     \includegraphics[scale=0.35]{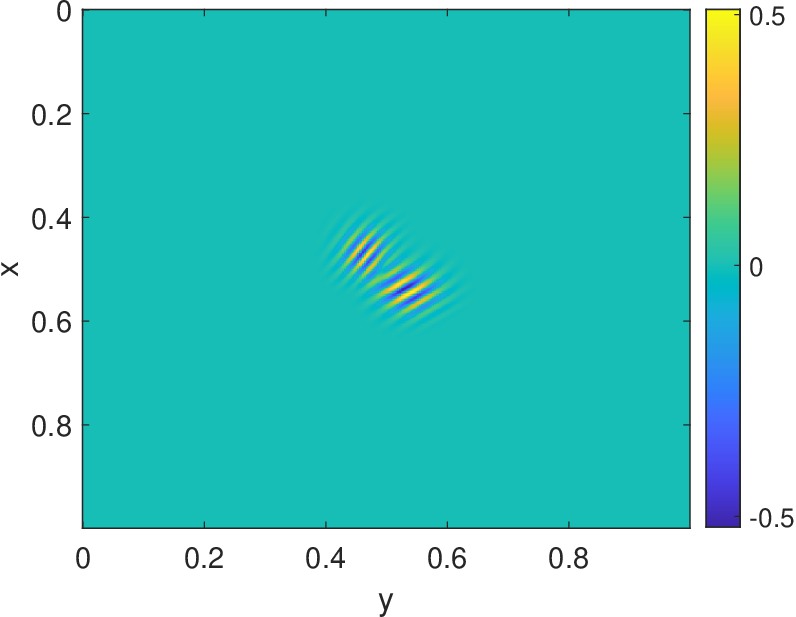}}
     \subfigure[]{
     \includegraphics[scale=0.35]{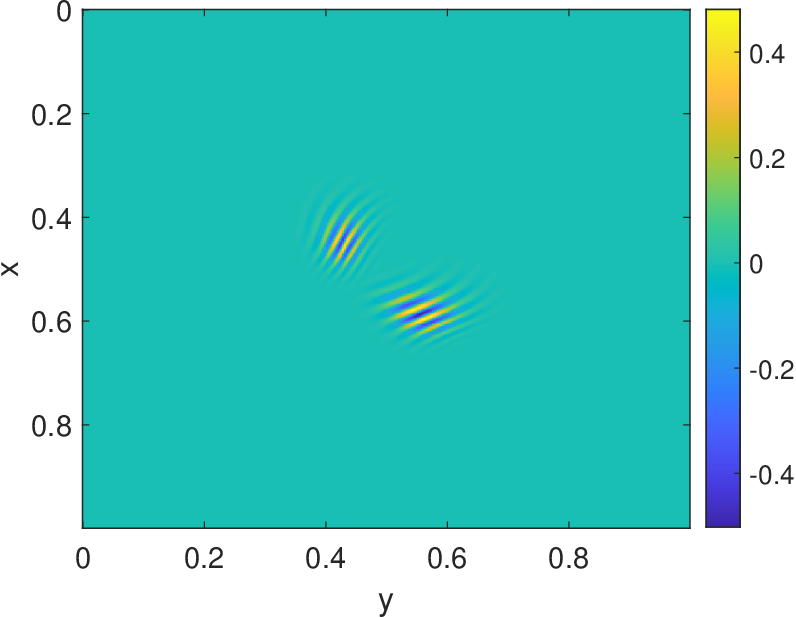}}\\
     \subfigure[]{
          \includegraphics[scale=0.35]{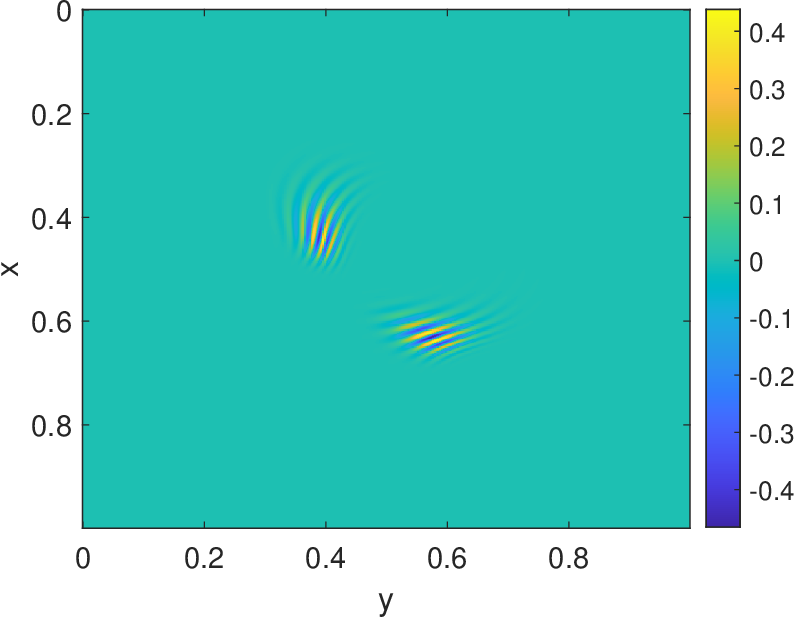}}
          \subfigure[]{
     \includegraphics[scale=0.35]{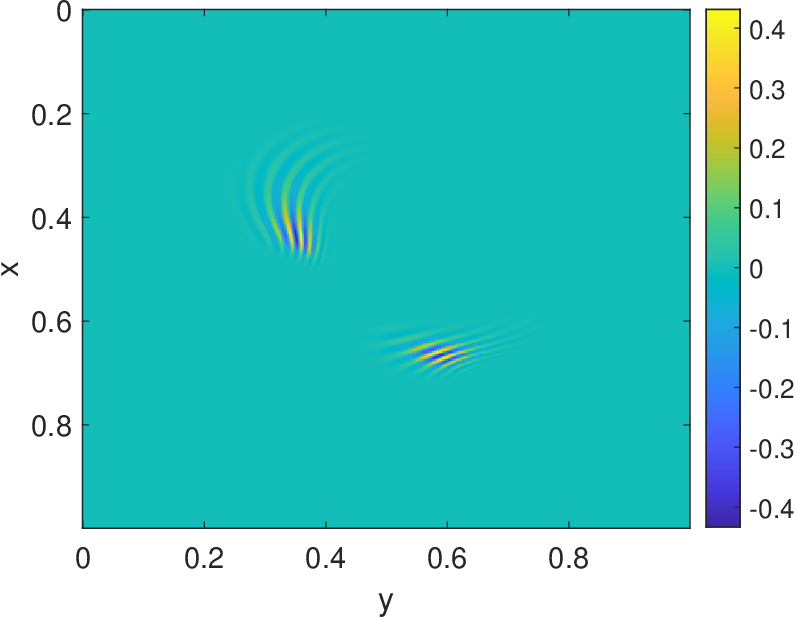}}
     \subfigure[]{
     \includegraphics[scale=0.35]{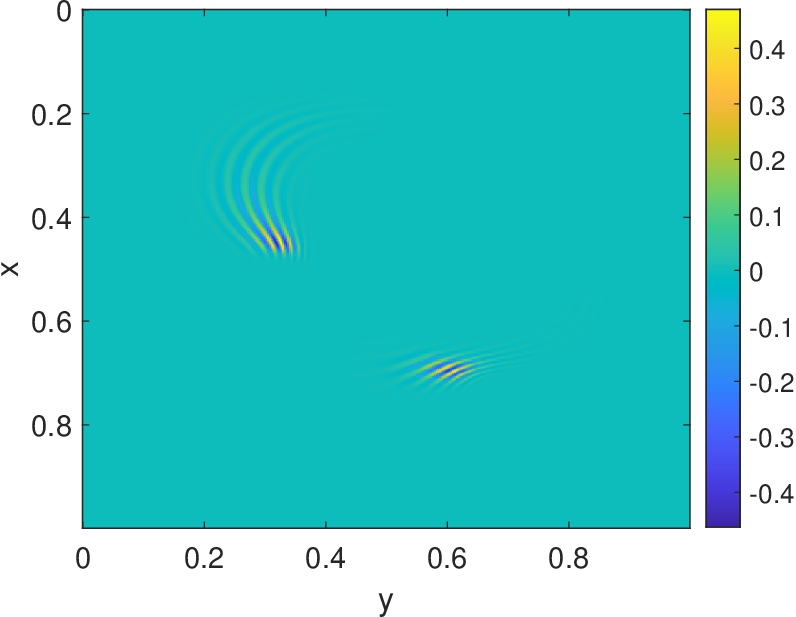}}\\
     \subfigure[]{
          \includegraphics[scale=0.35]{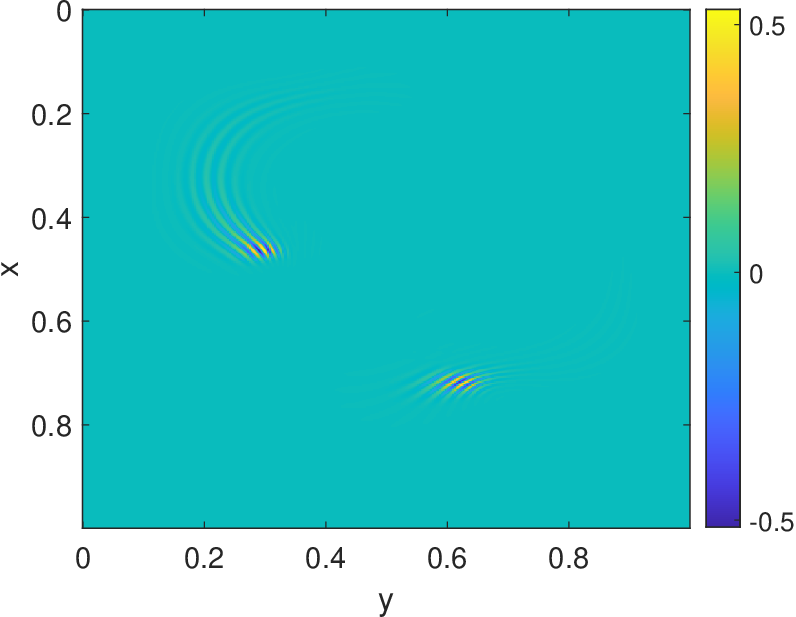}}
          \subfigure[]{
     \includegraphics[scale=0.35]{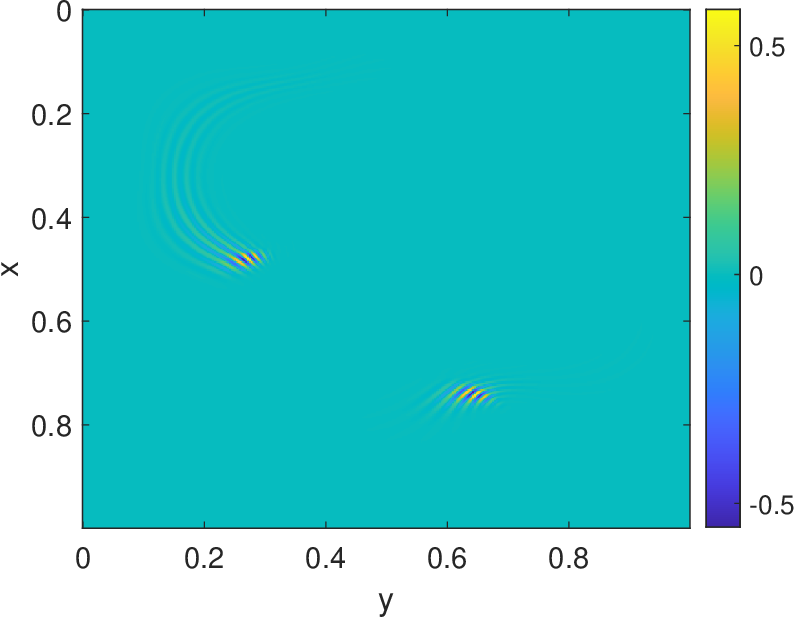}}
     \subfigure[]{
     \includegraphics[scale=0.35]{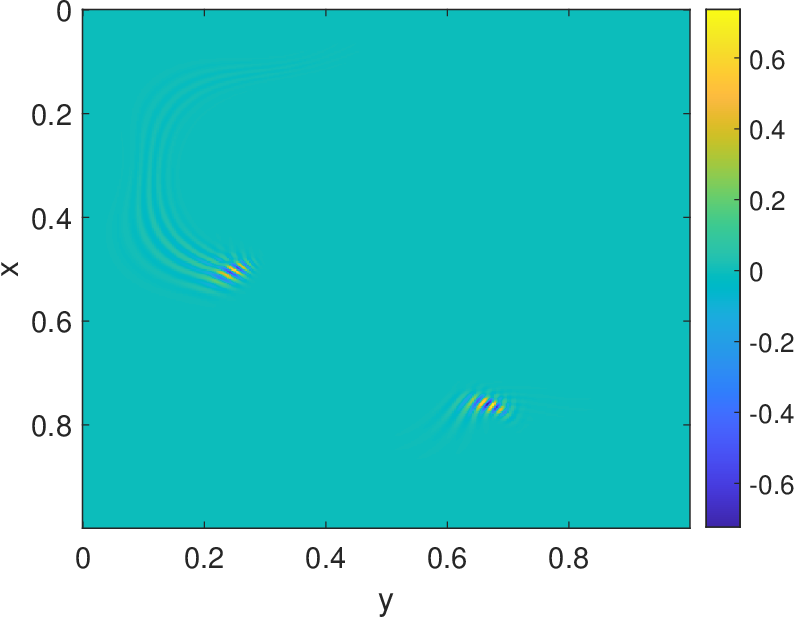}}\\
     \caption{Example 2: (a), (b), $\cdots$, (i): RT solutions at $T=0,0.1,\cdots,0.8$, respectively.}\label{figure18}
     \end{figure}

     Compared to the exact solution, the relative maximum error of the Hadamard  integrator at $T=0.8$ is $1.4\%$ and the relative $L^2$ error is $1.6\%$. To further compare the results at $T=0.8$, we overlay the exact solution (`-') on the RT solution (`o') along different lines in Figure \ref{figure19}. We can observe that the RT solution consistently matches the exact solution very well. As expected, the Hadamard integrator accurately propagates wavefields in a caustic-inducing medium.

          \begin{figure}[htbp]
     \centering
     \subfigure[]{
     \includegraphics[scale=0.35]{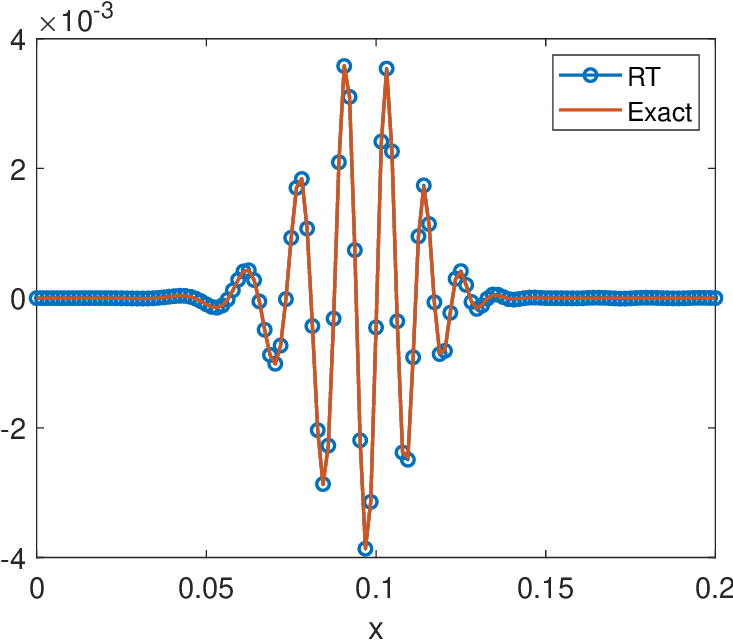}}
     \subfigure[]{
     \includegraphics[scale=0.35]{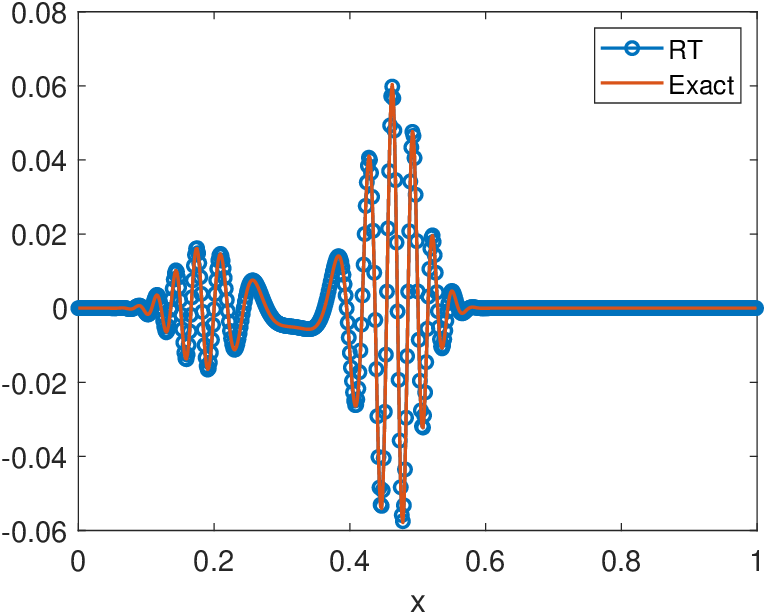}}
          \subfigure[]{
     \includegraphics[scale=0.35]{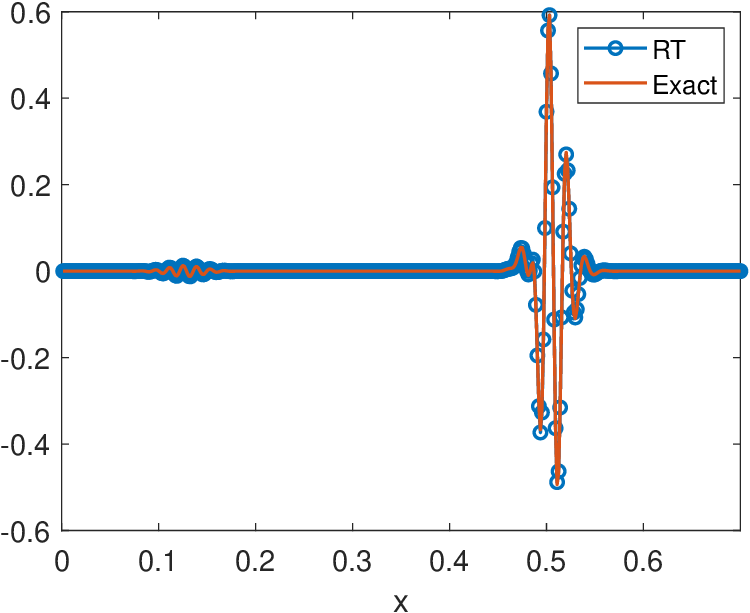}}\\
     \subfigure[]{
     \includegraphics[scale=0.35]{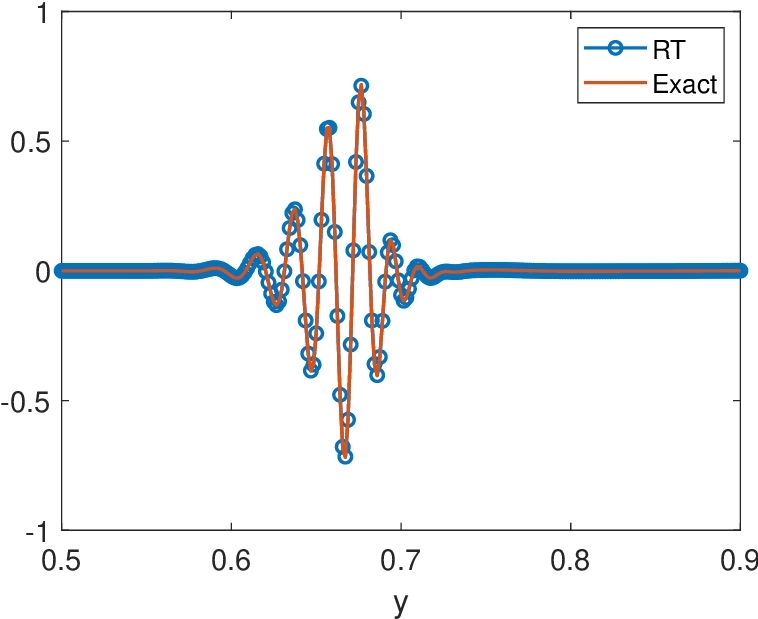}}
          \subfigure[]{
     \includegraphics[scale=0.35]{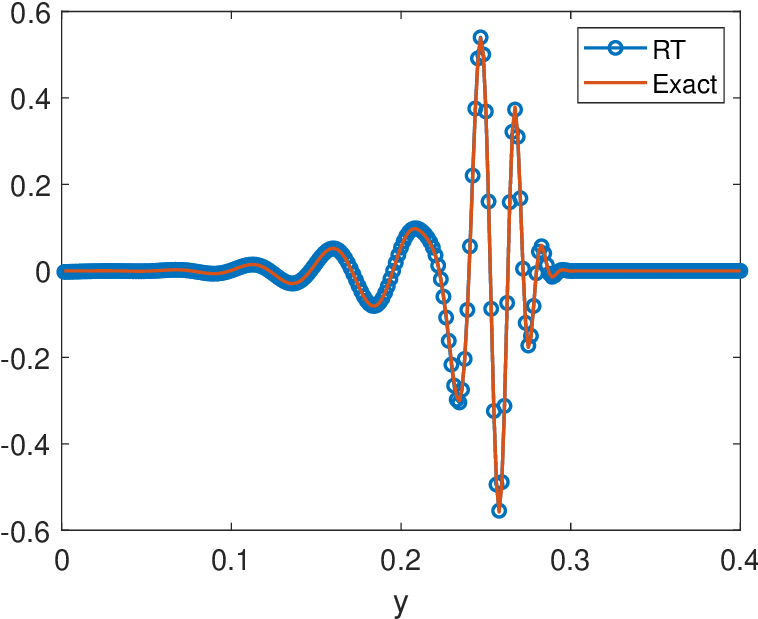}}
     \subfigure[]{
     \includegraphics[scale=0.35]{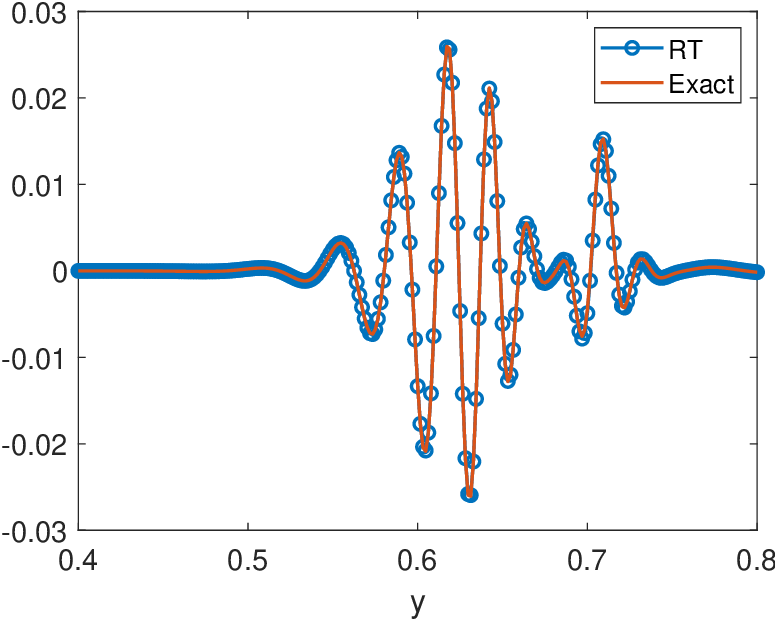}}
     \caption{Example 2. (a) a windowed slice at $y=0.4$; (b) a slice at $y=0.15$; (c) a windowed slice at $y=0.4$; (d) a windowed slice at $x=0.7625$; (e) a windowed slice at $x=0.5$; (f) a windowed slice at $x=0.8$;}\label{figure19}
     \end{figure}

     \subsection{Three-dimensional examples }
   \noindent  \textbf{Example 3.} We use the following setup.
\begin{itemize}
    \item[$\bullet$] $\rho=\frac{1}{(1+0.1\sin(2\pi x)\cos(2\pi y)\sin(2\pi z))^2},\;$ $\nu=1,\;$ and $c=1+0.1\sin(2\pi x)\cos(2\pi y)\sin(2\pi z).$
    \item[$\bullet$] $u(0,x,y,z)=\sin ( \pi \beta(x+y+z-1.5)) \exp \left(-200\left((x-0.5)^{2}+(y-0.5)^{2}+(z-0.5)^2\right)\right)$,\\ $u_t(0,x,y,z)=0.$
    \item[$\bullet$] The computational domain is $\Omega=[0,1]^3$, and the grid size used to discretize $\Omega$ is $h=\frac{1}{5\beta}$.
    \item[$\bullet$] The orders of the tensorized Chebyshev polynomials with respect to $[x_0,y_0,z_0,\theta_0,\xi_0]$ are $[15,15,15,15,15].$
    \item[$\bullet$] The number of Gaussian-quadrature nodes is $M_3=3\beta.$
\end{itemize}

The slice at $z=0.6$ of the velocity model and the rays starting from $\bx_0=[0.5,0.5,0.5]$ are shown in Figure \ref{figure20.2}; there is no caustic in this example.

     \begin{figure}[htbp]
     \centering
     \subfigure[]{
     \includegraphics[scale=0.385]{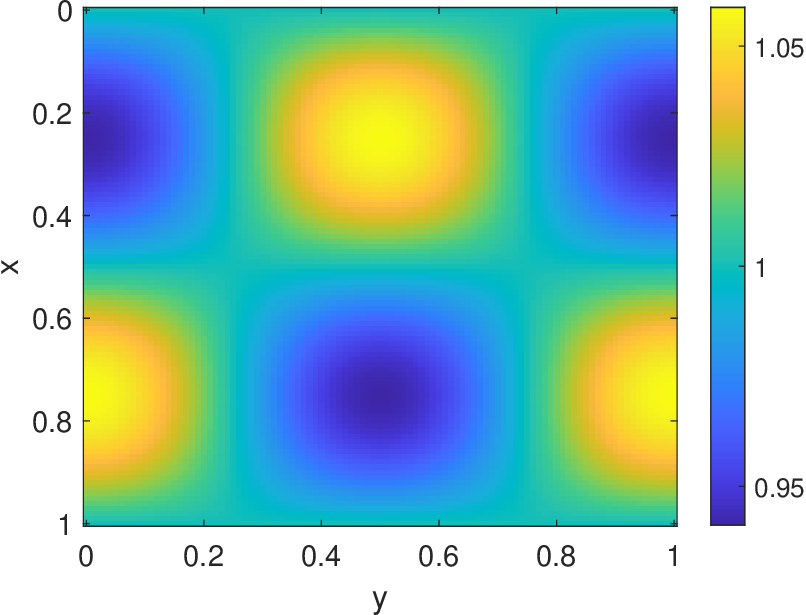}}
     \hspace{8mm}
     \subfigure[]{
     \includegraphics[scale=0.385]{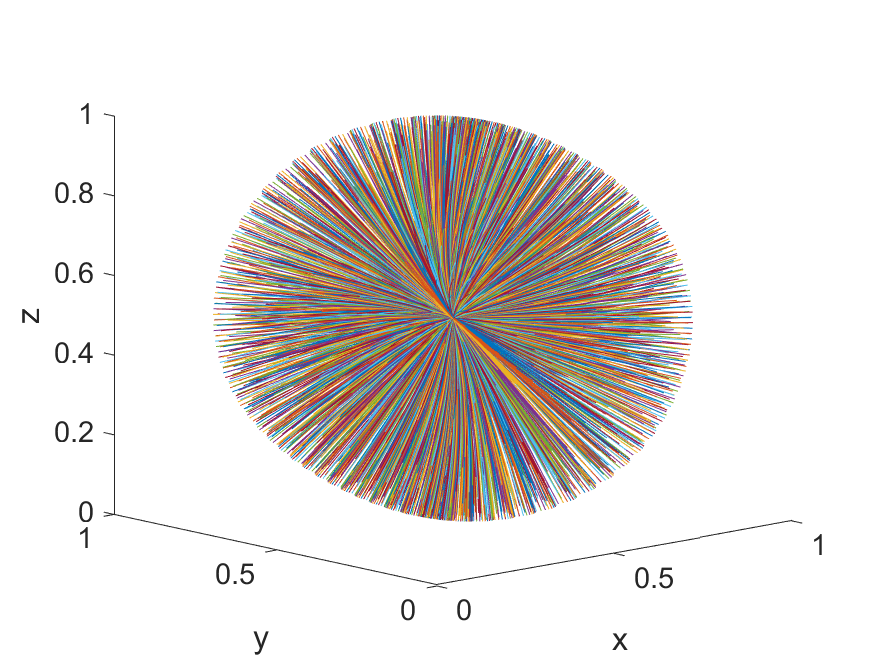}}
     \caption{Example 3. (a) The slice of velocity at $z=0.6$. (b) The rays with different take-off angles starting from the source $\bx_0=[0.5,0.5,0.5]$.}
     \label{figure20.2}
    \vspace{-3.75mm}
     \end{figure}

We set $\Delta t=0.1$ for the HKH propagator. To construct the low-rank representations, we solve the Hamiltonian system (\ref{ode2}) using the Runge–Kutta method (RK4) for $1000$ time steps, and we solve the eikonal and transport equations with mesh size $h_1=0.005$ in the cubic neighborhood (\ref{6.16}) with $H=0.15$.

We first set $\beta=16$, $h=\frac{1}{80}$, and $T=0.4$.
%then the relative $L^{\infty}$ error is $3.3\%$ and the relative $L^{2}$ error is $3.6\%$.
Figure \ref{figure20} shows slices and line comparisons of wavefields.
Figure \ref{figure20}(a) shows a $z$-section of the 3-D wavefield at $z=0.75$, and Figure \ref{figure20}(d) shows the line at $y=0.7$ and $z=0.75$, in which we compare the exact solution (`-') with the RT solution (`o').
Figure \ref{figure20}(b) shows an $x$-section of the 3-D wavefield at $x=0.45$, and Figure \ref{figure20}(e) shows the line at $x=0.45$ and $z=0.2$, in which we compare the exact solution (`-') with the RT solution (`o'). Figure \ref{figure20}(c) shows a $y-$ section of the 3-D wavefield at $y=0.15$, and Figure \ref{figure20}(f) shows a line at $x=0.3$ and $y=0.15$, in which we compare the exact solution (`-') with the RT solution (`o').

We next set $T=0.4$, $\beta=32$ and $h={1}/{160}$, and the results are shown in  Figure \ref{figure21}. Figure \ref{figure21}(a) shows a $z$-section of the 3-D wavefield at $z=0.75$, and Figure \ref{figure21}(d) shows a line at $y=0.7$ and $z=0.75$, in which we compare the exact solution (`-') with the RT solution (`o'). Figure \ref{figure21}(b) shows an $x$-section of the 3-D wavefield at $x=0.15$, and Figure \ref{figure21}(e) shows a line at $x=0.15$ and $z=0.3$, in which we compare the exact solution (`-') with the RT solution (`o'). All the RT solutions match well with the exact solutions. Figure \ref{figure21}(c) shows the $y$-section of the 3-D wavefield at $y=0.425$, and Figure  \ref{figure21}(f) shows a line at $x=0.3$ and $y=0.425$, in which we compare the exact solution (`-') with the RT solution (`o'). As shown, the RT solutions match the exact solutions very well.

          \begin{figure}[htbp]
     \centering
     \subfigure[]{
     \includegraphics[scale=0.35]{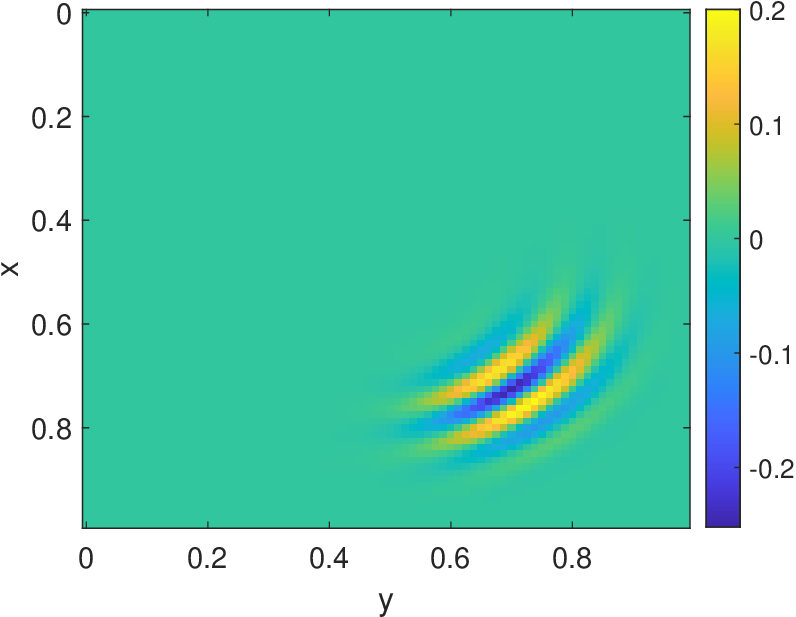}}
          \subfigure[]{
     \includegraphics[scale=0.35]{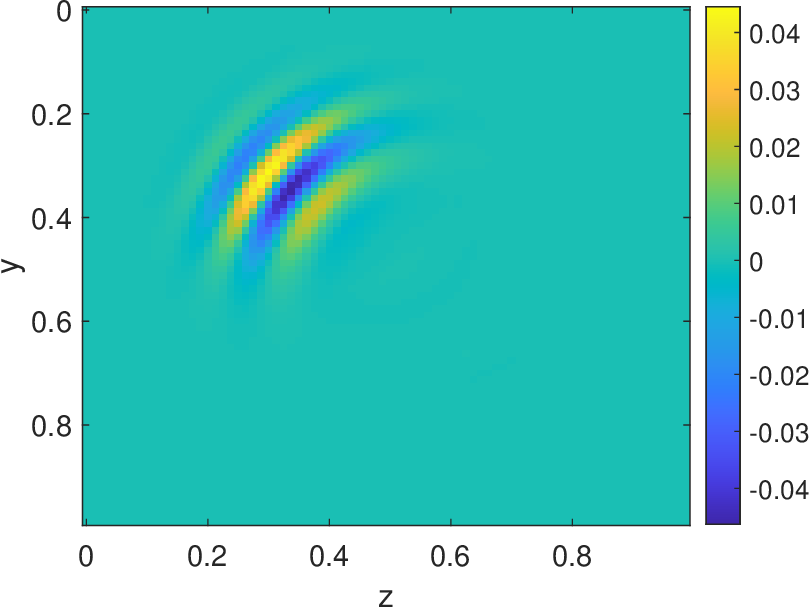}}
          \subfigure[]{
     \includegraphics[scale=0.35]{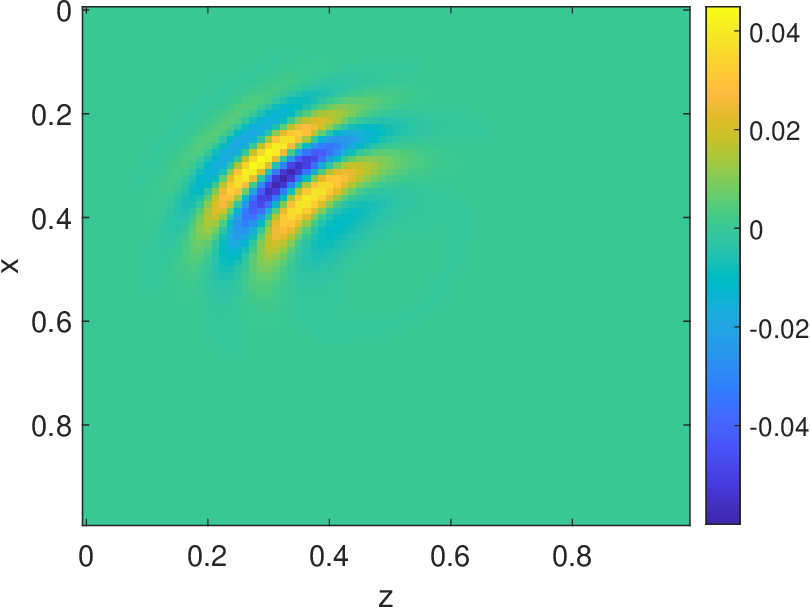}}\\
          \subfigure[]{
     \includegraphics[scale=0.35]{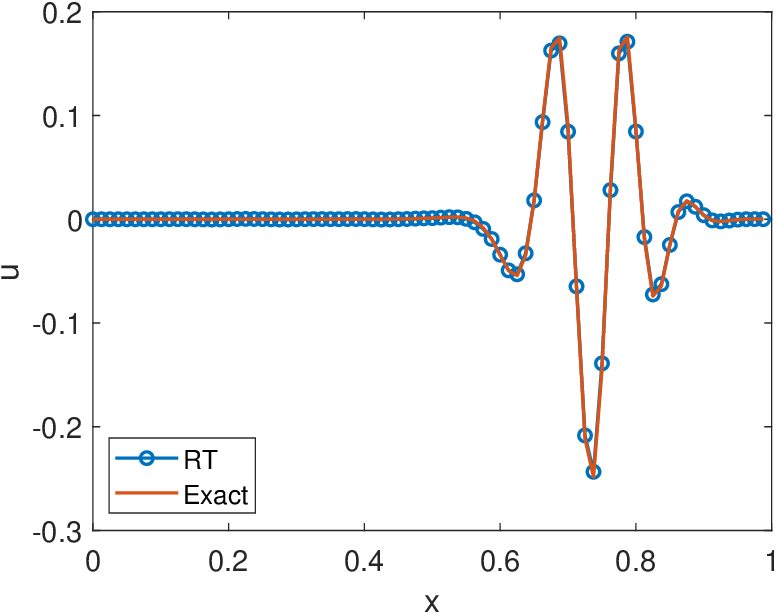}}
          \subfigure[]{
     \includegraphics[scale=0.35]{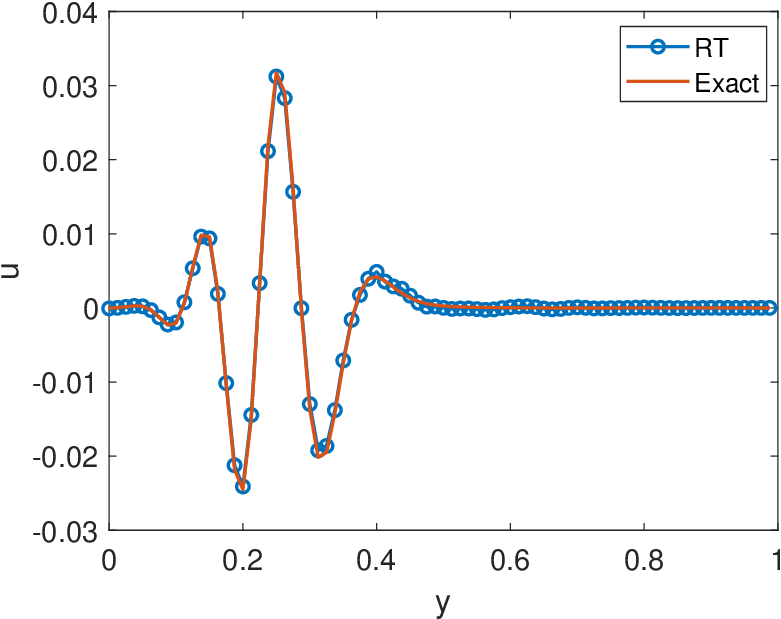}}
     \subfigure[]{
     \includegraphics[scale=0.35]{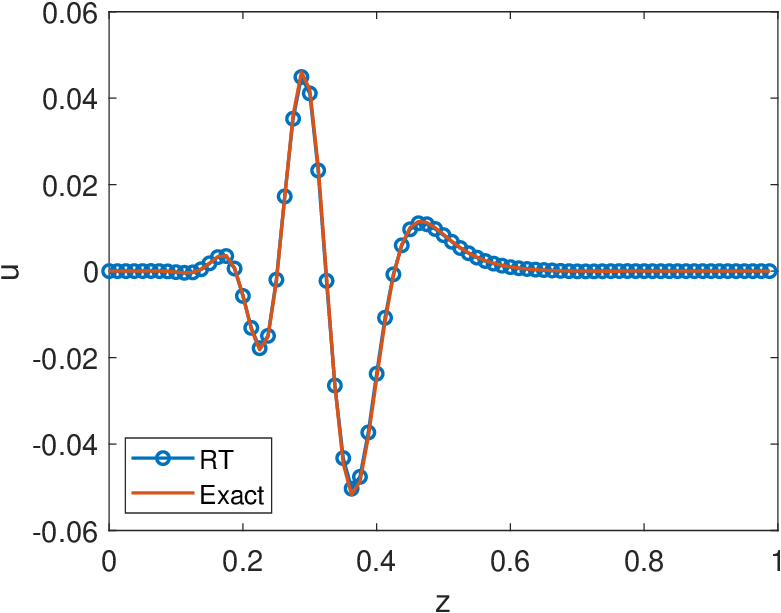}}\\
     \caption{Example 3. $\beta=16$, $T=0.4$. (a) a sectional slice at $z=0.75$;
 (b) a sectional slice at $x=0.45$; (c) a sectional slice at $y=0.15$; (d) comparison of the slices at $y=0.7$ and $z=0.75$; (e) comparison of the slices at $x=0.45$ and $z=0.2$; (f) comparison of the slices at $x=0.3$ and $y=0.15$. }\label{figure20}
     \end{figure}

          \begin{figure}[htbp]
     \centering
     \subfigure[]{
     \includegraphics[scale=0.35]{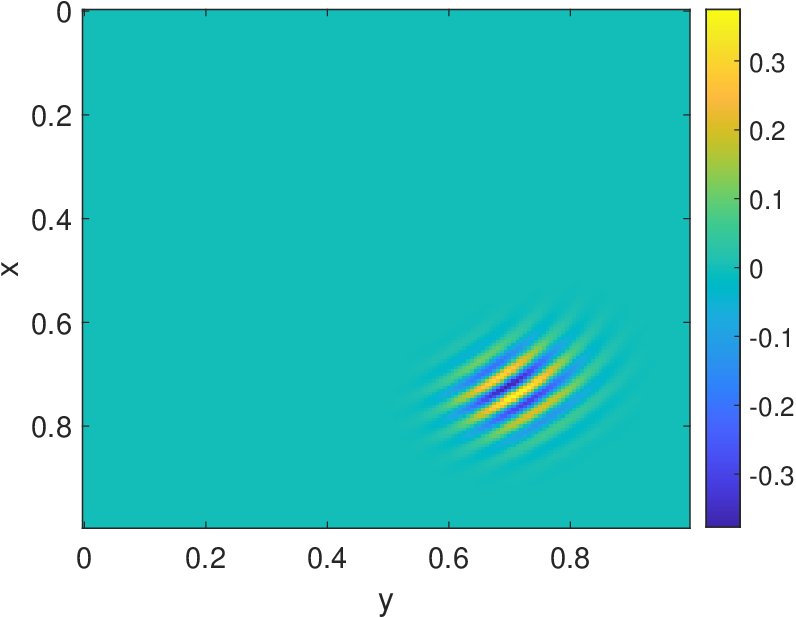}}
          \subfigure[]{
     \includegraphics[scale=0.35]{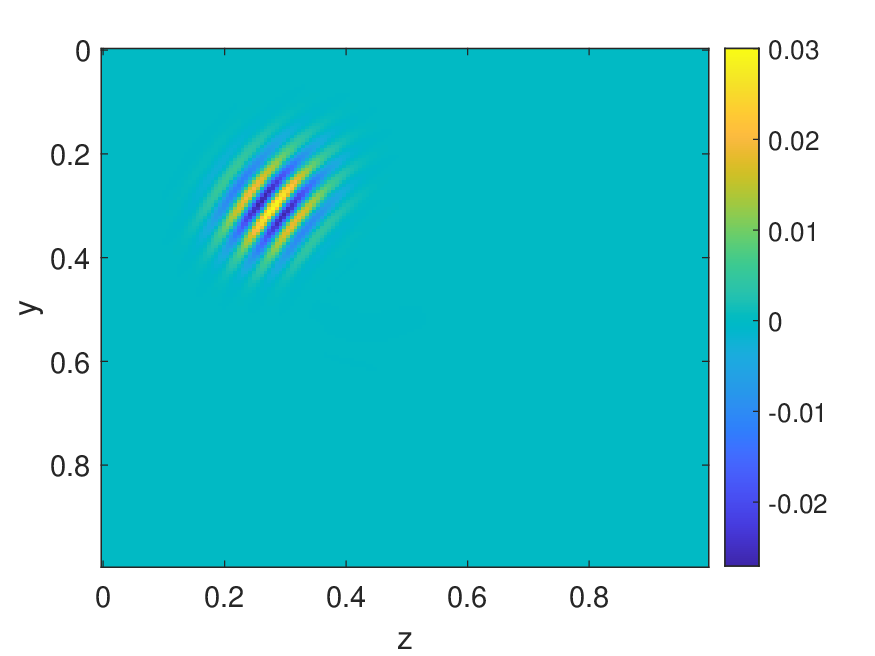}}
          \subfigure[]{
     \includegraphics[scale=0.35]{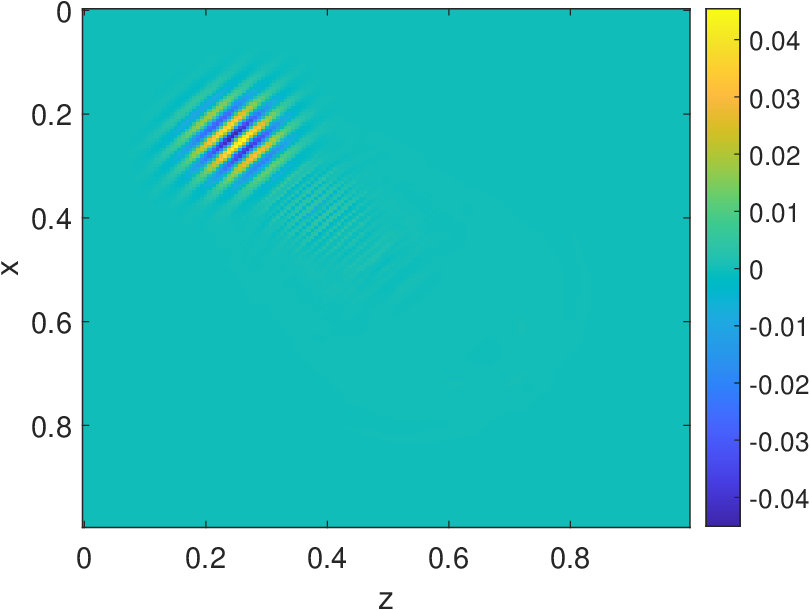}}\\
          \subfigure[]{
     \includegraphics[scale=0.35]{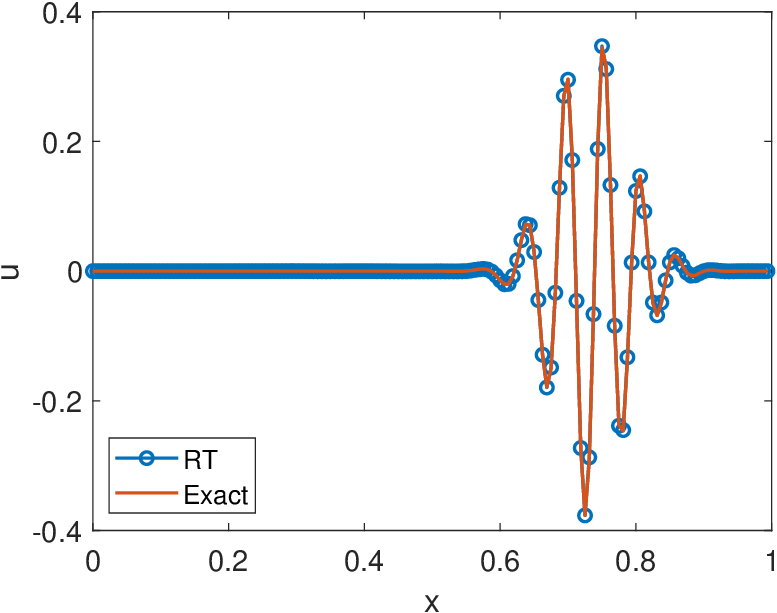}}
     \subfigure[]{
     \includegraphics[scale=0.35]{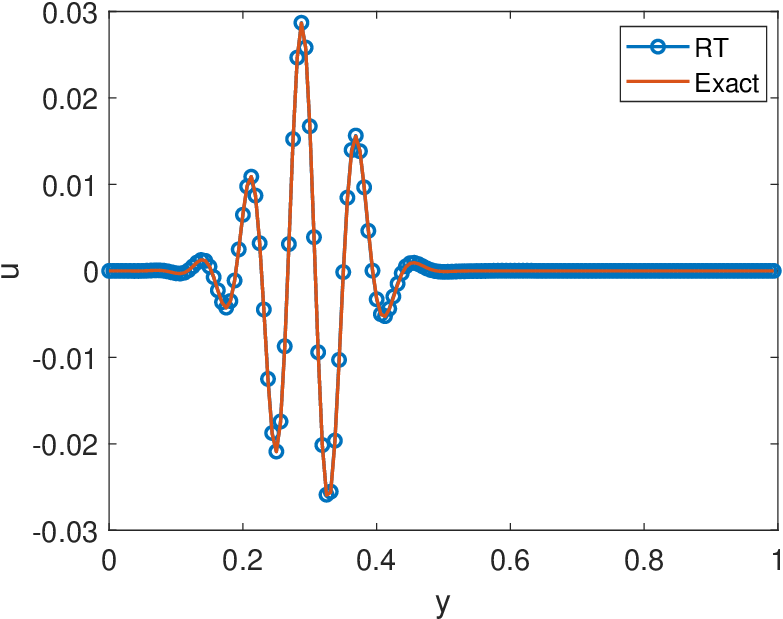}}
          \subfigure[]{
     \includegraphics[scale=0.35]{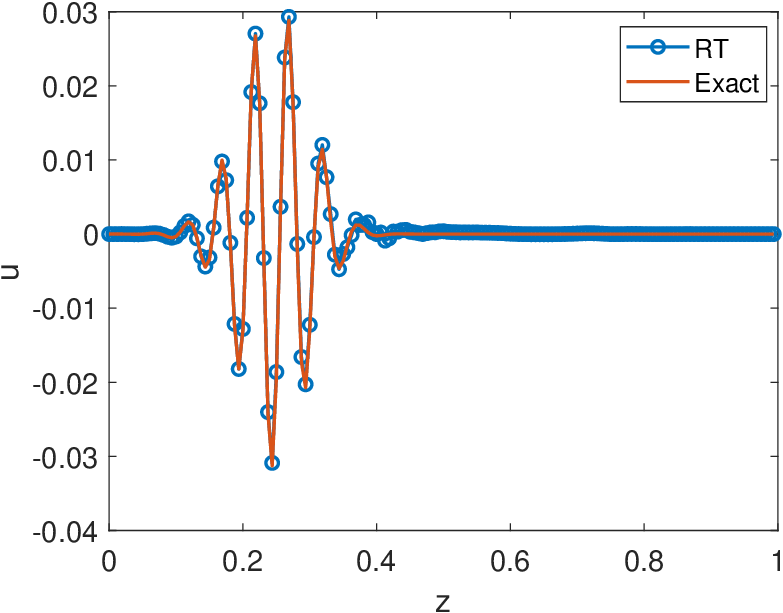}}
     \caption{Example 3. $\beta=32$, $T=0.4$. (a) a sectional slice at $z=0.75$; (b)  a sectional slice at $x=0.15$; (c) a sectional slice at $y=0.425$; (d) comparison of the slices at $y=0.7$ and $z=0.75$;  (e) comparison of the slices at $x=0.15$ and $z=0.3$; (f) comparison of the slices at $x=0.3$ and $y=0.425$.}\label{figure21}
     \end{figure}

\noindent\textbf{Example 4.} We use the following setup.
%This example is chosen to illustrate that our Hadamard integrator is able to handle caustics automatically and treat spatially overturning waves naturally.
\begin{itemize}
    \item $\rho=\left(\frac{1+\exp(-36(z-0.5))}{0.8+1.25\exp(-36(z-0.5))}\right)^2,$\; $\nu=1$,\; and \;$c=\frac{0.8+1.25\exp(-36(z-0.5))}{1+\exp(-36(z-0.5))}.$
    \item $u(0,x,y,z)=\sin(32\pi(x+y+z-1.5))\exp(-600((x-0.5)^2+(y-0.5)^2+(z-0.45)^2)),$\\
    $u_t(0,x,y,z)=0.$
    \item The computational domain is $\Omega=[0,1]^3$ and the grid size used to discretize $\Omega$ is $h=\frac{1}{160}.$
    \item $\Omega$ is divided into five parts to construct the low-rank representations,
    \begin{equation*}
        \tilde{\Omega}^{K}=\Omega\cap \{\bx :z\in [0.2(K-1),0.2K]\}.
    \end{equation*}
    \item The orders of the tensorized Chebyshev polynomials with respect to $(x_0,y_0,z_0,\theta_0,\xi_0)$ in different sub-regions are:
    \begin{itemize}
        \item $\tilde{\Omega}^1:[1,1,3,9,9]$;
        \item $\tilde{\Omega}^2:[1,1,15,15,15]$;
        \item $\tilde{\Omega}^3:[1,1,25,25,25]$;
        \item $\tilde{\Omega}^4:[1,1,15,15,15]$;
        \item $\tilde{\Omega}^5:[1,1,3,9,9]$.
    \end{itemize}
    \item The number of Gaussian-quadrature nodes is $M_3=96.$
\end{itemize}
Here we have chosen the $z$-dependent velocity $c$ to be analogous to a scaled, smoothed, and shifted Heaviside function in the $z$-direction. The slice of the velocity $c$ at $x=0.5$ and $y=0.5$ is shown in Figure \ref{figure20.1}(a), from which we can see that the velocity changes rapidly from $c=0.8$ to $c=1.25$  around $z=0.5$. We present the rays starting from $\bx_0=[0.5,0.5,0.45]$ in Figure \ref{figure20.1}(b). This velocity field produces overturning rays and a lot of caustics in $\Omega$. We will use this example to illustrate that the Hadamard integrator can not only handle caustics automatically but also treat spatially overturning waves naturally. We take $\Delta t=0.1$ to  construct the HKH propagator.

     \begin{figure}[htbp]
     \centering
     \subfigure[]{
     \includegraphics[scale=0.385]{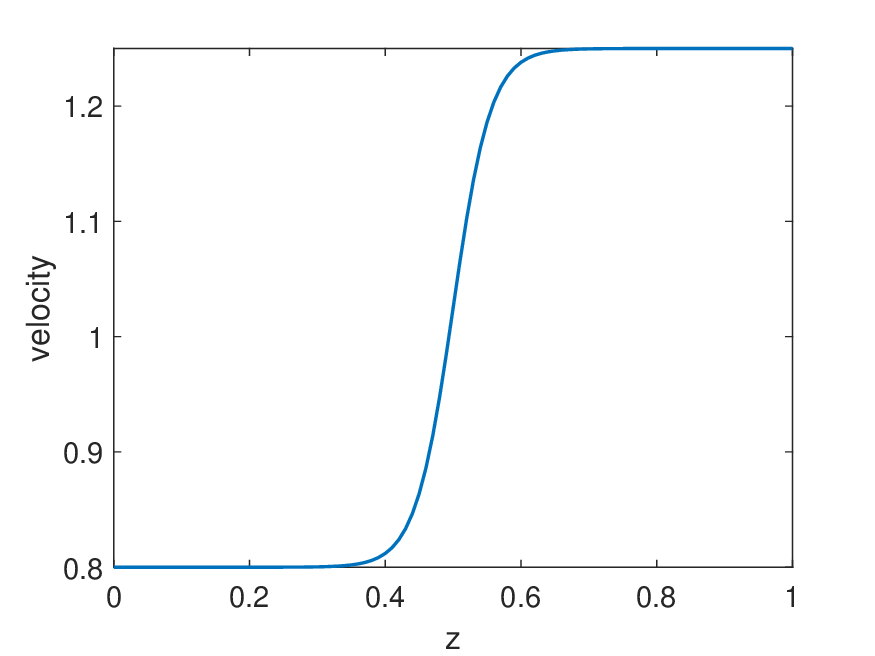}}
     \hspace{8mm}
     \subfigure[]{
     \includegraphics[scale=0.385]{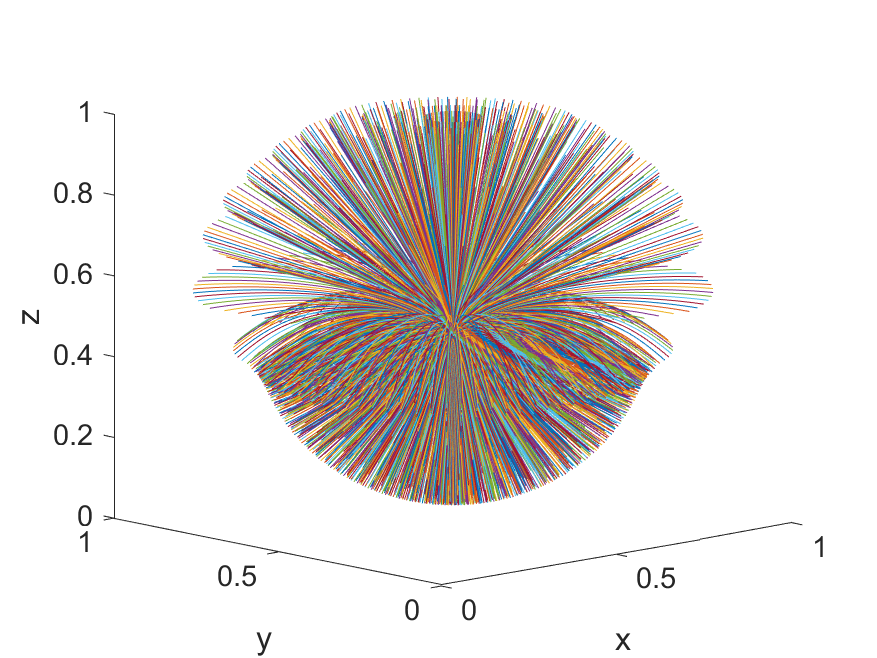}}
     \caption{Example 4. (a) The slice of velocity at $x=0.5$; (b) The rays with different take-off angles starting from the source $\bx_0=[0.5,0.5,0.45]$.}
     \label{figure20.1}
     \vspace{-3.75mm}
     \end{figure}
To construct the low-rank representations, we solve the Hamiltonian system (\ref{ode2}) using the Runge–Kutta method (RK4) for $1000$ time steps, and we solve the eikonal and transport equations with mesh size $h_1=0.005$ in a cubic neighborhood (\ref{6.16}) with $H=0.1$ for $\tilde{\Omega}^1$ and $\tilde{\Omega}^2$ and $H=0.15$ for $\tilde{\Omega}^3$, $\tilde{\Omega}^4$, and $\tilde{\Omega}^5$. We locate most of interpolation nodes in $\tilde{\Omega}^3$ to capture rapid changes in the velocity model.

In Figure \ref{example4.1}, we show the wavefields at $T=0.4$ in some caustic-free regions. Figure \ref{example4.1}(a) shows a $z$-section of the three-dimensional wavefield at $z=0.5$, and Figure \ref{example4.1}(d) shows slices at $y=0.8$ and $z=0.5$, in which we compare the exact solution (`-') with the RT solution (`o'). Figure \ref{example4.1}(b) shows an $x$-section of the three-dimensional wavefield at $x=0.28125$, and Figure \ref{example4.1}(e) shows slices at $x=0.28125$ and $z=0.3$, in which we compare the exact solution (`-')  with the RT solution (`o').  Figure \ref{example4.1}(c) shows a $y$-section of the three-dimensional wavefield at $y=0.25$, and Figure \ref{example4.1}(f) shows slices at $y=0.25$ and $z=0.35$, in which we compare the exact solution (`-') with the RT solution (`o'). As shown, The RT solutions match the exact solutions very well.

In Figure \ref{example4.2}, we show the wavefields at $T=0.4$ in some caustic-inducing regions. Figure \ref{example4.2}(a) shows an $x$-section of the three-dimensional wavefield at $x=0.625$, and Figure \ref{example4.2}(d) shows slices at $x=0.625$ and $y=0.8$, in which we compare the exact solution (`-') with the RT solution (`o'). Figure \ref{example4.2}(b) shows a $y$-section of the three-dimensional wavefield at $y=0.9375$, and Figure \ref{example4.2}(e) shows slices at $x=0.6$ and $y=0.9375$, in which we compare the exact solution (`-') with the RT solution (`o'). Figure \ref{example4.2}(c) shows a $y$-section of the three-dimensional wavefield at $y=0.75$, and Figure \ref{example4.2}(f) shows slices at $x=0.75$ and $y=0.75$, in which we compare the exact solution (`-') with the RT solution (`o').  The overturning waves refocus and bring about a significantly high amount of energy in wave motion in related regions in Figure \ref{example4.2}(a)(c). However, the Hadamard integrator still maintains the accuracy and the RT solutions match the exact solution very well. As demonstrated here, the Hadamard integrator can naturally propagate spatially overturning waves in time.
\begin{figure}[htbp]
     \centering
     \subfigure[]{
     \includegraphics[scale=0.35]{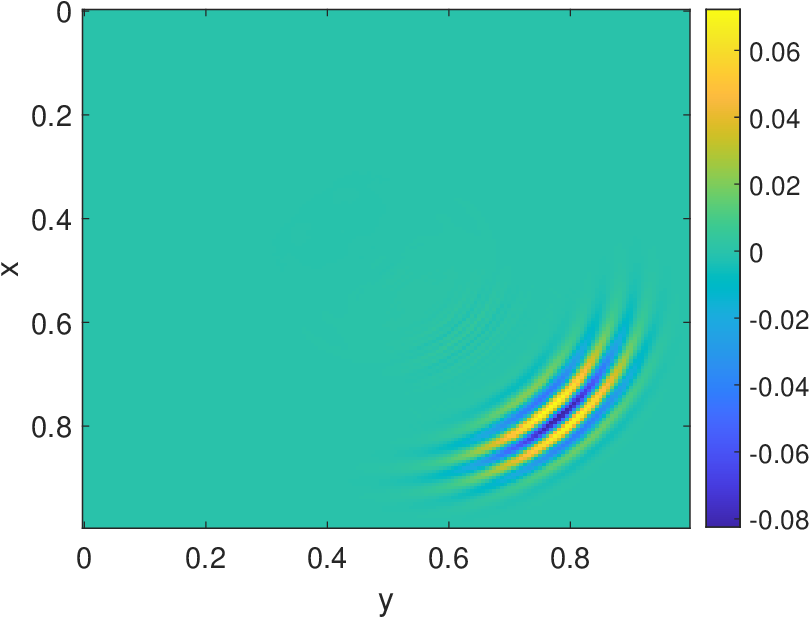}}
          \subfigure[]{
     \includegraphics[scale=0.35]{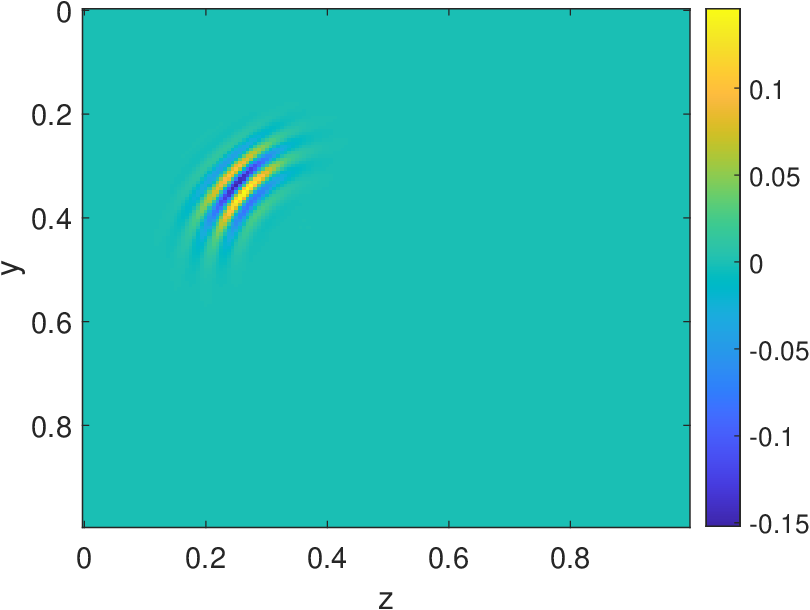}}
          \subfigure[]{
     \includegraphics[scale=0.35]{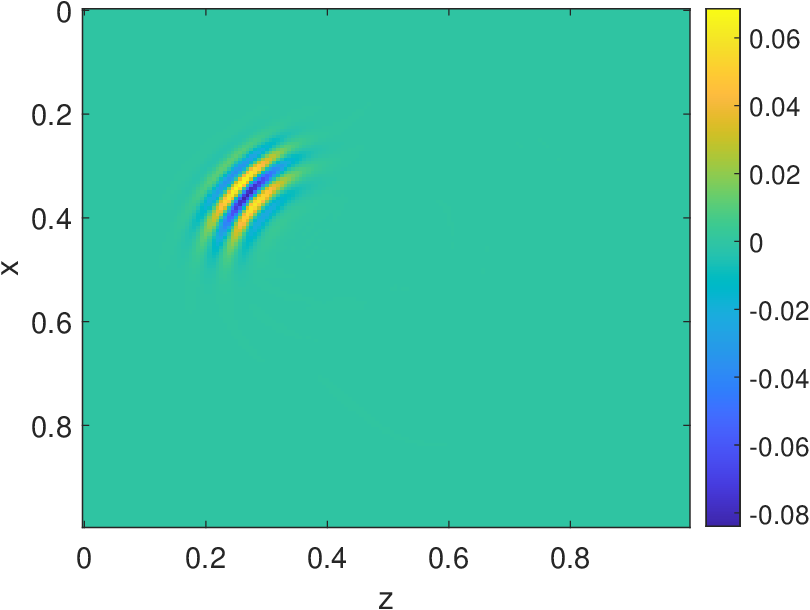}}\\
               \subfigure[]{
     \includegraphics[scale=0.35]{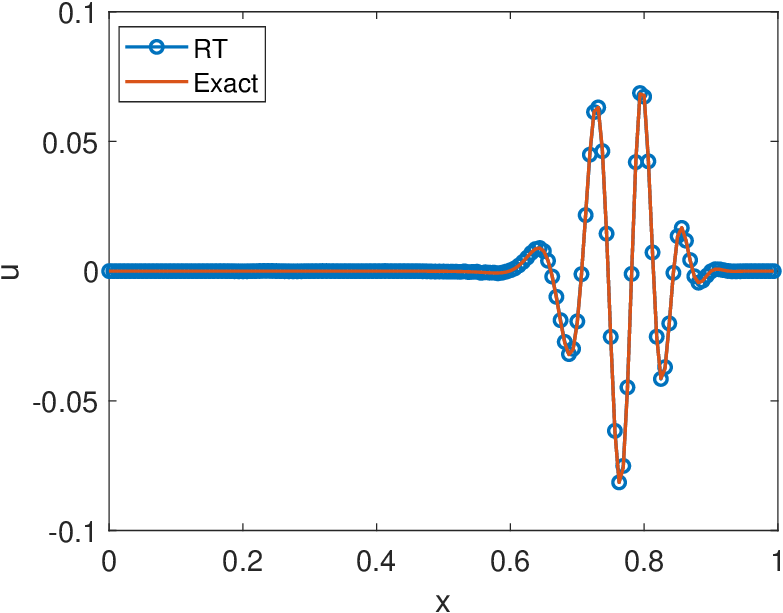}}
          \subfigure[]{
     \includegraphics[scale=0.35]{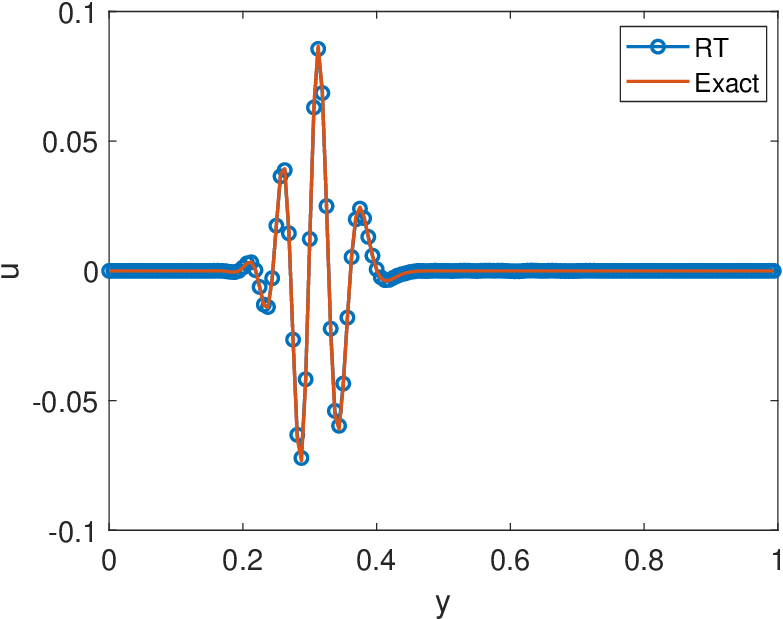}}
     \subfigure[]{
     \includegraphics[scale=0.35]{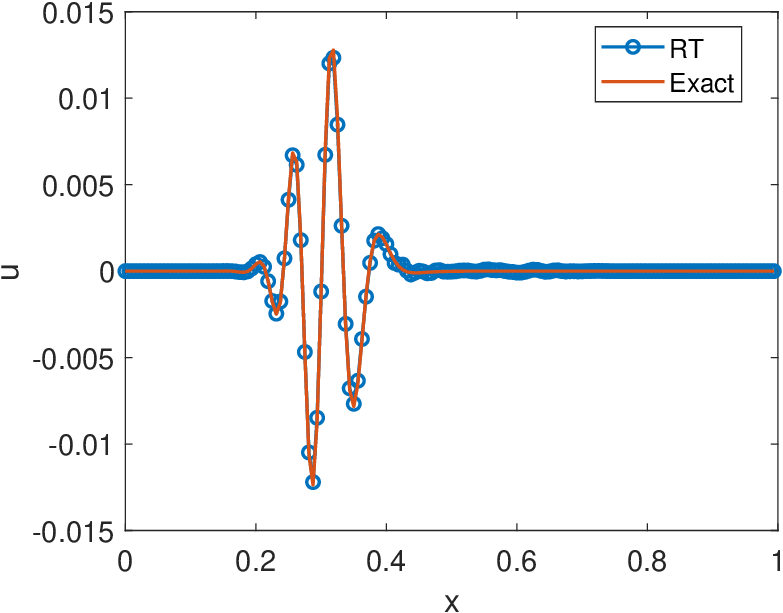}}\\
     \caption{Example 4. $T=0.4$. (a): a sectional slice at $z=0.5$; (b): a sectional slice at $x=0.28125$; (c): a sectional slice at $y=0.25$;(d):comparison of the slices at $y=0.8$ and $z=0.5$; (e):comparison of the slices at $x=0.28125$ and $z=0.3$;(f):comparison of the slices at $y=0.25$
     and $z=0.35$; }\label{example4.1}
     \end{figure}

\begin{figure}[htbp]
     \centering
     \subfigure[]{
     \includegraphics[scale=0.35]{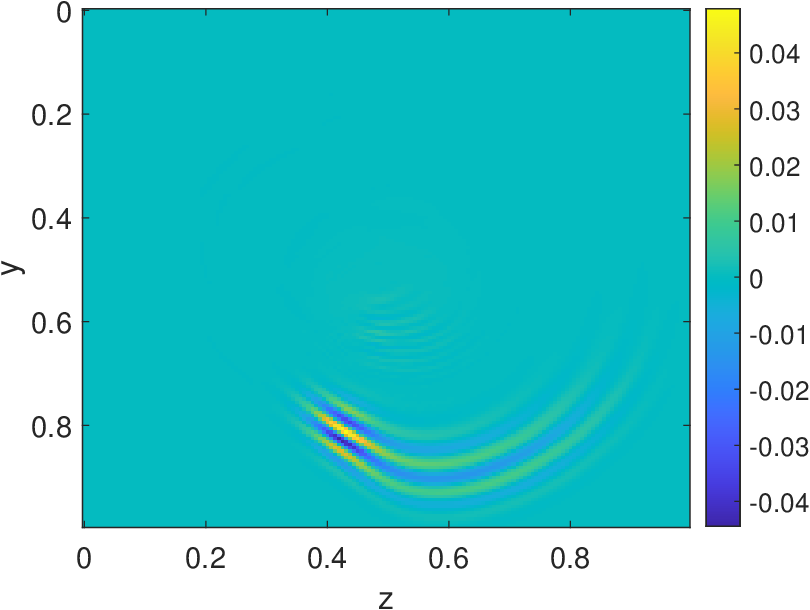}}
          \subfigure[]{
     \includegraphics[scale=0.35]{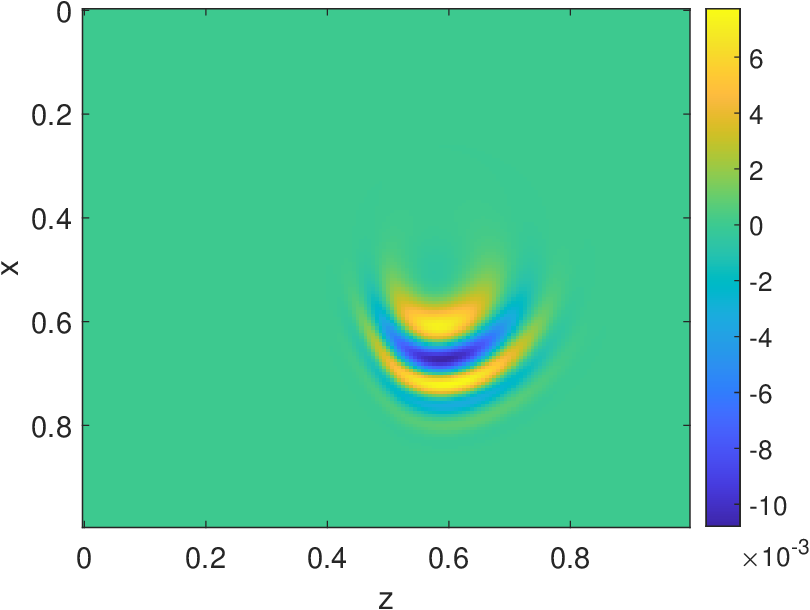}}
          \subfigure[]{
     \includegraphics[scale=0.35]{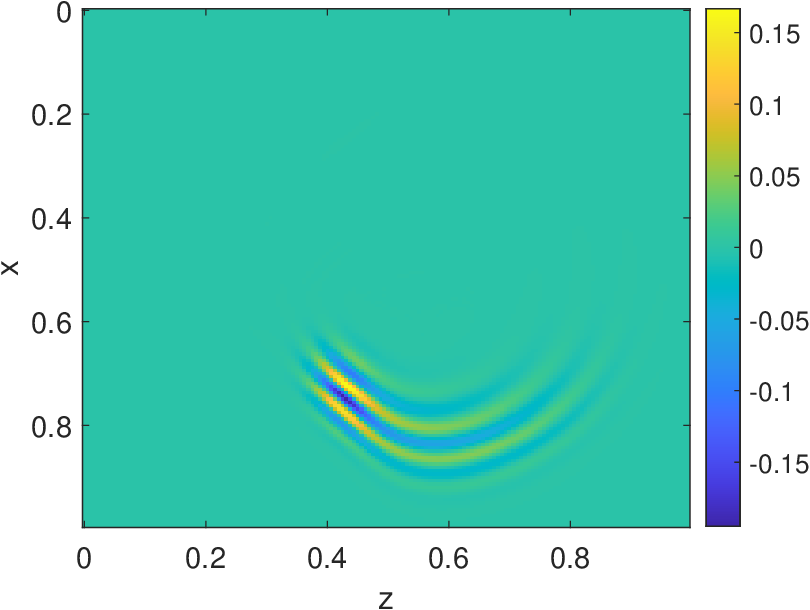}}\\
               \subfigure[]{
     \includegraphics[scale=0.35]{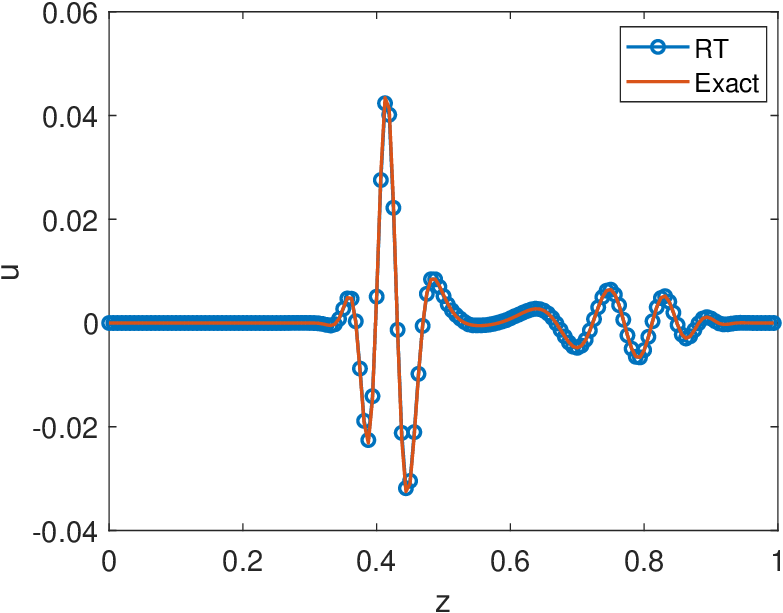}}
          \subfigure[]{
     \includegraphics[scale=0.35]{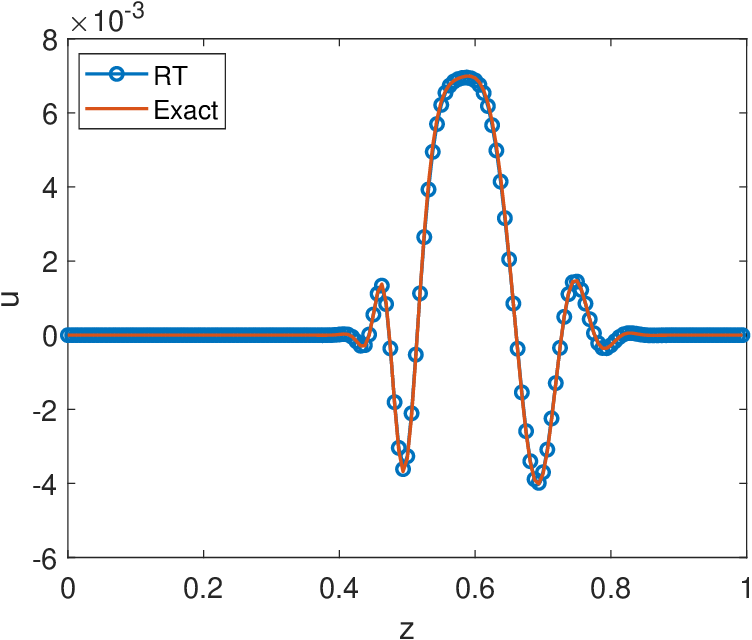}}
     \subfigure[]{
     \includegraphics[scale=0.35]{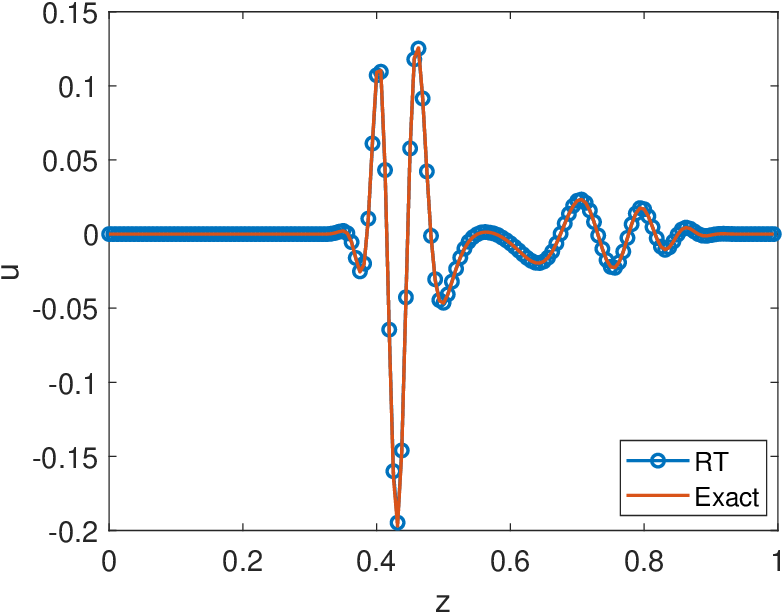}}\\
     \caption{Example 4. $T=0.4$. (a): a sectional slice at $x=0.625$; (b):a sectional slice at $y=0.9375$; (c): a sectional slice at $y=0.75$; (d): comparison of the slices at $x=0.625$ and $y=0.8$; (e): comparison of the slices at $x=0.6$ and $y=0.9375$; (f): comparison of the slices at $x=0.75$ and $y=0.75$;}\label{example4.2}
     \vspace{-5mm}
     \end{figure}
%%%%%%%%%%%%%%%%%%%%%%%%%%%%%%%%%%%%%%%%%%%%%%%%%%%%%%%%%%%%%%%%%%%%%%%%
 \section{Conclusion}
 Based on the Kirchhoff-Huygens representation and Hadamard's ansatz, we developed an original Hadamard integrator for solving time-dependent wave equations with highly oscillatory initial conditions. We derived the Lagrangian formulations via ray tracing and constructed low-rank representations for the wavefront locations and the Hadamard ingredients to accelerate the application of the integrator. By judiciously choosing a medium-dependent time step, the Hadamard integrator can propagate wavefields beyond caustics implicitly and advance spatially overturning wave naturally. Both two-dimensional and three-dimensional numerical examples illustrated the accuracy and performance of the new integrator. Applying this new integrator to seismic and medical imaging is an ongoing work.

\section*{Acknowledgement}
Cheng's research was supported by NSFC 11971121, 12241103 and the Sino-German Mobility Programme (M-0187) by the Sino-German Center for Research Promotion. Qian's research is partially supported by NSF 2012046, 2152011, and 2309534.

 \setcounter{equation}{0}
\renewcommand\theequation{A.\arabic{equation}} %（A1）
\begin{appendices}
\section{Jacobian of geodesic polar transformation}
    We now give the exact expression of the Jacobian of geodesic polar transformation. Consider
\beq{jacobian4}
\Int{\displaystyle{\cal T}}{} \nab^2\tau\,\dee\,V\,=\,\Int{\displaystyle{\cal T}}{} \nab\bcdot\nab\tau\,\dee\,V
\,=\,\Int{\displaystyle{\cal S}}{} \nab\tau\bcdot\,\bn\,\dee\,S\,.
\eeq
Here $\cal T$ is a volume bounded by a segment of a ray tube cut at one end by the surface $\tau = t_0$ and at the other by the surface $\tau=t_1$, and  $\bn$ is the outward unit normal to the surface $\cal S$, the boundary of $\cal T$.  $\cal S$ consists of the curvilinear tube of rays and the two ends consist of patches of the wavefronts $\tau=t_0$ and $\tau=t_1$, where $t_0<t_1$.  The tangents to the rays are parallel to $\nab\tau$ and so $\nab\tau\bcdot\,\bn = 0$ on the tube of rays, whereas $\bn$ is normal to the wave fronts and therefore parallel to $\nab \tau$ and so $\nab\tau\bcdot\,\bn= |\nab\tau|=1/c$ on $\tau=t_1$ and $\nab\tau\bcdot\,\bn= -|\nab\tau|=-1/c$ on $\tau=t_0$.  Putting this together we get
\beq{jacobian5}
\Int{\displaystyle{\cal T}}{} \nab^2\tau\,\dee\,V
\,=\,\Int{\displaystyle{\cal S}_1}{} |\nab\tau|\,\dee\,S
-\Int{\displaystyle{\cal S}_0}{} |\nab\tau|\,\dee\,S=\,\Int{\displaystyle{\cal S}_1}{} \fr{1}{c}\,\dee\,S
-\Int{\displaystyle{\cal S}_0}{}\fr{1}{c}\,\dee\,S
\,.
\,,
\eeq
where ${\cal S}_0$ and ${\cal S}_1$ are the patches cut out by the tube of rays on wavefronts $\tau=t_0$ and $\tau=t_1$.
%Recalling that $|\nab\tau|=1/c$ we see that \rf{jacobian5} reduces to
%\beq{jacobian6}
%\Int{\displaystyle{\cal T}}{} \nab^2\tau\,\dee\,V
%\,=\,\Int{\displaystyle{\cal S}_1}{} \fr{1}{c}\,\dee\,S
%-\Int{\displaystyle{\cal S}_0}{}\fr{1}{c}\,\dee\,S
%\,.
%\eeq
Utilizing the geodesic polar transformation, we have
%Let us now change variables of integration to $\bom$ and $\tau$, where $\bom$ is the initial direction of the ray.  First, we note that the volume element may be written as
%\beq{jacobian7}
%\dee\,s\,\dee\,S\,=\,c\,\dee\,\tau\,\dee\,S
%\eeq
%and $\dee\,S$ is the element of surface cut out on wavefronts $\tau\,=\,\mbox{const.}$ by the tube of rays.  Let
%  Then \rf{jacobian6} may be written as
\beq{jacobian8}
\Int{\displaystyle{\hat{\Omega}}}{}\,\dee\,\bom\Int{\displaystyle{t_0}}{t_1}
c\,\Big|\fr{\partial S}{\partial \bom}\Big|\nab^2\tau\,\dee\,\tau
\,=\,\Int{\displaystyle{\hat{\Omega}}}{}\,\dee\,\bom \,\fr{1}{c}\,\Big|\fr{\partial S}{\partial \bom}\Big|\Big]_{t_0}^{t_1}
\,.
\eeq
Here we have supposed that the tube of rays consists of the bundle of rays having take-off angles $\bom\in\hat{\Omega}$, where $\hat{\Omega}$ is a patch on the unit sphere. But $\hat{\Omega}$ is arbitrary and so we may equate the integrands with respect to $\bom$ to get
\beq{jacobian9}
\Int{\displaystyle{t_0}}{t_1}
c\,\Big|\fr{\partial S}{\partial \bom}\Big|\nab^2\tau\,\dee\,\tau
\,=\,\fr{1}{c}\,\Big|\fr{\partial S}{\partial \bom}\Big|\Big]_{t_0}^{t_1}
\,.
\eeq
On differentiating with respect to $t_1$ and dropping the subscript ${}_1$, we obtain the ordinary differential equation
\beq{jacobian10}
\fr{\dee}{\dee\,t}\,\Big(\fr{1}{c}\,\Big|\fr{\partial S}{\partial \bom}\Big|\Big)
\,=\,c\,\Big|\fr{\partial S}{\partial \bom}\Big|\nab^2\tau\,,
\eeq

Now let $J=\left|\frac{\partial S}{\partial \boldsymbol{\omega}}\right|$.
According to (\ref{jacobian10}), we have
    \begin{equation}\label{2.62}
      \nab^2 \tau=J^{-1} \frac{d}{d s}\left(c^{-1} J\right)=J^{-1} c \nab \tau \cdot \nab\left(c^{-1} J\right),
    \end{equation}
where
    %\begin{equation}\label{2.63}
      $\frac{d}{d s} \equiv c \nab \tau \cdot \nab.$
   % \end{equation}
 Letting $A_0^2=cJ^{-1}$, relation (\ref{2.62}) implies that $A_0^2$ satisfies the standard transport equation
\begin{equation}\label{2.64}
  \nab \cdot\left(A_0^2 \nab \tau\right)=0.
\end{equation}
Comparing (\ref{2.64}) with (\ref{SASWE2.17}), we get
\begin{equation}\label{2.65}
  \frac{c}{J}=A_0^2=\gamma \rho c^2 u_0^2\,,
\end{equation}
where $\gamma$ is a constant to be determined. Taking $\tau\rightarrow 0$ in (\ref{2.29}), we have
\begin{equation}\label{2.66}
  c_0^{m-1}=\frac{1}{\gamma \rho_0 c_0 v_0^2(\mathbf{0})}=\frac{4 \rho_0^2 \pi^{m-1}}{\gamma \rho_0 c_0 n_0^{2 m}},
\end{equation}
%where $c_0=c(\boldsymbol{0}),\rho_0=\rho(\boldsymbol{0})$.
Thus
\begin{equation}\label{2.67}
  \gamma=4 \rho_0 c_0^m \pi^{m-1},
\end{equation}
which means that
\begin{equation}\label{2.68}
  \left|\frac{\partial S}{\partial \omega}\right|=\frac{\tau^{m-1}}{4 \rho_0 c_0^m \pi^{m-1} \rho c v_0^2}.
\end{equation}
\end{appendices}
%%%%%%%%%%%%%%%%%%%%%%%%%%%%%%%%%%%%%%%%%%%%%%%%%%%%%%%%%%%%%%%%%
\bibliographystyle{plain}
\bibliography{myref}

%%%%%%%%%%%%%%%%%%%%%%%%%%%%%%%%%%%%%%%%%%%%%%%%%%%%%%%%%%%%%%%

\end{document}